\documentclass[12pt,leqno]{amsart}
\usepackage{mathtools,amsmath}
\usepackage{amssymb}
\usepackage{amsthm}
\usepackage{epsf}
\usepackage{graphicx}
\usepackage{esint} 
\usepackage{dsfont}
\usepackage[pagebackref]{hyperref} 
\hypersetup{pdfpagemode=FullScreen,  colorlinks=true} 
\usepackage{enumerate} 
\usepackage{xcolor,bbm,easybmat,bm}

\usepackage{amsfonts}
\usepackage{amsrefs}
\usepackage{verbatim}
\usepackage{booktabs}
\usepackage{color}
\usepackage{graphicx,float}

\usepackage{pgf,tikz}
\usetikzlibrary{cd,arrows,babel,shapes,positioning}
\usepackage{amsmath,amsfonts}
\usepackage{xcolor}
\usepackage{tikz} 
\usetikzlibrary{positioning, decorations.pathmorphing, arrows.meta, matrix}

\usepackage[letterpaper,hmargin=1in,vmargin=1.25in]{geometry}

\numberwithin{equation}{section}

\newcommand\footnoteref[1]{\protected@xdef\@thefnmark{\ref{#1}}\@footnotemark}
\makeatother

\newcommand{\EPS}{\varepsilon}
\newcommand{\Div}{\mathrm{div}}

\DeclareMathOperator{\supp}{supp}

\newcommand{\BMO}{{\rm BMO}}

\newcommand{\R}{\mathbb R}

\newcommand{\norm}[1]{\left\Vert#1\right\Vert}

\newcommand{\Loc}{\mathrm{loc}}

\newcommand{\AP}{A_{\Vert}}

\newcommand{\Ann}{\mathfrak{A}}

\newcommand{\sq}[1]{\widetilde{#1}}

\def\YYint#1#2#3{{\setbox0=\hbox{$#1{#2#3}{\iint}$}
		\vcenter{\hbox{$#2#3$}}\kern-.51\wd0}}
 
\theoremstyle{plain}
\newtheorem{thm}[equation]{Theorem}
\newtheorem{lemma}[equation]{Lemma}
\newtheorem{cor}[equation]{Corollary}
\newtheorem{prop}[equation]{Proposition}

\newtheorem{defin}[equation]{Definition}

\numberwithin{equation}{section}

\theoremstyle{plain}
\newcommand{\refthm}[1]{\emph{Theorem \ref{#1}}}
\newcommand{\reflemma}[1]{\emph{Lemma \ref{#1}}}
\newcommand{\refprop}[1]{\emph{Proposition \ref{#1}}}
\newcommand{\refcor}[1]{\emph{Corollary \ref{#1}}}

\theoremstyle{definition}
\newcommand{\refdef}[1]{\emph{Definition \ref{#1}}}

\theoremstyle{remark}
\newcommand{\refrem}[1]{\textit{Remark \ref{#1}}}

\theoremstyle{definition}

\theoremstyle{remark}
\newtheorem{rem}[equation]{Remark}

\allowdisplaybreaks[4]

\begin{document}

\title[Regularity problem for elliptic operators with BMO antisymmetric part]{The $L^p$ regularity problem for elliptic operators with transversally independent coefficients and BMO antisymmetric part}


\author[Pasquariello]{Marcus Pasquariello}
\address{Department of Mathematics, 
 Brown University, RI, US}
\email{marcus\_pasquariello@brown.edu}

\author[Pipher]{Jill Pipher}
\address{Department of Mathematics, 
 Brown University, RI, US}
\email{jill\_pipher@brown.edu}

\author[Ulmer]{Martin Ulmer}
\address{Department of Mathematics, 
 Brown University, RI, US}
\email{martin\_ulmer@brown.edu}


\begin{abstract} 
We study the solvability of both Poisson-Dirichlet and Regularity boundary value problems for a class of elliptic divergence form operators whose matrix splits into a bounded symmetric component and an unbounded antisymmetric component. Specifically, the unbounded antisymmetric component is assumed to have coefficients in the function space of bounded mean oscillation (BMO), which lies outside the scope of the classical De Giorgi-Nash-Moser theory. The Poisson-Dirichlet problem is shown to be solvable in a range dual to the solvability of the adjoint Dirichlet problem on chord arc domains. We apply this to obtain the first positive results on the Regularity problem (for $t$-independent matrices) in the upper half space, contributing to a growing body of work on this class of elliptic/parabolic operators.

\end{abstract}

\maketitle



\tableofcontents

\section{Introduction}

\noindent
In the present work, we consider linear elliptic operators of the type \(L:=\mathrm{div}(A\nabla \cdot)\). The matrix  \(A(x,t)\in\mathbb{R}^{(n+1)\times(n+1)}\) has measurable coefficients, is not assumed to be symmetric, and is uniformly elliptic, i.e. there exists \(\lambda_0>0\) with
    \[\lambda_0|\xi|^2\leq \xi^T A(x,t)\xi\qquad\textrm{ for every }\xi\in \mathbb{R}^{n+1}, \textrm{ and a.e. }(x,t)\in\Omega.\]
   
    We write $A^s$ and $A^a$ for the symmetric and antisymmetric part of $A$ respectively. We assume that the symmetric part is bounded, i.e. $\Vert A^s\Vert_{L^\infty}\leq\Lambda_0$, while the antisymmetric part is assumed merely to belong to $BMO(\Omega)$, i.e. for every component $(A^a)_{ij}$, $1\leq i,j\leq n+1$,
    \[\Vert (A^a)_{ij}\Vert_{BMO(\Omega)}:=\sup_{B\subset \Omega}\fint_B \big|(A^a)_{ij}(y)-\fint_B(A^a)_{ij}(z)dz\big|dy\leq \Lambda_0<\infty,\]
    where the supremum is taken over all balls $B\subset\Omega$.
    The domain $\Omega$ where such operators will be considered is assumed to be
    a {\it chord-arc domain} (see \refdef{def:CAD}), a stronger condition than being a {\it uniform} domain, which interacts nicely with the $BMO$ space: 
    Jones showed in \cite{jones_extension_1980} that a function $(A^a)_{ij}\in BMO(\Omega)$ can be extended to a function $(A^a)_{ij}\in BMO(\mathbb{R}^{n+1})$ if and only if $\Omega$ is a uniform domain (see \refdef{def:uniformDomain}). That is, Jones's result allows one to extend the antisymmetric matrix to all of $\mathbb{R}^{n+1}$, still satisfying the $BMO$ condition. 

    The study of elliptic operators with antisymmetric part in $BMO(\Omega)$ originated in \cite{seregin_divergence-free_2012}, who were interested in incompressible flows with low regularity drift term \(c\), i.e. the solutions to the PDE
    \[-\Delta u+c\cdot \nabla u=0.\]
    and its parabolic analog. Using the fact that the divergence-free 
    vector field \(c\) is the divergence of an
    antisymmetric tensor, they incorporated the drift term into the matrix as an
    antisymmetric addition.
   They further realized that the
   low regularity assumption of \(c\in BMO^{-1}\) implies that $L$ gives rise to a well-defined bilinear form, and much more: they proved the Harnack inequality and De Giorgi-Nash-Moser estimates such as interior regularity results of H\"older continuity for solutions. That is, one can rewrite the PDE in the form \(Lu=\mathrm{div}(A\nabla u)=0\), with an antisymmetric part \(A^a\in BMO(\Omega)\). In \cite{li_lp_2019}, the boundary elliptic measure associated to such an operators was studied and its relationship with the Green's function established. This is the context of the present work. This class of elliptic operators (and their parabolic counterparts) have drawn considerable attention in recent years; many results for elliptic/parabolic operators with bounded matrices \(A\) have been extended to the case with antisymmetric part in \(BMO\) (see \cites{escauriaza_hofmann_kato_2018,dong_kim_fundamental_2018,qian_xi_parabolic_2019,hofmann_lp_2022, hofmann_dirichlet_2022,Ataei2025, DNU2024, li_pipher_bmo}). 

    Taken together, these results suggest that much if not all of the current body of knowledge concerning boundary value problems for elliptic/parabolic operators could be extended to this low regularity setting,  previously outside the scope of the De Giorgi-Nash-Moser theory. The references cited above include the
    resolution of the Kato conjecture in both the elliptic and parabolic settings, the solvability of the $L^p$ Dirichlet problem with matrices independent of the transverse variable, and the development of the perturbation theory for such operators. 

    The first main result of this paper concerns the Poisson-Dirichlet problem.  
   To state this result, we set
    \[
    Y^{1,2}(\Omega)=\{u\in L^{2^*}(\Omega):\nabla u\in L^2(\Omega)\},
    \qquad Y^{1,2}_0=\overline{C_c^\infty(\Omega)}^{\norm{\nabla\cdot}_2},
    \]
    and point to the definitions of the square function \(\mathcal{S}\), the nontangential maximal function \(\mathcal{N}\), and the mean-valued Carleson function \(\mathcal{C}\) in Section~\ref{subsection:DefandNot}.

\begin{thm}[Square function bounds for the Poisson Dirichlet problem]\label{thm:MainThm1}
    Let $\Omega\subset\mathbb{R}^{n+1}$ be an CAD domain, and let $L:=\mathrm{div}(A\nabla\cdot)$ be an elliptic operator with bounded symmetric part $A^s$ and antisymmetric part $A^a\in BMO$. If $(\operatorname{D})_{p'}^{L^*}$ is solvable for some $p\in(1,\infty)$, then for every $H,F\in L^\infty_c(\Omega)$ the solution $v\in Y_0^{1,2}(\Omega)$ to 
    \begin{align}\begin{cases} L^*v=H-\mathrm{div}(F)&\textrm{in }\Omega\\v=0& \textrm{on }\partial\Omega\end{cases}\label{eq:PDEForvPD}\end{align}
    satisfies the square function bound
    \begin{align}\Vert \mathcal{S}(v)\Vert_{L^q(\partial\Omega)}&\leq C\big(\Vert \mathcal{N}(v)\Vert_{L^q(\partial\Omega)}+\Vert \sq{\mathcal{C}}_{2}(\delta H)\Vert_{L^{q}(\partial\Omega)}+ \Vert \sq{\mathcal{C}}_{2}(F)\Vert_{L^{q}(\partial\Omega)}\big),\label{eq:SquareFunctionEstimate}
    \end{align}

    \noindent for all $1<q<\infty$. Here the constant $C$ depends on $\lambda_0,\Lambda_0, n, q,$ and the chord arc domain constants. Moreover,
    \begin{align}
    \Vert \sq{\mathcal{N}}_{2^*}(v)\Vert_{L^{p'}(\partial\Omega)}&\lesssim\Vert \sq{\mathcal{C}}_{2_*}(\delta H)\Vert_{L^{p'}(\partial\Omega)} + \Vert \sq{\mathcal{C}}_{2}(F)\Vert_{L^{p'}(\partial\Omega)},\textrm{ and}\label{eq:PD2*inequalityElliptic}
    \\
    \Vert \mathcal{N}(v)\Vert_{L^{p'}(\partial\Omega)}&\lesssim\Vert \sq{\mathcal{C}}_{r'}(\delta H)\Vert_{L^{p'}(\partial\Omega)} + \Vert \sq{\mathcal{C}}_{q'}(F)\Vert_{L^{p'}(\partial\Omega)}\label{eq:PDinequalityElliptic}
    \end{align}
    hold for any \(q\in [1,\frac{n+1}{n}), r\in [1,\frac{n+1}{n-1})\).
\end{thm}
The article \cite{mourgoglou_solvability_2025} defines the \textit{Poisson Dirichlet problem} $(\operatorname{PD})_{p'}^{L^*}$ to be \textit{solvable} if for all solutions $v$ to \eqref{eq:PDEForvPD} and $H,F\in L^\infty_c(\Omega)$ the estimate \eqref{eq:PD2*inequalityElliptic} holds. Moreover in \cite{mourgoglou_solvability_2025}*{Theorem 1.11}, they show for operators with bounded antisymmetric part that solvability of the Dirichlet problem $(\operatorname{D})_{p'}^{L^*}$ implies both \eqref{eq:PD2*inequalityElliptic} and \eqref{eq:PDinequalityElliptic}. We note that \eqref{eq:PD2*inequalityElliptic} and \eqref{eq:PDinequalityElliptic} are estimates that do not follow from each other, since the left and right hand of \eqref{eq:PDinequalityElliptic} are larger than the left and right hand side of \eqref{eq:PD2*inequalityElliptic}, respectively. In fact, estimates of the type of \eqref{eq:PDinequalityElliptic} will be more useful for our application to the regularity problem.

The proof of \cite{mourgoglou_solvability_2025}*{Theorem 1.11} only uses tools that have been established for matrices with either bounded or BMO antisymmetric part. Hence, no further adaptations are needed for the more general class of operator with BMO antisymmetric part, and a verbatim proof allows us to conclude that if \(L:=\mathrm{div}(A\nabla\cdot)\) is an elliptic operator with antisymmetric part in BMO, then we also have the estimates \eqref{eq:PD2*inequalityElliptic} and \eqref{eq:PDinequalityElliptic} for all \(H,F\in L^\infty_c(\Omega)\) where \(v\in Y^{1,2}_0(\Omega)\) is the unique solution of \(Lv=H-\mathrm{div}(F)\). Thus, only \eqref{eq:SquareFunctionEstimate} requires proof in this new context. In Section~\ref{sec:ProofMainThm1}, we add some details to explain this assertion.

The square function bound \eqref{eq:SquareFunctionEstimate}, however, is an entirely new contribution for operators with antisymmetric part in \(BMO\). Interestingly, combining \eqref{eq:SquareFunctionEstimate} and \eqref{eq:PDinequalityElliptic} gives us a bound on the square function of $v$ that only depends on the Carleson function of the interior data, i.e.
    \[\Vert \mathcal{S}(v)\Vert_{L^{p'}(\partial\Omega)}\lesssim\Vert \sq{\mathcal{C}}_{r'}(\delta H)\Vert_{L^{p'}(\partial\Omega)} + \Vert \sq{\mathcal{C}}_{q'}(F)\Vert_{L^{p'}(\partial\Omega)}.\]

Square functions are a central tool in the theory of elliptic boundary value problems, and are used to study absolute continuity of the harmonic measure with respect to the surface measure, establish solvability of the Dirichlet or regularity problem, or in the study of perturbation theory.

In our second main result, we apply the 
square function bound \eqref{eq:SquareFunctionEstimate} from \refthm{thm:MainThm1} to prove solvability of the regularity problem when \(A\) has antisymmetric part in \(BMO\) and is independent in transversal or \(t\)-direction. Here we will
restrict the domain to the upper half space \(\Omega=\mathbb{R}^{n+1}_+=\mathbb{R}^n\times(0,\infty)\) and call the last component \(t\). The study of \(t\)-independent operators, i.e. operators with \(A(x,t)=A(x)\), began in \cite{jerison_kenig_1981}, where Jerison and Kenig observed that Dahlberg's theorem on the solvability of the $L^2$-Dirichlet problem for the Laplacian on a Lipschitz graph domain could be obtained by flattening of the boundary and solving the \(L^2\)-Dirichlet problem for an elliptic operator \(L:=\mathrm{div}(A\nabla \cdot)\) on the upper half space. Here, \(A\) would be determined by the Lipschitz function of the original domain, and is symmetric, bounded and \(t\)-independent. Under these conditions, (\cite{jerison_kenig_1981}) presented a ``Rellich" type identity that showed \(t\)-independence of \(A\) implies solvability for the \(L^2\)-Dirichlet problem. On the other hand, for $2 \times 2$ nonsymmetric matrices and on the upper half plane \(\Omega=\mathbb{R}\times(0,\infty)\), \cite{kenig_new_2000} showed that the harmonic measure lies in the Muckenhoupt class \(A_\infty(\sigma)\) of the surface measure \(\sigma\), and provided examples showing that this result was sharp. Thus we only have solvability in this general class of the \(L^p\)-Dirichlet problem for \textit{some} \(p<\infty\). The two-dimensional $A_\infty$ result for nonsymmetric matrices was extended much later to all dimensions (\cite{hofmann_square_2015}). The higher dimensional approach required new ideas using the tools developed in the resolution of the Kato conjecture. Notably, square function estimates played a prominent role in that proof. 

In the present work, we shift our focus to the regularity problem (see \refdef{def:L^pregularityProblem}). Solvability of the regularity problem for \(t\)-independent operators with nonsymmetric, bounded matrices was addressed in \cite{hofmann_regularity_2015} using layer potentials and the \(L^p\) Kato theory. Here we extend this result to operators with antisymmetric part in \(BMO\).

\begin{thm}[Solvability of regularity problem for $t$-independent operators]\label{thm:MainThm2}
    Let $\Omega:=\mathbb{R}_+^{n+1}$ be the upper half space, and let $L:=\mathrm{div}(A\nabla\cdot)$ be an elliptic operator with bounded symmetric part $A^s$ and antisymmetric part $A^a\in BMO$. If the matrix is $t$-independent, i.e. $A(x,t)=A(x)$, then the solvability range of the regularity problem is dual to the solvability range of the dual Dirichlet problem, i.e. $(\operatorname{R})_p^L$ is solvable if and only if $(\operatorname{D})_{p'}^{L^*}$ is solvable.
\end{thm}

Our proof follows a different route from that of \cite{hofmann_regularity_2015} and avoids the use of layer potentials. It can therefore be viewed as an alternative proof for \(t\)-independent operators with bounded matrices and real coefficients. We outline the main steps of the proof in Section~\ref{sec:ProofMainThm1}. As in \cite{hofmann_square_2015} and \cite{hofmann_regularity_2015}, the ellipticized heat semigroup \(\mathcal{P}_tf:=e^{-t^2L}f\) is a key object for us. (See also \cites{hofmann_dirichlet_2022, ulmer_solvability_2026, ulmer_equivalence_2025}.) In the resolution of the Dirichlet problem, several \(L^2\) area function bounds involving \(\mathcal{P}_tf\) were used, while the regularity problem requires \(L^p\) area function bounds involving \(\mathcal{P}_tf\). We prove these \(L^p\) area function bounds for operators with antisymmetric part in \(BMO\) in the appendix~\ref{sec:AreaFctBounds}. Approaches to solvability of the Dirichlet problem can take advantage of reducing the problem to proving
\(\omega\in A_\infty(\sigma)\) instead of dealing with the nontangential maximal function bound directly. This avenue is not available in dealing with the regularity problem, and one must prove the nontangential maximal function bound directly. Our strategy is to bound  \(\sq{\mathcal{N}}(\nabla (u-\mathcal{P}_tf))\) by duality, and use the \(L^p\) area function bounds for \(\mathcal{P}_tf\). This approach makes use of the solvability of the Poisson-Dirichlet problem, and specifically the square function bound \eqref{eq:SquareFunctionEstimate}.

The structure of this article is as follows. Section~\ref{section:DefandNot} establishes definitions and notation, Section~\ref{sec:GreenAndSolutionPrelims} collects important preliminaries about solutions and the Green's function. Section~\ref{sec:ProofSquare} establishes the square function bound and proves \refthm{thm:MainThm1}, and Section~\ref{sec:ProofRegProblem} provides the proof of the regularity problem \refthm{thm:MainThm2}. Lastly, the \(L^p\) area function estimates are found in appendix \ref{sec:AreaFctBounds}, and a second appendix section~\ref{appendix:technicalLemmas} is dedicated to providing rigorous justifications for some of the technical steps that, in a bounded coefficient setting, are fairly standard but require some adaptation for our setting.

\medskip
{\bf Acknowledgements.} AI was only used to audit final drafts, and it identified typos and minor gaps.

\section{Definitions and Notation}\label{section:DefandNot}

\subsection{Definitions and Notation}\label{subsection:DefandNot}
    As is customary, we denote a boundary ball centered at \(x\in\partial\Omega\) with radius \(r\) by \(\Delta(x,r):=B(x,r)\cap\partial\Omega\), and a corresponding Carleson region as \(T(\Delta(x,r)):=B(x,r)\cap\Omega\). The surface measure on \(\partial\Omega\) is the \(n\)-dimensional Hausdorff measure restricted on the boundary of the domain and denoted by \(\sigma:=\mathcal{H}^{n}|_{\partial\Omega}\).
  
    A nontangential cone with aperture \(\alpha>0\) is given by \[\Gamma_\alpha(x):=\{y\in\Omega; |x-y|< (1+\alpha)\delta(y)\},\]
    and \(\Gamma_\alpha^\tau(x):=\Gamma_\alpha(x)\cap\{y\in\Omega; \delta(y)\leq \tau\}\) denotes the at height \(\tau>0\) truncated cone. For a point in the boundary \(x\in\partial\Omega\), we define the mean-valued nontangential maximal function as \[\sq{\mathcal{N}}_{\hat{c},\alpha,q}(F)(x):=\sup_{y\in\Gamma_\alpha(x)}\Big(\fint_{B(y,\hat{c}\delta(y)/2)}|F(z)|^qdz\Big)^{1/q}. \]
    Here \(1\leq q<\infty\) and \(0<\hat{c}<2\). When  \(q=\infty\), we recover the classical nontangential maximal function \(\mathcal{N}_{\alpha}(F)(x):=\sup_{y\in\Gamma_\alpha(x)}|F(y)|\). From real variable arguments (\cite{fefferman_hpspaces_1972}) it is known that the \(L^p\) norms of the nontangential maximal function with different apertures \(\alpha\) or averaging constant \(\hat{c}\) are comparable, and we tend to drop the subscripts $\alpha$ and $\hat{c}$. For mean-valued nontangential maximal functions on 
    chord-arc domains (CAD), this result can be found in \cite{mourgoglou_solvability_2025}*{Lemma 2.1}. The definition of CAD, and other geometric properties of domains

    There is also a tent space duality between the averaged Carleson function and the mean-valued nontangential maximal function, where the averaged Carleson function is defined by 
    \begin{align}\sq{\mathcal{C}}_{\hat{c},q}(f)(x):=\sup_{r>0}\frac{1}{|\Delta(x,r)|}\int_{T(\Delta(x,r))}\big(\fint_{B(y,\hat{c}\delta(y)/2)}|f(z)|^qdz\big)^{1/q}dy\label{recall:Tildecarelsonfunctiondef},\end{align}
    
  for \(x\in \partial\Omega=\mathbb{R}^n\), \(1\leq q<\infty\), and \(0<\hat{c}<1\). The averaged Carleson functions have comparable \(L^p\) norms when the averaging constant \(\hat{c}\) varies (see \cite{mourgoglou_solvability_2025}*{Lemma 2.2}), and so we tend to drop the subscripts $\hat{c}$.

    Carleson's theorem shows that the nontangential maximal function and the Carleson function are quasi dual to each other - see \cite{carleson_interpolations_1962}*{Theorem 2}, and also some special cases in \cite{coifman_new_1985} and \cite{hytonen_carleson_2013}. On a CAD, Carleson's theorem can be found in \cite{mourgoglou_solvability_2025}{Prop. 2.4} and the definition of CAD along with other geometric properties of domains is below in subsection \ref{Geom}.

    \begin{prop}[Carleson's theorem]\label{prop:NontangnetialAndCarelsonAreDual}
        Let \(1<p,q<\infty\). Then for each \(u\in L^q_{loc}(\Omega), h\in L^{q'}_{loc}(\Omega)\)
        \[\left|\int_\Omega u h dx\right|\lesssim \Vert \sq{\mathcal{N}}_{q}(u)\Vert_{L^p(\partial\Omega)} \Vert \sq{\mathcal{C}}_{q'}(h)\Vert_{L^{p'}(\partial\Omega)},\]
        and
        \[\Vert \sq{\mathcal{N}}_q(u)\Vert_{L^p}\lesssim \sup_{\substack{H:H\in L^{q'}_{loc}(\Omega)\\ \Vert \sq{\mathcal{C}}_{q'}(H)\Vert_{L^{p'}(\partial\Omega)}=1}}\int_\Omega uH dx.\]
    \end{prop}

    As in \cite{coifman_new_1985}, we define the area function 
    \begin{align}\mathcal{A}_{\alpha}(F)(x):=\Big(\int_{\Gamma_\alpha(x)}\frac{|F(y)|^2}{\delta(y)^{n+1}}dy\Big)^{1/2},\label{recall:areafunctiondef}\end{align}
    and the square function as \(\mathcal{S}(F)(x):=\mathcal{A}(t\nabla F)(x)\). 
    We will also use the Carleson function from \cite{coifman_new_1985} which is defined for \(x\in\partial\Omega\) by
    \begin{align}\mathcal{C}(f)(x):=\Big(\sup_{r>0}\frac{1}{|\Delta(x,r)|}\int_{T(\Delta(x,r))}\frac{|f(y,t)|^2}{t}dydt\Big)^{1/2}.\label{recall:carelsonfunctiondef}\end{align}

    We now introduce a decomposition of the domain into Whitney cubes which we are going to use throughout the article.

    According to \cite{christ_tb_1990}*{Section 3}, we can decompose \(\Omega\) into a family of cubes \(\mathcal{D}:=\mathcal{D}(\Omega):=\{Q_\beta^k\}\) such that for certain fixed constants \(1/2<c_0<1<c_1<3/2\)
    \begin{enumerate}
        \item \(\forall \beta,k \exists x_{\beta,k}\in\Omega: B(x_{\beta,k},c_02^k)\subset Q_\beta^k\subset B(x_{\beta,k},c_12^k)\);
        \item \(\sum_{\beta,k}\chi_{Q_\beta^k}\equiv\chi_{\Omega} \);
        \item \(\frac{1}{8}l(Q_\beta^k)\leq\mathrm{dist}(Q_\beta^k, \partial\Omega)\leq 8 l(Q_\beta^k)\).
    \end{enumerate}
    The set \(\mathcal{D}_k:=\{Q_\beta^k;\beta\}\) denotes all Whitney cubes of generation \(k\) and we can enlarge every Whitney cube by the factor \(2\) and obtain a fixed constant \(N\), such that \(\chi_\Omega\leq\sum_{Q_\beta^k\in\mathcal{D}}\chi_{Q_\beta^k}\leq\sum_{Q_\beta^k\in\mathcal{D}}\chi_{2Q_\beta^k}\leq N\chi_\Omega \). We will often use \(W\) to denote a Whitney cube from \(\mathcal{D}\).

    Furthermore, fix $\alpha>1$, a truncation parameter $\tau\in (0,\infty]$, and boundary set $E\subset\partial\Omega$. We set the nested sawtooth regions over set $E$ as \(\sq{\Omega}^{(\tau,\alpha)}_i(E)\coloneqq\bigcup_{x\in E}\Gamma_{(i+1)\alpha}^{(i+1)\tau}(x)\) for \(i=0,1,2,3,4,5\), and define 
    \(\sq{\mathcal{D}}^{(\tau,\alpha)}_{i,j}(E):=\{W;W\in\mathcal{D}(\Omega),  W\cap \sq{\Omega}^{(\tau,\alpha)}_j(E)\setminus \sq{\Omega}^{(\tau,\alpha)}_i(E)\neq \emptyset \}\) 
    as the set of Whitney cubes that cover the area between two sawtooth regions 
    \(\sq{\Omega}^{(\tau,\alpha)}_j\) and \(\sq{\Omega}^{(\tau,\alpha)}_i\) where \(j>i\).

    Along with the Whitney cube decomposition, we are also going to make use of a smooth version thereof. For that, we define for $W\in \mathcal{D}(\Omega)$ a function $\phi_W$ with the properties that
    \begin{enumerate}\label{Section:Defphi_W}
    \item \(\phi_W\in C_0^\infty(\frac{3}{2}W)\)
    \item \(0\leq \phi_W\leq 1\) 
    \item \(\Vert\nabla \phi_W\Vert_{L^\infty}\approx\frac{1}{\mathrm{diam}(W)}\)
    \item \[\sum_{W\in \mathcal{D}(\Omega)}\phi_W\equiv \chi_\Omega;\]
     \item \begin{equation}\label{eq:SumOverPhiWOnSawtooth=1}\sum_{W\in \sq{\mathcal{D}}_{i-1,j+1}^{(\tau,\alpha)}(\Omega)}\phi_W\equiv 1\qquad \textrm{on } \sq{\Omega}_j^{(\tau,\alpha)}\setminus\sq{\Omega}_i^{(\tau,\alpha)};\textrm{ and}\end{equation}
     \item\[\sum_{W\in \sq{\mathcal{D}}_{i,j}^{(\tau,\alpha)}}\phi_W\equiv 0 
        \qquad\textrm{on }\Omega\setminus(\sq{\Omega}^{(\tau,\alpha)}_{j+1}\setminus\sq{\Omega}^{(\tau,\alpha)}_{i-1}).\]
    \end{enumerate}

    Since we are often going to work locally on Whitney cubes, it is convenient to introduce a discretized version of the area function as in \cite{kenig_neumann_1995}. For \(x\in \partial\Omega\), we set
    \[\mathcal{T}_{\alpha,q}(f)(x):=\sum_{W\in\mathcal{D};W\cap\Gamma_\alpha(x)\neq \emptyset}l(W)\Big(\fint_W |f(y)|^qdy\Big)^{1/q}.\]

    A useful connection between the averaged Carleson function and the discretized version of the area function is
    \begin{align}\Vert \sq{\mathcal{C}}_q(f)\Vert_{L^p(\partial\Omega)}\approx \Vert \mathcal{T}_q(f)\Vert_{L^p(\partial\Omega)}.\label{eq:CarlAreaTcomparable}\end{align}
    Although not explicitly stated, this follows from ideas of \cite{coifman_new_1985}. The proof is a straightforward combination of \reflemma{lemma:appendixAcompC}, \reflemma{lemma:appendixDiffApComp} and \reflemma{lemma:AppendixAcompT}.

Let us define the Dirichlet boundary value problem with boundary data in \(L^p\).

\begin{defin}[\((\operatorname{D})_p^L\)]\label{def:L^pDirichletProblem}
    We say the \textit{\(L^p\) Dirichlet boundary value problem} is solvable for \(L\) if for all boundary data for \(f\in C^\infty_c(\partial\Omega)\cap L^p(\partial\Omega)\) the solution \(u\) to the continuous Dirichlet problem
    \[\begin{cases}
        Lu=0&\Omega,
        \\
        u=f&\partial\Omega,
    \end{cases}\]
    satisfies
    \begin{align}\Vert \mathcal{N}(u)\Vert_{L^p(\partial\Omega)}\lesssim\Vert f\Vert_{L^p(\partial\Omega)},\label{eq:DPinequalityElliptic}\end{align}
    where the implied constants are independent of \(u\) and \(f\). In this case we also write that \((\operatorname{D})^L_p\) holds for \(L\).
\end{defin}

Furthermore, we call a Borel function \(g:\partial\Omega\to\mathbb{R}\) a \textit{Haj\l{}asz upper gradient} of \(f:\partial\Omega\to\mathbb{R}\) if
\[|f(X)-f(Y)|\leq |X-Y|(g(X)+g(Y))\qquad\textrm{for a.e. }X,Y\in \partial\Omega.\]
We denote the collection of all Haj\l{}asz upper gradients of \(f\) as \(\mathcal{D}(f)\) and define \(\dot{L}_1^p(\partial\Omega)\) by all \(f\) with
\[\Vert f\Vert_{\dot{L}^p_1(\partial\Omega)}:= \inf_{g\in\mathcal{D}(f)}\Vert g\Vert_{L^p(\partial\Omega)}<\infty.\]
This space is also called \textit{homogeneous Haj\l{}asz Sobolev space}. In the case of a flat boundary we have that \(\dot{L}_1^p(\mathbb{R}^n)=\dot{W}^{1,p}(\mathbb{R}^n)\) with comparable norms (see \cite{mourgoglou_regularity_2023}). 

We can define the Regularity problem with boundary data in \(\dot{L}_1^p(\partial\Omega)\) but in our case, on the Lipschitz graph domain, this could also be done equivalently with \(\dot{W}^{1,p}(\partial\Omega)\).

\begin{defin}[\((\operatorname{R})_p^L\)]\label{def:L^pregularityProblem}
    We say the \textit{\(L^p\) Regularity boundary value problem} is solvable for \(L\) if for all boundary data for \(f\in C^\infty_c(\partial\Omega)\cap \dot{L}^p_{1}(\partial\Omega)\) the solution \(u\) to the classical Dirichlet problem
    \[\begin{cases}
        Lu=0&\Omega,
        \\
        u=f&\partial\Omega,
    \end{cases}\]
    satisfies
    \begin{align}\Vert \sq{\mathcal{N}}(\nabla u)\Vert_{L^p(\partial\Omega)}\lesssim\Vert f\Vert_{\dot{L}_{1}^p(\partial\Omega)},\label{eq:RPinequalityElliptic}\end{align}
    where the implied constants are independent of \(u\) and \(f\). In this case we also write that \((R)^L_p\) holds for \(L\).
\end{defin}

\subsection{Geometric preliminaries}\label{Geom}

We work on a set $\Omega\subset\mathbb{R}^{n+1}$, which we call a \textit{domain} if $\Omega$ is open and connected. The Hausdorff dimension of the boundary of $\Omega$ is assumed to be $n$ and we denote the Hausdorff measure on the boundary by $\sigma:=\mathcal{H}^{n}|_{\partial\Omega}$. 

\begin{defin}[Ahlfors regular]
 If there exists a constant $C_0>0$ such that for every $x\in\partial\Omega$ and $0<r<\mathrm{diam}(\partial\Omega)$
\[C_0^{-1}r^n\leq \sigma(\partial\Omega\cap B(x,r))\leq C_0r^n,\]
then we say that the boundary of the domain is $n$\textit{-Ahlfors regular}.   
\end{defin}

\begin{defin}[Corkscrew condition]
    We say that a domain satisfies the \textit{interior corkscrew condition} if there exists constants $c_1>0$ such that for every $x\in\partial\Omega$ and $0<r<2\mathrm{diam}(\Omega)$ there exists a ball $B\subset B(x,r)\cap \Omega$ with $\mathrm{diam}(B)\geq c_1 r$. The center point of the ball $B$ is also called a \textit{corkscrew ball}, and its center point is called \textit{corkscrew point} of $\Delta=\Delta(x,r)$ and is denoted by $A(\Delta)$ or $A_\Delta$.
    
    Analogously, we say that a domain satisfies the \textit{exterior corkscrew condition} if there exists constants $c_1>0$ such that for every $x\in\partial\Omega$ and $0<r<2\mathrm{diam}(\Omega)$ there exists a ball $B\subset B(x,r)\cap \mathbb{R}^{n+1}\setminus\Omega$ with $\mathrm{diam}(B)\geq c_1 r$.

\end{defin}

\begin{defin}[Harnack chain condition]
	\( \Omega \) is said to satisfy the \emph{Harnack chain condition} if there is a constant \( c_0 \)
	such that for each \( \rho > 0, \; \Lambda \geq 1, \; x_1,x_2 \in \Omega \) 
	    with \( \delta(x_j) \geq \rho \) and \( |x_1 - x_2| \leq \Lambda \rho \),
	    there exists a chain of open balls \( B_1,\dots,B_N \subset \Omega \) 
	    with \( N \lesssim_{\Lambda} 1, \;  x_1 \in B_1, \; x_2 \in B_N, \; B_i \cap B_{i+1} \neq \emptyset \)
	    and \( c_0^{-1} r(B_i) \leq \mathrm{dist}(B_i,\partial \Omega) \leq c_0 r(B_i) \).
	The chain of balls is called a \emph{Harnack chain}.
\end{defin}

\begin{defin}[Uniform domain]\label{def:uniformDomain}
    The domain $\Omega$ is a uniform domain, if it satisfies the interior corkscrew condition and the Harnack chain condition.
\end{defin}

\begin{defin}[NTA]\label{def:NTA}
	The domain $\Omega$ is a NTA domain (nontangentially accessible domain), if it is a uniform domain and satisfies the exterior corkscrew condition.
\end{defin}

\begin{defin}[CAD]\label{def:CAD}
A domain \(\Omega\subset\mathbb{R}^{n+1}\) is called CAD domain (chord arc domain),
    if \(\Omega\) is a NTA set and $\partial\Omega$ is $n$-Ahlfors regular. 
\end{defin}

\subsection{Some useful results specific to matrices with BMO-antisymmetric parts}

We begin by recalling a few results specific to antisymmetric matrices that have been discovered and used in earlier works.

\begin{prop}\label{prop:ContantAntiSymmMatrix}
    Let \( b \) be a constant anti-symmetric matrix
        and let \( u \in W^{1,2}(E) \) and \( v \in W_0^{1,2}(E) \), with \( E \subset \Omega \) measurable.
    Then
    \[ \int_{E} b\nabla u \cdot \nabla v = 0. \]
\end{prop}

\begin{prop}\label{prop:DivCurl}
    Let $1<p<\infty$. Let $u\in \dot{W}^{1,p}(\mathbb{R}^n)$ and $v\in \dot{W}^{1,p'}(\mathbb{R}^n)$, then for any $1\leq i,j\leq n$, $\partial_i u\partial_jv-\partial_ju\partial_iv\in \mathcal{H}^1(\mathbb{R}^n)$ and
    \[\Vert \partial_i u\partial_jv-\partial_ju\partial_iv\Vert_{\mathcal{H}^1(\mathbb{R}^{n})}\lesssim\Vert \nabla u\Vert_{L^p(\mathbb{R}^{n})}\Vert \nabla v\Vert_{L^{p'}(\mathbb{R}^{n})}.\]
\end{prop}

\begin{prop}\label{prop:H1normofproduct}
    Let $1<p<\infty$. Let $u\in \dot{W}^{1,p}(\mathbb{R}^n)$ and $v\in \dot{W}^{1,p'}(\mathbb{R}^n)$, then for any $1\leq i\leq n$, $\partial_i (uv)\in \mathcal{H}^1(\mathbb{R}^n)$ and
    \[\Vert \partial_i (uv)\Vert_{\mathcal{H}^1(\mathbb{R}^{n})}\lesssim\Vert  u\Vert_{L^p(\mathbb{R}^{n})}\Vert \nabla v\Vert_{L^{p'}(\mathbb{R}^{n})}+\Vert \nabla u\Vert_{L^{p}(\mathbb{R}^{n})}\Vert  v\Vert_{L^{p'}(\mathbb{R}^{n})}.\]
\end{prop}

The following higher integrability result for the gradient of solutions to the Poisson-Dirichlet problem has not appeared yet in the literature for matrices with antisymmetric part in BMO.

\begin{lemma}\label{lemma:HigherIntegrabilityForGradient}
    Let \(u\in W^{1,2}_{loc}(\Omega)\) be a weak solution to \(Lu=H+\mathrm{div}(F)\). Then there exists \(\EPS>0\) such that if $B(x,r)$ is a ball with $B(x,4r)\subset\Omega$, then
    \begin{align}\Big(\fint_{B(x,r)}|\nabla u|^{2+\EPS}dy\Big)^{\frac{1}{2+\EPS}} & \lesssim \Big(\fint_{B(x,2r)}|\nabla u|^{2}dy\Big)^{\frac{1}{2}} \\
        & + \Big(\fint_{B(x,2r)}|\delta H|^{2+\EPS}dy\Big)^{\frac{1}{2+\EPS}} + \Big(\fint_{B(x,2r)}|F|^{2+\EPS}dy\Big)^{\frac{1}{2+\EPS}},\label{eq:RevBoundaryHölder1}\end{align}
    and for every $s\in (1,\frac{n+1}{n})$ and $c\in\mathbb{R}$
    \begin{align}
    \int_{B(x,r)}|\nabla u|^2 dy&\lesssim \int_{B(x,2r)}|\delta H|^2 dy+\int_{B(x,2r)}|F|^2 dy +  \frac{C_\sigma}{r^{2-2(n+1)/s'}}\Big(\int_{B(x,2r)}|u-c|^{\frac{2s}{2-s}} dy\Big)^{\frac{2-s}{s}}.\label{eq:RevBoundaryHölder2}
    \end{align}
    Moreover, if \(x\in \partial\Omega\) and \(r>0\), then \begin{equation}\label{eq:RevBoundaryHölder3}\begin{split}\Big(\fint_{T(\Delta(x,r))}|\nabla u|^{2+\EPS}dy\Big)^{\frac{1}{2+\EPS}}\lesssim \Big(\fint_{T(\Delta(x,2r))}|\nabla u|^{2}dy\Big)^{\frac{1}{2}} + \Big(\fint_{T(\Delta(x,2r))}|\delta H|^{2+\EPS}dy\Big)^{\frac{1}{2+\EPS}}
    \\
    + \Big(\fint_{T(\Delta(x,2r))}|F|^{2+\EPS}dy\Big)^{\frac{1}{2+\EPS}}.\end{split}
    \end{equation}
\end{lemma}
The proof of this lemma is very similar to the proof for elliptic operators with bounded antisymmetric part as in the original work of \cite{meyer_gradient_1963} and follows essentially from Caccioppoli inequality, Sobolev inequality and Gehring's lemma (see for instance \cite[Prop. 1.1 in Chapter 5]{giaquinta_multiple_1983}). Given that some adaptations are required to deal with the unboundedness of the antisymmetric matrix, a proof is provided in  Appendix~\ref{sec:AppendixHigherIntProof}.

\section{Properties of the Green's function and solutions}\label{sec:GreenAndSolutionPrelims}

In this section, we define the Green's function in our setting and record some of its properties. We will closely follow the approaches of \cite{hofmann_greens_function_2007}, \cite{li_pipher_bmo}, and \cite{mourgoglou_greens_function_2023}, but adapt it to the present setting of operators with BMO antisymmetric part on an unbounded domain. Throughout this section $\Omega$ will be an NTA domain, and
\[
    2^*=\frac{2(n+1)}{n-1},\qquad 2_*=\frac{2(n+1)}{n+3}.
\]
We write
\[
    \mathcal{B}[u,v]:=\int_\Omega A\nabla u\cdot \nabla v\,dx,
    \qquad u,v\in Y^{1,2}_0(\Omega),
\]
for the bilinear form associated with \(L=-\Div(A\nabla \cdot)\).
The form is understood in the usual sense for matrices whose symmetric part is bounded and whose antisymmetric part has bounded \(\BMO\) seminorm; the antisymmetric part is handled by the standard div-curl/Hardy--\(\BMO\) estimate (see \cite{li_lp_2019}).
Thus \(\mathcal{B}\) is bounded on \(Y^{1,2}_0(\Omega)\times Y^{1,2}_0(\Omega)\), and ellipticity gives
\[
    \mathcal{B}[u,u]\geq \lambda_0 \|\nabla u\|_{L^2(\Omega)}^2.
\]
Here and below the implicit constants depend only on \(n\), \(\lambda_0\), and \(\Lambda_0\).

\begin{prop}\label{prop:green-variational-poisson}
Let \(H\in L^{2_*}(\Omega)\) and \(F\in L^2(\Omega)\).
There exists a unique weak solution \(u\in Y^{1,2}_0(\Omega)\) to
\begin{equation}\label{eq:green-variational-identity}
    \begin{cases}
        Lu = H-\Div F & \text{in }\Omega,\\
        u=0 & \text{on }\partial\Omega
    \end{cases}
\end{equation}
Moreover,
\begin{equation}\label{eq:green-variational-estimate}
    \|u\|_{Y^{1,2}(\Omega)}
    \lesssim
    \|H\|_{L^{2_*}(\Omega)}+\|F\|_{L^2(\Omega)}.
\end{equation}
\end{prop}

\begin{proof}
The existence of a weak solution follows from the boundedness and coercivity of the bilinear form $\mathcal{B}$ and the Lax-Milgram theorem. Moreover, we have that 
\begin{align*}
    \|u\|_{Y^{1,2}(\Omega)}^2 & \lesssim \|\nabla u\|_{L^2(\Omega)}^2 \\
    & \lesssim \mathcal{B}[u, u] \\
    & = \int_{\Omega} H u \, dx + \int_{\Omega} F \nabla u \, dx \\
    & \leq \|H\|_{L^{2_*}(\Omega)}\|u\|_{L^{2^*}(\Omega)} + \|F\|_{L^2(\Omega)}\|\nabla u\|_{L^2(\Omega)} \\
    & \leq \|u\|_{Y^{1,2}(\Omega)}(\|H\|_{L^{2_*}(\Omega)}+\|F\|_{L^2(\Omega)})
\end{align*}
and the result follows.
\end{proof}

\begin{lemma}[Local estimates for BMO-antisymmetric coefficients]\label{Lemma:LocalEstimatesForBMO}
Let \(A=A^s+A^a\), where \(A^s\) is real, bounded, symmetric, and uniformly elliptic, and where \(A^a\) is real, antisymmetric, and has finite \(\BMO(\Omega)\) seminorm.
Then \(L=-\Div(A\nabla\cdot)\) and \(L^*=-\Div(A^T\nabla\cdot)\) satisfy the following standard local estimates.
If \(u\) is a weak solution of \(Lu=0\) in \(B_{2r}\subset\Omega\), then
\[
    \fint_{B_r}|\nabla u|^2\lesssim r^{-2}\fint_{B_{2r}}|u|^2
    \qquad\text{and}\qquad
    \sup_{B_r}|u|\lesssim \left(\fint_{B_{2r}}|u|^2\right)^{1/2}.
\]
If \(u\geq 0\) is a weak solution of \(Lu=0\) in \(B_{2r}\subset\Omega\), then
\[
    \sup_{B_r}u\lesssim \inf_{B_r}u,
    \qquad
    \left(\fint_{B_r}u^p\right)^{1/p}\lesssim_p \inf_{B_r}u
    \quad\text{for every }1\leq p<\infty.
\]
Moreover, if \(D\Subset\Omega\) is a bounded subdomain and \(v\in W^{1,2}_{\Loc}(D)\cap C(\overline D)\) is a weak subsolution of \(Lv=0\), then
\[
    \sup_D v\leq \sup_{\partial D}v^+.
\]
The same estimates hold for \(L^*\). These estimates are proved in \cite{li_lp_2019} for bounded domains with BMO antisymmetric part.
The arguments are local, except for the maximum principle, which is used here only on bounded subdomains.
Thus, the same estimates carry over to unbounded domains \(\Omega\).
\end{lemma}

\begin{thm}[Green's function]\label{thm:green-function-construction}
For every \(y\in \Omega\) there exists a function \(G(\cdot,y)\), called the Green's function for \(L\) with pole at \(y\), with the following properties.
\begin{enumerate}
    \item \(G(\cdot,y)\geq 0\) a.e. in \(\Omega\), \(G(\cdot,y)\in W^{1,1}_{\Loc}(\Omega)\), and
    \[
        G(\cdot,y)\eta\in Y^{1,2}_0(\Omega\setminus B(y,r))
        \qquad\text{for every }r>0,\textrm{ and }\eta\in C^\infty(\Omega),\eta\equiv 0\textrm{ on }B(y,r).
    \]
    \item If \(H\in L^\infty_c(\Omega)\) and \(F\in L^\infty_c(\Omega;\R^{n+1})\), then the solution \(u\in Y^{1,2}_0(\Omega)\) of \(L^*u=H-\Div F\) is given by
    \begin{equation}\label{eq:green-representation-adjoint}
        u(y)=\int_\Omega G(x,y)H(x)\,dx
        +\int_\Omega \nabla_xG(x,y) F(x)\,dx
    \end{equation}
    for a.e. \(y\in\Omega\).
    \item If \(\widehat G\) is another function satisfying the preceding representation property, then
    \[G(\cdot,y)=\widehat G(\cdot,y)\]
    for a.e. point of \(\Omega\setminus\{y\}\).
\end{enumerate}
In addition, for every \(r>0\),
\begin{equation}\label{eq:green-y12-away}
    \|G(\cdot,y)\|_{Y^{1,2}(\Omega\setminus B(y,r))}
    \lesssim r^{\frac{1-n}{2}}.
\end{equation}
If \(0<r<\delta(y):=\operatorname{dist}(y,\partial\Omega)\), then
\begin{align}
    \|G(\cdot,y)\|_{L^p(B(y,r))}
    &\lesssim_p r^{1-n+\frac{n+1}{p}},
    &&1\leq p<\frac{n+1}{n-1},\label{eq:green-local-Lp}\\
    \|\nabla G(\cdot,y)\|_{L^p(B(y,r))}
    &\lesssim_p r^{-n+\frac{n+1}{p}},
    &&1\leq p<\frac{n+1}{n}.\label{eq:green-grad-local-Lp}
\end{align}
Moreover, for every \(t>0\),
\begin{align}
    |\{x\in\Omega:G(x,y)>t\}|
    &\lesssim t^{-\frac{n+1}{n-1}},\label{eq:green-weak-bound}\\
    |\{x\in\Omega:|\nabla_xG(x,y)|>t\}|
    &\lesssim t^{-\frac{n+1}{n}}.\label{eq:green-grad-weak-bound}
\end{align}
Finally,
\begin{equation}\label{eq:green-pointwise-upper}
    G(x,y)\lesssim |x-y|^{1-n}
    \qquad\text{for a.e. }x,y\in\Omega,\ x\neq y.
\end{equation}
\end{thm}

\begin{proof}
The proof of this theorem follows exactly as in the proof of Theorem 6.1 in \cite{mourgoglou_greens_function_2023}.
\end{proof}

\begin{cor}[Adjoint Green's function]\label{cor:adjoint-green-function}
The adjoint operator \(L^*=-\Div(A^T\nabla\cdot)\) has a Green's function \(G^*\) satisfying the same conclusions as \(G\).
Moreover,
\begin{equation}\label{eq:green-adjoint-symmetry}
    G^*(x,y)=G(y,x)
    \qquad\text{for a.e. }(x,y)\in\Omega\times\Omega,\ x\neq y.
\end{equation}
Consequently, if \(u\in Y^{1,2}_0(\Omega)\) solves \(Lu=H-\Div F\) with \(H,F\in L^\infty_c(\Omega)\), then
\begin{equation}\label{eq:green-representation-L}
    u(x)=\int_\Omega G^*(y,x)H(y)\,dy
    +\int_\Omega \nabla_yG^*(y,x)\cdot F(y)\,dy
\end{equation}
for a.e. \(x\in\Omega\).
\end{cor}

\begin{proof}
The proof of this theorem follows exactly as in the proof of Theorem 6.1 in \cite{mourgoglou_greens_function_2023}.
\end{proof}

\begin{defin}[Elliptic measure]\label{def:elliptic-measure}
For \(x\in\Omega\), the elliptic measure \(\omega_L^x\) is the Borel probability measure on \(\partial\Omega\) such that, for every \(f\in C_c(\partial\Omega)\), the Perron solution \(u\) of
\[
    Lu=0\quad\text{in }\Omega,\qquad u=f\quad\text{on }\partial\Omega
\]
satisfies
\[
    u(x)=\int_{\partial\Omega} f\,d\omega_L^x.
\]
We write $\omega_L^x$ and \(\omega_{L^*}^x\) for the elliptic measure of $L$ or the adjoint operator $L^*$ respectively. If clear from context, we might also drop the subscript.
\end{defin}

The following estimate goes back to Bourgain's work in \cite{bourgain_hausdorff_1987}, and can be found in a geometric setting that includes ours in \cite[Lemma 11.21]{heinonen_nonlinear_1993}. For elliptic operators with BMO antisymmetric part, no particular adaptations are required, as observed in \cite[Prop. 6.2]{li_lp_2019}.
\begin{lemma}[Bourgain's estimate]\label{lemma:BourgainsEstimate}
    Let $\Omega\subset\mathbb{R}^{n+1}$ be a CAD, and $L$ and elliptic operator with BMO antisymmetric part. Then there exists $c>0$ such that for all $x\in\partial\Omega, r>0$
    \[\omega_L^{A(\Delta(x,r))}(\Delta(x,r))\geq c.\]
\end{lemma}

\begin{lemma}[Green's function and elliptic measure]\label{lem:green-elliptic-measure-comparison}
Let \(\Delta = \Delta(Q, r) \subset\partial\Omega\) and let \(A_\Delta\) be the corkscrew point above \(\Delta\).
If \(x\in\Omega\setminus B(Q,4r)\), then
\begin{equation}\label{eq:green-measure-comparison}
    C^{-1}r^{\,n-1}G(x,A_\Delta)
    \leq
    \omega_L^x(\Delta)
    \leq
    Cr^{\,n-1}G(x,A_\Delta).
\end{equation}
The analogous estimate with \(G\) and \(\omega_*^x\) also holds.
Equivalently, if \(y\in\Omega\), \(y'\in\partial\Omega\) satisfies \(|y-y'|=\delta(y)\), \(\Delta_{y'}=\Delta(y',\delta(y))\), and \(x\notin B(y',4\delta(y))\), then
\begin{equation}\label{eq:green-over-delta-measure-density}
    \frac{G^*(x,y)}{\delta(y)}
    \approx
    \frac{\omega_L^x(\Delta_{y'})}{\sigma(\Delta_{y'})}.
\end{equation}
\end{lemma}

\begin{proof}
We start by proving the first bound in  \eqref{eq:green-measure-comparison}. Consider the two functions
\[u(x) = r^{n-1}G(x, A_\Delta), \quad v(x) = \omega_L^x(\Delta).\]
By definition of the corkscrew point, 
\[\text{dist}(B(A_\Delta,  c_1r) , \partial \Omega) \gtrsim r.\]  We will show that $u \leq Cv$ on the region $\sq{\Omega} := \Omega \setminus B(A_\Delta, c_1 r)$. On this region, both $u$ and $v$ solve $Lu = Lv = 0$. We will apply the maximum principle to $u - Cv$ on $\sq{\Omega}$ where $C$ is a suitably large constant to be chosen later. The boundary of $\sq{\Omega}$ consists of two parts: 
\[\partial\sq{\Omega} = \partial \Omega \cup \partial B(A_\Delta, c_1r).\]
On $\partial \Omega$, $u = 0$ so we trivially have $u - Cv \leq 0$ for any $C > 0$. On $\partial B(A_\Delta, c_1r)$, we know that 
\[u(x) \lesssim r^{n-1}|x - A_\Delta|^{1 - n} \lesssim 1\]
and by Bourgain's estimate \reflemma{lemma:BourgainsEstimate} and Harnack \reflemma{Lemma:LocalEstimatesForBMO}
\[1\lesssim \omega^{A_{\Delta}}(\Delta) \lesssim \omega^x(\Delta) = v(x).\]
Hence \(u(x)\leq C v(x)\) also on $\partial B(A_\Delta, c_1r)$, if $C$ is chosen sufficiently large. 
Since 
\[\limsup_{|x| \to \infty}(u - Cv) \leq \limsup_{|x| \to \infty}u = 0,\]
by the maximum principle, we have that $u \leq Cv$ on $\sq{\Omega}$ as wanted. The proof of the second inequality can be found in \cite{li_pipher_bmo}.

\end{proof}

\begin{cor}\label{cor:green-measure-doubling}
The measures \(\omega_L^x\) and \(\omega_{L^*}^x\) are doubling on surface balls whose dilates do not contain the pole.
In particular, if $Q \in \partial \Omega$ and \(x\notin B(Q,8r)\), then for $\Delta = \Delta(Q, r)$ we have
\[
    \omega_L^x(2\Delta)\lesssim \omega_L^x(\Delta),
    \qquad
    \omega_{L^*}^x(2\Delta)\lesssim \omega_{L^*}^x(\Delta).
\]
\end{cor}

\begin{proof}
The proof follows from \reflemma{Lemma:LocalEstimatesForBMO} and \reflemma{lem:green-elliptic-measure-comparison}:
\[\omega^x(2\Delta) \approx r^{n - 1}G(x, A_{2\Delta}) \approx r^{n-1}G(x, A_\Delta) \approx \omega^x(\Delta)\]
and similarly for $\omega_{L^*}^x$.
    
\end{proof}

\begin{defin}[Muckenhoupt and Reverse Hölder classes]
    For $p>1$, we say $\omega_L\in B_p(\sigma)$, if there exists $C>0$ such that for all boundary balls $\Delta\subset\partial\Omega$ and all $x\in \Omega$ with $k^x:=\frac{\partial\omega_{L}^x}{\partial\sigma}$
    \[\Big(\fint_\Delta (k^x)^p\,d\sigma\Big)^{1/p}\leq C\fint_\Delta k^x\,d\sigma\]
    holds. Furthermore, we define $A_\infty(\sigma):=\bigcup_{p>1}B_p(\sigma)$.
\end{defin}

\begin{lemma}\label{lemma:IntOver1/delta}
    Let $\Omega\subset\mathbb{R}^{n+1}$ be a CAD. Fix a small $0<\tau\leq 1$, some $\alpha>2$ and $E\subset\Delta$ nonempty, where $\Delta\subset\partial\Omega$ is a boundary ball. Let $\sq{\Omega}_j:=\sq{\Omega}_j^{(\tau \mathrm{diam}(\Delta),\alpha)}(E)$ as defined in Section \ref{Section:Defphi_W}. Then
    \begin{align}
        \label{eq:GeoLemmaInt1}\int_{\sq{\Omega}_{1}\setminus\sq{\Omega}_0}\frac{1}{\delta}\, dx\leq C(n,\alpha, C_0,\tau)\sigma(\Delta),
    \end{align}
    holds. If additionally, $y\in \Omega$ is such that $\mathrm{dist}(y,\Delta)>C(\alpha,C_0,c_1)\mathrm{diam}(\Delta)$, then
    \begin{align}\label{eq:GeoLemmaInt2}\int_{\sq{\Omega}_{1}\setminus\sq{\Omega}_0}\frac{G(y,x)}{\delta(x)^2}\, dx \leq C(n,\lambda_0,\Lambda_0,\alpha, C_0,c_1,c_0\tau) \int_{3\Delta}\mathcal{M}[k^y\chi_{3\Delta}] d\sigma,
    \end{align}
    where $k^y=\frac{d\omega^y}{d\sigma}$ is the Radon-Nikodym derivative, and $C_0$, $c_1$, and $c_0$ come from the Ahlfors regularity of $\partial\Omega$,interior corkscrew condition of $\Omega$, and the Harnack chain condition, respectively.
\end{lemma}

\begin{proof}
    First, we set \(R=\tau\mathrm{diam}(\Delta)\).
Since \(\tau\leq1\), we have \(R\leq \mathrm{diam}(\Delta)\).

First, the geometric observation. If \(x\in \sq{\Omega}_1\setminus \sq{\Omega}_0\) and \(\delta(x)\leq R\), then
\begin{align}
|x-p|\geq(1+\alpha)\delta(x)
\qquad\text{for every }p\in E,
\label{eq:proofgeolemma2}
\end{align}
because otherwise \(x\in\sq{\Omega}_0\). On the other hand, membership in \(\sq{\Omega}_1\) gives some \(p_x\in E\) such that
\begin{align}
|x-p_x|<(1+2\alpha)\delta(x).
\label{eq:proofgeolemma3}
\end{align}
Choose \(x^*\in\partial\Omega\) with \(|x-x^*|=\delta(x)\). Equations \eqref{eq:proofgeolemma2} and \eqref{eq:proofgeolemma3} imply
\begin{align}
\alpha\delta(x)
\leq\operatorname{dist}(x^*,E)
<(2+2\alpha)\delta(x).
\label{eq:proofgeolemma4}
\end{align}
Set
\[
s=\frac{R}{100(1+\alpha)},
\]
then we can split the integrals \eqref{eq:GeoLemmaInt1} and \eqref{eq:GeoLemmaInt2} into close and away parts, i.e.
\begin{align}\label{eq:awaypart1}
        \int_{\sq{\Omega}_{1}\setminus\sq{\Omega}_0}\frac{1}{\delta}\, dx=\int_{(\sq{\Omega}_{1}\setminus\sq{\Omega}_0)\cap \{x\in\Omega:\delta(x)<s\}}\frac{1}{\delta}\, dx + \int_{(\sq{\Omega}_{1}\setminus\sq{\Omega}_0)\cap \{x\in\Omega:\delta(x)\geq s\}}\frac{1}{\delta}\, dx
    \end{align}
and 
\begin{align}\label{eq:awaypart2}\int_{\sq{\Omega}_{1}\setminus\sq{\Omega}_0}\frac{G(y,x)}{\delta(x)^2}\, dx =\int_{(\sq{\Omega}_{1}\setminus\sq{\Omega}_0)\cap \{x\in\Omega:\delta(x)< s\}}\frac{G(y,x)}{\delta(x)^2}\, dx+\int_{(\sq{\Omega}_{1}\setminus\sq{\Omega}_0)\cap \{x\in\Omega:\delta(x)\geq s\}}\frac{G(y,x)}{\delta(x)^2}\, dx.
\end{align}
Since \(|(\sq{\Omega}_{1}\setminus\sq{\Omega}_0)\cap \{x\in\Omega:\delta(x)\geq s\}|\lesssim \mathrm{diam}(\Delta)^{n+1}\), we obtain for the away part in \eqref{eq:awaypart1}
\begin{align}\label{eq:awayPart1Bound}
\int_{(\sq{\Omega}_{1}\setminus\sq{\Omega}_0)\cap \{x\in\Omega:\delta(x)\geq s\}}\frac{1}{\delta}\, dx\leq \int_{(\sq{\Omega}_{1}\setminus\sq{\Omega}_0)\cap \{x\in\Omega:\delta(x)\geq s\}}\frac{1}{s}\, dx\lesssim_{\alpha,n,C_0} \frac{1}{\tau}\sigma(\Delta)\lesssim_\tau\sigma(\Delta).
 \end{align}
Furthermore, for the away part in \eqref{eq:awaypart2}, by \reflemma{Lemma:LocalEstimatesForBMO} and \reflemma{lem:green-elliptic-measure-comparison} for \(y\in \Omega\) with \(\mathrm{dist}(y,\Delta)>C(\alpha, C_0,c_1,c_0)\mathrm{diam}(\Delta)\), where \(C(\alpha, C_0,c_1,c_0)>0\) is chosen sufficiently large, 
\begin{equation}\label{eq:awayPart2Bound}\begin{split}\int_{(\sq{\Omega}_{1}\setminus\sq{\Omega}_0)\cap \{x\in\Omega:\delta(x)\geq s\}}\frac{G(y,x)}{\delta(x)^2}\, dx\lesssim_{\tau,\alpha,C_0,c_1} \int_{(\sq{\Omega}_{1}\setminus\sq{\Omega}_0)\cap \{x\in\Omega:\delta(x)\geq s\}}\frac{G(y,A_\Delta)}{\delta(x)^2}\, dx
\\
\lesssim \frac{\omega^y(\Delta)}{\sigma(\Delta)}\int_{(\sq{\Omega}_{1}\setminus\sq{\Omega}_0)\cap \{x\in\Omega:\delta(x)\geq s\}}\frac{\mathrm{diam}(\Delta)}{\delta(x)^2}\, dx\lesssim \int_{3\Delta}\mathcal{M}[k^y\chi_{3\Delta}] d\sigma.
\end{split}\end{equation}

\noindent Next, let us turn to bounding the close parts in \eqref{eq:awaypart1} and \eqref{eq:awaypart2}, respectively.
Recall the Whitney decomposition from Section~\ref{Section:Defphi_W}, and set
\[
\mathcal D_{\mathrm{low}}:
=\bigl\{W\in\sq{\mathcal D}_{0,1}:
W\cap\{x\in \sq{\Omega}_1\setminus\sq{\Omega}_0|\delta(x)<s\}\neq\varnothing\bigr\}.
\]
Then 
\[\{x\in \sq{\Omega}_1\setminus\sq{\Omega}_0|\delta(x)<s\}\subset\bigcup_{W\in\mathcal{D}_{\mathrm{low}}}W.\]
For every \(W\in\mathcal D_{\mathrm{low}}\), fix
\[
x_W\in W\cap (\sq{\Omega}_1\setminus\sq{\Omega}_0),\qquad\textrm{with} \qquad\delta(x_W)<s,
\]
and choose \(x_W^*\in\partial\Omega\) with \(|x_W-x_W^*|=\delta(x_W)\). Define
\[
\Delta_W=B(x_W^*,\delta(x_W)/10)\cap\partial\Omega.
\]
By the properties of the Whitney cubes and the Ahlfors regularity of \(\partial\Omega\), we obtain
\begin{align}
\sigma(\Delta_W)\approx \delta(x_W)^n\approx l(W)^n.
\label{eq:proofgeolemma6}
\end{align}
Moreover, \eqref{eq:proofgeolemma4} gives, for every \(z\in \Delta_W\),
\begin{align}
(\alpha-\tfrac1{10})\delta(x_W)
\leq\operatorname{dist}(z,E)
\leq(2+2\alpha+\tfrac1{10})\delta(x_W).
\label{eq:proofgeolemma7}
\end{align}
These balls lie in \(3\Delta\). Indeed, using the point \(p_{x_W}\) from \eqref{eq:proofgeolemma3} and denoting the center point of \(\Delta\) as \(p_0\) and the radius as \(r\),
\[
|x_W^*-p_0|
\leq |x^*_W-x_W|+|x_W-p_{x_W}|+|p_{x_W}-p_0|
<r+(2+2\alpha)\delta(x_W).
\]
Consequently, for \(z\in \Delta_W\), the choice of \(s\) gives
\begin{align}
|z-p_0|
<r+(2+2\alpha+\tfrac1{10})s<3r.
\label{eq:proofgeolemma8}
\end{align}
The same calculation shows that
\begin{align}
\Delta(z,\delta(x_W))\subset3\Delta
\qquad\textrm{for all }z\in \Delta_W.
\label{eq:proofgeolemma9}
\end{align}
We next prove
\begin{align}
\sum_{W\in\mathcal D_{\mathrm{low}}}\chi_{\Delta_W}
\leq C(n,\alpha).
\label{eq:proofgeolemma10}
\end{align}
Fix \(z\in\partial\Omega\), and consider cubes \(W\in \mathcal{D}_{\mathrm{low}}\) for which \(z\in \Delta_W\). If \(\operatorname{dist}(z,E)=0\), \eqref{eq:proofgeolemma7} shows that there are no such cubes. Otherwise, \eqref{eq:proofgeolemma7} gives
\[
\frac{\mathrm{dist}(z,E)}{2+2\alpha+\tfrac1{10}}
\leq \delta(x_W)
\leq\frac{\mathrm{dist}(z,E)}{\alpha-\tfrac1{10}}.
\]
Hence their side lengths are comparable, with constants depending on \(\alpha\). Also,
\begin{align}\label{eq:proofgeolemmaNew}
|x_W-z|\leq |x_W-x^*_W|+|x_W^*-z|
\leq\tfrac{11}{10}\delta(x_W).
\end{align}
This places all these cubes \(W\) in a ball of radius comparable to their common scale. Since the Whitney cubes are disjoint and each has volume bounded below by a fixed multiple of that scale to the power \(n+1\), their number is uniformly bounded. This proves \eqref{eq:proofgeolemma10}.

Now let us return to the close part in \eqref{eq:awaypart1}.
It follows from \eqref{eq:proofgeolemma6}, and \eqref{eq:proofgeolemma10} that
\[
\begin{aligned}
\sum_{W\in\mathcal D_{\mathrm{low}}}
\int_W\frac{dx}{\delta(x)}
&\lesssim
\sum_{W\in\mathcal D_{\mathrm{low}}}\delta(x_W)^n
\lesssim
\sum_{W\in\mathcal D_{\mathrm{low}}}\sigma(\Delta_W)
\lesssim \sigma(3\Delta)
\lesssim \sigma(\Delta),
\end{aligned}
\]
which combined with \eqref{eq:awayPart1Bound} implies \eqref{eq:GeoLemmaInt1}.
For the close part in the second integral \eqref{eq:awaypart2}, let
\(x\in W\) for \(W\in \mathcal{D}_{\mathrm{low}}\). Furthermore, let \(y\in \Omega\) such that $\mathrm{dist}(y,\Delta)>C(\alpha,C_0,c_1)\mathrm{diam}(\Delta)$, where \(C(\alpha,C_0,c_1)\) is chosen such that the following argument using \reflemma{lem:green-elliptic-measure-comparison} and the doubling property \refcor{cor:green-measure-doubling} apply.   Then, by \eqref{eq:proofgeolemmaNew}

\begin{align}
\frac{G(y,x)}{\delta(x)}
\lesssim
\frac{\omega^y(\Delta(x^*,\delta(x)))}
{\delta(x)^n}
\lesssim
\frac{\omega^y(\Delta(z,C\delta(x_W)))}{\delta(x_W)^n}
\lesssim
\frac{\omega^y(\Delta(z,\delta(x_W)))}{\delta(x_W)^n}.\label{eq:g/delta1}
\end{align}
 By \eqref{eq:proofgeolemma9}, the final ball is contained in \(3\Delta\), so Ahlfors regularity of \(\partial\Omega\) gives
\begin{align}
\frac{\omega^y(\Delta(z,\delta(x_W)))}{\delta(x_W)^n}
\lesssim
\fint_{\Delta(z,\delta(x_W))}k^y\,d\sigma
\leq \mathcal M[k^y\chi_{3\Delta}](z).\label{eq:g/delta2}
\end{align}
Consequently,
\[
\begin{aligned}
\int_W\frac{G(y,x)}{\delta(x)^2}\,dx
&\lesssim \delta(x_W)^n\inf_{z\in \Delta_W}\mathcal M[k^y\chi_{3\Delta}](z)
&\lesssim\int_{\Delta_W}\mathcal M[k^y\chi_{3\Delta}]\,d\sigma.
\end{aligned}
\]
Summing and using \eqref{eq:proofgeolemma10},
\begin{align}
\sum_{W\in\mathcal D_{\mathrm{low}}}
\int_W\frac{G(y,x)}{\delta(x)^2}\,dx
\lesssim\int_{3\Delta}\mathcal{M}[k^y\chi_{3\Delta}]\,d\sigma.
\label{eq:proofgeolemma16}
\end{align}
Combining this with \eqref{eq:awayPart2Bound} implies \eqref{eq:GeoLemmaInt2}.
\end{proof}

\section{Solvability of the Poisson Dirichlet problem}\label{sec:ProofSquare}

\subsection{Square function bound}

In this section we prove that solvability of $(\operatorname{D})_{p'}^{L^*}$ implies a square function bound of the type \eqref{eq:SquareFunctionEstimate}. Throughout we assume that \(L\) has antisymmetric part in BMO. 

\begin{lemma}\label{lemma:SquareFunctionEstimate}
Let $\Omega\subset\mathbb{R}^{n+1}$ be a CAD, let \(v\in Y_0^{1,2}(\Omega)\) be the solution to \(L^*v=H-\mathrm{div}(F)\) for \(H,F\in L^\infty_c(\Omega)\) and let \(\omega_{L^*}\in A_\infty(\sigma)\). Then for every \(1<q<\infty\)
\begin{align}\Vert \mathcal{S}(v)\Vert_{L^q(\partial\Omega)}&\leq C_q\big(\Vert \mathcal{N}(v)\Vert_{L^q(\partial\Omega)}+\Vert \mathcal{T}_{2}(\delta H)\Vert_{L^{q}(\partial\Omega)}+ \Vert \mathcal{T}_{2}(F)\Vert_{L^{q}(\partial\Omega)}\big).
\end{align}
\end{lemma}
To prove \reflemma{lemma:SquareFunctionEstimate}, we use the following "good \(\lambda\) inequality". For that, we first assume without loss of generality that \(0\in \partial\Omega\) and define a truncated version of the square function by
\[\mathcal{S}_{\alpha,\rho}(v)^2(x)\coloneqq \int_{\Gamma_\alpha(x)\cap K_\rho}|\nabla v(y)|^2\delta(y)^{1-n}(y)dy,\quad\textrm{ where }\quad K_\rho:=B(0,\frac{1}{\rho})\cap \{z:\delta(z)\geq \rho\}.\]
Furthermore, for \(\beta>0\) and \(x\in\partial\Omega\), set
\[D_\beta:=\{\mathcal{N}_{8\alpha}(v)+\sq{\mathcal{N}}_{8\alpha,2}(\delta|\nabla v|) + \mathcal{T}_{8\alpha,2}(\delta H)(x)+\mathcal{T}_{8\alpha,2}(F)>\beta\},\quad\textrm{ and }\quad q(\beta,x):=\mathcal{M}(\chi_{D_\beta}). \]
We prove the following "good \(\lambda\) inequality".

\begin{lemma}\label{lemma:GoodLambdaIneq}
     Let $\Omega\subset\mathbb{R}^{n+1}$ be a CAD, let \(v\in Y_0^{1,2}(\Omega)\) be the solution to \(L^*v=H-\mathrm{div}(F)\) for \(H,F\in L^\infty_c(\Omega)\), and let \(\omega_{L^*}\in A_\infty(\sigma)\).
    For all apertures \(\alpha>2\) there exists \(C,\eta,\kappa>0\) such that for all sufficiently small \(\gamma>0\) and all \(\lambda,\rho>0\)
    \begin{align*}\sigma(\{x\in\partial\Omega:\mathcal{S}_{\alpha,\rho}(v)(x)>2\lambda; q(\gamma\lambda,x)
    \leq \kappa\})
    \leq C\gamma^\eta\sigma(\{\mathcal{S}_{8\alpha,\rho}(v)>\lambda\}).\end{align*}
\end{lemma}

First, let us prove \reflemma{lemma:SquareFunctionEstimate} assuming \reflemma{lemma:GoodLambdaIneq}.
\begin{proof}[Proof of \reflemma{lemma:SquareFunctionEstimate}]
By the definition of \(q(\gamma\lambda,x)\) and \(L^1-L^{1,\infty}\) boundedness of the maximal function, we have
\begin{align*}
\sigma(\{q(\gamma\lambda,\cdot)>\kappa\})
&=\sigma(\{\mathcal M_\sigma(\chi_{ D_{\gamma\lambda}})>\kappa\})
\lesssim\|\chi_{D_{\gamma\lambda}}\|_{L^1(\sigma)}
\\
&\lesssim\sigma(\{\mathcal{N}_{8\alpha}(v)+\sq{\mathcal{N}}_{8\alpha,2}(\delta|\nabla v|) + \mathcal{T}_{8\alpha,2}(\delta H)(x)+\mathcal{T}_{8\alpha,2}(F)>\gamma\lambda\}).
\end{align*}
Thus, the good $\lambda$ inequality \reflemma{lemma:GoodLambdaIneq} implies 
\begin{align*}
    \int_{\partial\Omega}\mathcal{S}_{\alpha,\rho}(v)^qd\sigma \lesssim C_q\Big(\int_{\partial\Omega}\mathcal{N}_{8\alpha}(v)^qd\sigma+\int_{\partial\Omega}\sq{\mathcal{N}}_{8\alpha}(\delta|\nabla v|)^qd\sigma+\int_{\partial\Omega}\mathcal{T}_{8\alpha,2}(\delta H)^qd\sigma
    \\
    +\int_{\partial\Omega}\mathcal{T}_{8\alpha,2}(F)^qd\sigma+\gamma^\eta\int_{\partial\Omega}\mathcal{S}_{8\alpha,\rho}(v)^qd\sigma\Big).
\end{align*}
By \reflemma{lemma:HigherIntegrabilityForGradient}, we can bound $\sq{\mathcal{N}}_{8\alpha}(\delta|\nabla v|)\lesssim \sq{\mathcal{N}}_{9\alpha}(v)+\mathcal{T}_{9\alpha,2}(F) + \mathcal{T}_{9\alpha,2}(\delta H)$ pointwise. Furthermore, recall that the $L^q$ norms of the nontangential maximal function are comparable for different apertures $\alpha$. Moreover, the $L^q$ norms of the square function and the discretized area function $\mathcal{T}$ are also comparable for different apertures $\alpha$. Its proof goes back to \cite{coifman_new_1985}, and for area functions on CAD we can refer to \cite[Prop. 4.5]{milakis_harmonic_2013} and \refcor{cor:appendixDiffApCompT}. Thus, we get
\begin{align*}
    \int_{\partial\Omega}\mathcal{S}_{\alpha,\rho}(v)^qd\sigma \lesssim C_q\Big(\int_{\partial\Omega}\mathcal{N}_{\alpha}(v)^qd\sigma+\int_{\partial\Omega}\mathcal{T}_{\alpha,2}(\delta H)^qd\sigma
    \\
    +\int_{\partial\Omega}\mathcal{T}_{\alpha,2}(F)^qd\sigma+\gamma^\eta\int_{\partial\Omega}\mathcal{S}_{\alpha,\rho}(v)^qd\sigma\Big).
\end{align*}
Choosing $\gamma$ sufficiently small allows us to hide the last term on the left hand side. Taking the limit in \(\rho\) to zero finishes the proof of \reflemma{lemma:SquareFunctionEstimate}.

\end{proof}

Hence, it remains to prove \reflemma{lemma:GoodLambdaIneq}.

\begin{proof}[Proof of \reflemma{lemma:GoodLambdaIneq}]

We begin by considering a boundary Whitney cube decomposition \(\{\Delta_j\}_j\) of \(\{\mathcal{S}_{\alpha,\rho}(v)>\lambda\}\subset\partial\Omega\). According to \cite{christ_tb_1990} and as in Section~\ref{Section:Defphi_W} for $\Omega$, we can ask for the following properties:
\begin{enumerate}
\item $\Delta_j$ is a boundary ball, i.e. there exists $x_j\in \partial\Omega, r_j>0$ such that $\Delta_j=B(x_j,r_j)\cap\partial\Omega$;
\item $\frac{1}{4}\mathrm{dist}(\Delta_j,\{\mathcal{S}_{\alpha,\rho}(v)\leq\lambda\})\leq\mathrm{diam}(\Delta_j)\leq 4\mathrm{dist}(\Delta_j,\{\mathcal{S}_{\alpha,\rho}(v)\leq\lambda\})$; and
\item $\chi_{\{\mathcal{S}_{\alpha,\rho}(v)>\lambda\}}\leq\sum_j \chi_{\Delta_j}\leq N\chi_{\{\mathcal{S}_{\alpha,\rho}(v)>\lambda\}}$, where $N\in\mathbb{N}$ is a fixed constant.
\end{enumerate}

Let us now fix such a boundary Whitney cube $\Delta_j$ and drop the subscript $j$. We define 
\[\mathcal{E} :=\{x\in\partial\Omega: \mathcal{S}_{\alpha,\rho}(v)(x)>2\lambda, q(\gamma\lambda,x)\leq \kappa\},\]
and \(E:=E_\Delta:=\Delta\cap\mathcal{E}\). 
We claim that it is enough to prove that there exists a $C>0$ independent of $\Delta$ and $\gamma$ such that
\begin{align}\omega_{L^*}^{X_\Delta}(E)\leq C \gamma^2\omega_{L^*}^{X_\Delta}(\Delta),\label{eq:ClaimInGoodLambdaInequality}\end{align}
where $X_\Delta$ is a point that is chosen as follows:
If \(\Omega\) is bounded, we can set \(X_\Delta\) to be a ``center point" of \(\Omega\), i.e a point with \(\mathrm{diam}(\Omega)\lesssim \delta(X_\Delta)\). If \(\Omega\) is unbounded, we can set \(X_\Delta\) to be a corkscrew point of the enlargement \(4c_0^{-1}\Delta\). This choice of $X_\Delta$ is to guarantee that tools like \reflemma{lem:green-elliptic-measure-comparison} and the doubling property of the harmonic measure can be applied with uniform constants independent of $\Delta$ or $X_\Delta$.

Assuming \eqref{eq:ClaimInGoodLambdaInequality}, \(\omega_{L^*}\in A_\infty(\sigma)\) implies that there exist uniform constants $\sq{C}>0$ and $\eta>0$ such that 
\[\frac{\sigma(E)}{\sigma(\Delta)}\leq \sq{C}\Big(\frac{\omega_{L^*}^{X_\Delta}(E)}{\omega_{L^*}^{X_\Delta}(\Delta)}\Big)^\eta,\]
which implies
\(\sigma(E_\Delta)\leq \sq{C}\gamma^\eta\sigma(\Delta)\), and the finite overlap of the Whitney cubes allows to sum over them to get
\begin{align*}
    \sigma(\mathcal{E})&\leq\sum_{j}\sigma(E_{\Delta_j})
    \leq \sum_{j}\gamma^\eta\sigma(\Delta_j)\lesssim \sq{C}\gamma^\eta\sigma(\{\mathcal{S}_{8\alpha,\rho}(v)>\lambda\}),
\end{align*}
which proves the "good \(\lambda\) inequality". Thus, it remains to prove \eqref{eq:ClaimInGoodLambdaInequality}.

\smallskip

To prove \eqref{eq:ClaimInGoodLambdaInequality}, we first abbreviate notation by writing the corresponding elliptic measure of \(L^*\) to the point \(X_\Delta\) by \(\sq{\omega}=\omega_{L^*}^{X_\Delta}\), the Green's function as \(\sq{G}=G_{L^*}(X_\Delta,\cdot)\), and $\sq{k}:=\frac{d\sq{\omega}}{d\sigma}$. Moreover, set $l:=\mathrm{diam}(\Delta)$. Let us also recall the definition of the Whitney cube decompositions, the sawtooth regions and the definitions of the $\phi_W$ from Section~\ref{Section:Defphi_W}.

Lemma 1 in \cite{dahlberg_area_1984} is formulated on Lipschitz domains but the proof holds verbatim on CAD and thus depends on the CAD constants. Using this, we have for every \(\tau>0\) and all sufficiently small \(\gamma>0\) 
    that the truncated square function is bounded below by \(\lambda/2\) on \(E\), i.e. that
\begin{align}\mathcal{S}_{\alpha,\rho}^{\tau l}(v)^2(x)\coloneqq \int_{\Gamma_\alpha^{\tau l}(x)\cap K_\rho}|\nabla v(y)|^2\delta(y)^{1-n}dy>\frac{\lambda^2}{4}\label{eq:TruncatedSquareFctBoundedBelow}\end{align}
    holds for all points \(x\in E\). We choose $\tau=1$.

Now, \eqref{eq:TruncatedSquareFctBoundedBelow} and Fubini yields
\begin{align*}
    \sq{\omega}(E)&\lesssim \frac{1}{\lambda^2}\int_{E} \mathcal{S}_{\alpha,\rho}^{l}(v)^2d\sq{\omega}
    \leq\frac{1}{\lambda^2}\int_{E}\Big(\int_{\sq{\Omega}_1^{(l,\alpha)}(E)}|\nabla v|^2\delta^{-n+1}\chi_{\Gamma_\alpha^{ l}(x)}dy\Big)d\sq{\omega}(x)
    \\
    &=\frac{1}{\lambda^2}\int_{\sq{\Omega}_1^{(l,\alpha)}(E)}|\nabla v|^2\delta^{-n+1}\Big(\int_{E}\chi_{\Gamma_\alpha^{l}(x)}d\sq{\omega}(x)\Big)dy.
\end{align*}
From here on, all sawtooth regions $\sq{\Omega}_j^{(\tau,\alpha)}$ are with respect to the good set $E$, and we drop the argument.
Let \(y^*\in\partial\Omega\) denote a point satisfying \(|y-y^*|=\delta(y)\), then \reflemma{lem:green-elliptic-measure-comparison} and the doubling property of the elliptic measure implies that
\begin{align*}
    &\frac{1}{\lambda^2}\int_{\sq{\Omega}_2^{(l,\alpha)}}|\nabla v|^2\delta^{-n+1}\Big(\int_{E}\chi_{\Gamma^{l}(x)}d\sq{\omega}(x)\Big)dy
    \lesssim \frac{1}{\lambda^2}\int_{\sq{\Omega}_2^{(l,\alpha)}}|\nabla v|^2\delta^{-n+1}\sq{\omega}(\Delta(y^*,2\alpha\delta(y)))dy
    \\
    &\lesssim \frac{1}{\lambda^2}\int_{\sq{\Omega}_2^{(l,\alpha)}}|\nabla v|^2\delta^{-n+1}\sq{\omega}(\Delta(y^*,\delta(y)))dy\lesssim \frac{1}{\lambda^2}\int_{\sq{\Omega}_2^{(l,\alpha)}}|\nabla v|^2\sq{G}dy.
\end{align*}
Next, we choose  a smooth cut-off \(\psi\) with the following properties
\begin{enumerate}
    \item \(\psi\in C^\infty(\Omega)\),
    \item \(0\leq \psi\leq 1\),
    \item \(\psi\equiv 1\) on \(\sq{\Omega}_2^{(l,\alpha)}\) and \(\psi\equiv 0\) outside of  \(\sq{\Omega}_3^{(l,\alpha)}\),
    \item \(|\nabla \psi(y) |\lesssim \frac{1}{\delta(y)}\) for all \(y\in \sq{\Omega}_3^{(l,\alpha)}\setminus\sq{\Omega}_2^{(l,\alpha)}\).
\end{enumerate}
Utilizing the ellipticity condition, and the cutoff $\psi$, we obtain
\begin{align*}
    \frac{1}{\lambda^2}\int_{\sq{\Omega}_2^{(l,\alpha)}}|\nabla v|^2\sq{G} dy
    &\lesssim \frac{1}{\lambda^2}\int_{\Omega} A^*\nabla v\cdot \nabla v \sq{G}\psi dy.
\end{align*}
We apply the following algebraic identity to the integrand to split it into four pieces:
\begin{align*}
    A^* \nabla v \cdot \nabla v \sq{G} \psi & = A^*\nabla v\cdot\nabla(v\sq{G}\psi) -\frac12 A\nabla \sq{G}\cdot\nabla(v^2\psi) \\
    & -\frac12 A^a\nabla\psi\cdot\nabla(v^2\sq{G}) +\frac12 A^s\nabla\psi\cdot
    \bigl(v^2\nabla \sq{G}-\sq{G}\nabla(v^2)\bigr) \\
    & =: T_1 + T_2 + T_3 + T_4.
\end{align*}
The first term is handled by the PDE for $v$ to convert the integral into an integral for the data $H,F$. The second term is handled by the PDE for $\sq{G}$. The third term is handled by a dyadic decomposition and the BMO condition on $A^a$. The fourth term is handled by the boundedness of $A^s$. We provide the details below.

\textbf{Bound for $T_2$:} We will show that $\int_\Omega T_2 \, dy = 0$ by using the PDE for $\sq{G}$. To do this, we just need to show that $v^2 \psi \in W^{1,2}_0(\Omega)$. First, we note that since \(|v|\lesssim \sup_E \mathcal{N}(v)\lesssim \lambda\gamma\) on \(\supp(\psi)\), we have that
\[|\nabla (v^2\psi)|\lesssim |\nabla v| |v|+|v|^2|\nabla \psi| \lesssim (\lambda \gamma) |\nabla v| + (\lambda \gamma) \frac{|v|}{\delta}.\]
Using the Hardy inequality from \cite{Lewis_1988}\footnote{The hypothesis from \cite{Lewis_1988} require uniform fatness of the complement of the domain \(\Omega\), which holds even for 1-sided chord-arc domains.
}, we have $\|v / \delta\|_{L^2(\Omega)} \lesssim \|\nabla v\|_{L^2(\Omega)}$ and so \(v^2\psi\in W_0^{1,2}(\Omega)\), and since the pole of the Green's function \(\sq{G}\) lies away from the support of \(\psi\), the integral over \(T_2\) vanishes.

\textbf{Bound for $T_1$:} For the integral over \(T_1\), we note that 
\[|\nabla (v\sq{G}\psi)|\lesssim |\nabla v|+ \gamma\lambda|\nabla \sq{G}|\chi_{\mathrm{supp}(\psi)}+\frac{|v|}{\delta}.\] 
Thus, using the same reasoning as in the bound for $T_2$, \(v\sq{G}\psi\in W^{1,2}_0(\Omega)\), and we can apply the PDE for \(v\) and obtain
\begin{align*}
    \frac{1}{\lambda^2}\int_\Omega T_1 \, dy&=\frac{1}{\lambda^2}\int_{\Omega}Hv\sq{G}\psi+ F \cdot\nabla v \sq{G}\psi + F \cdot\nabla \sq{G} v\psi + F \cdot\nabla\psi v \sq{G} \,dy
\end{align*}

Since for all \(x\in E\) we have \(\kappa\geq q(\gamma\lambda,x)=\mathcal{M}(\chi_{D_{\gamma\lambda}})(x)\), it follows that 
\[\frac{\sigma(B(x,r)\setminus D_{\gamma\lambda})}{\sigma(B(x,r))}\geq 1-\kappa\qquad \textrm{ for all }r>0.\]
Thus \(x\) is a global \(1-\kappa\)-density point of the surface measure on \(\partial \Omega\setminus D_{\gamma\lambda}\). For $\kappa$ sufficiently small, we can apply
\cite[Lemma~3.14]{milakis_harmonic_2013} and obtain
\begin{align*}
    \frac{1}{\lambda^2}\int_\Omega T_1 \, dy & \lesssim\frac{1}{\lambda^2}\int_{8\Delta\setminus D_{\gamma\lambda}}\int_{\Gamma^{4l}_{4\alpha}(P)} \delta^{-n}|H v \,\sq{G}\psi| \, dy d\sigma(P)
      \\ 
      & + \frac{1}{\lambda^2}\int_{8\Delta\setminus D_{\gamma\lambda}}\int_{\Gamma^{4l}_{4\alpha}(P)} \delta^{-n}|F \cdot\nabla v \,\sq{G}\psi|\, dy d\sigma(P) 
    \\
    &\quad+ \frac{1}{\lambda^2}\int_{8\Delta\setminus D_{\gamma\lambda}}\int_{\Gamma^{4l}_{4\alpha}(P)}\delta^{-n}|F \cdot\nabla \sq{G}\, v\psi|\, dy d\sigma(P) \\ 
    & + \frac{1}{\lambda^2}\int_{8\Delta\setminus D_{\gamma\lambda}}\int_{\Gamma^{4l}_{4\alpha}(P)}\delta^{-n}|F\cdot \nabla\psi\, v \sq{G}| \, dy d\sigma(P)
    \\
    &=:J_1+J_2+J_3+J_4.
\end{align*}
Here, the restriction of the outer integrals from \(\partial\Omega\setminus D_{\gamma\lambda}\) to \(8\Delta\cap D_{\gamma\lambda}\) follows from the compact support of the integrand. 

By \eqref{eq:g/delta1} and \eqref{eq:g/delta2}, for \(W\cap\Gamma_{4\alpha}^{4l}(P)\neq\emptyset\) it holds that \(\frac{\sup_{y\in 2W}\sq{G}(y)}{l(W)}\lesssim \mathcal{M}[\sq{k}\chi_{8\Delta}](P)\) and $2W\subset \Gamma_{8\alpha}(P)$. Combining this with Hölder and \eqref{eq:BoundOFM[k]byomega} gives
\begin{align*}
    J_2&\lesssim \frac{1}{\lambda^2}\int_{8\Delta\setminus D_{\gamma\lambda}}\sq{\mathcal{N}}_{8\alpha}(\delta|\nabla v|)\sum_{\substack{W\in \mathcal{D}(\Omega),\\ W\cap \Gamma^{4l}_{4\alpha}(P)\neq \emptyset}}\frac{\sup_{y\in W}\sq{G}(y)}{l(W)}l(W)\Big(\fint_{W}|F|^2dy\Big)^{1/2} d\sigma(P)
    \\
    &\lesssim \frac{1}{\lambda^2}\int_{8\Delta\setminus D_{\gamma\lambda}} \mathcal{M}[\sq{k}\chi_{8\Delta}]\mathcal{T}_
    {8\alpha,2}(F)\sq{\mathcal{N}}_{8\alpha}(\delta|\nabla v|) d\sigma(P)\lesssim \gamma^2\int_{8\Delta} \mathcal{M}[\sq{k}\chi_{8\Delta}] d\sigma(P)\lesssim \gamma^2\sq{\omega}(\Delta).
\end{align*}

Similarly by Caccioppoli, we obtain 

\begin{align*}
    J_3&\lesssim \frac{1}{\lambda^2}\int_{8\Delta\setminus D_{\gamma\lambda}}\sum_{\substack{W\in \mathcal{D},\\ W\cap \Gamma^{4l}_{4\alpha}(P)\neq \emptyset}}l(W)\Big(\fint_{W}|F|^2dy\Big)^{1/2}\Big(\fint_{W}|\nabla \sq{G}|^2dy\Big)^{1/2}\mathcal{N}_{8\alpha}(v)d\sigma(P)
    \\
    &\lesssim \frac{1}{\lambda^2}\int_{8\Delta\setminus D_{\gamma\lambda}}\sum_{\substack{W\in \mathcal{D},\\ W\cap \Gamma^{4l}_{4\alpha}(P)\neq \emptyset}}l(W)\Big(\fint_{W}|F|^2dy\Big)^{1/2}\Big(\fint_{2W}\frac{|\sq{G}|^2}{l(W)^2}dy\Big)^{1/2}\mathcal{N}_{8\alpha}(v)d\sigma(P)
    \\
    &\lesssim\frac{1}{\lambda^2}\int_{8\Delta\setminus D_{\gamma\lambda}} \mathcal{T}_{8\alpha,2}(F)\mathcal{M}[\sq{k}\chi_{8\Delta}]\mathcal{N}_{8\alpha}(v) d\sigma\lesssim \gamma^2\int_{8\Delta} \mathcal{M}[\sq{k}\chi_{8\Delta}] d\sigma\lesssim \gamma^2\sq{\omega}(\Delta),
\end{align*}

and also

\begin{align*}
    J_4&\lesssim  \frac{1}{\lambda^2}\int_{8\Delta\setminus D_{\gamma\lambda}}\sum_{\substack{W\in \mathcal{D},\\ W\cap \Gamma^{4l}_{4\alpha}(P)\neq \emptyset}}\frac{\sup_{y\in W}\sq{G}(y)}{l(W)}l(W)\Big(\fint_{W}|F|^2dy\Big)^{1/2}\mathcal{N}_{8\alpha}(v)d\sigma(P)
    \\
    &\lesssim \frac{1}{\lambda^2}\int_{8\Delta\setminus D_{\gamma\lambda}} \mathcal{M}[\sq{k}\chi_{8\Delta}]\mathcal{T}_{8\alpha,2}(F)\mathcal{N}_{8\alpha}(v) d\sigma\lesssim \gamma^2\int_{8\Delta} \mathcal{M}[\sq{k}\chi_{8\Delta}] d\sigma\lesssim \gamma^2\sq{\omega}(\Delta).
\end{align*}

Lastly,

\begin{align*}
    J_1&\lesssim  \frac{1}{\lambda^2}\int_{8\Delta\setminus D_{\gamma\lambda}}\sum_{\substack{W\in \mathcal{D},\\ W\cap \Gamma^{4l}_{4\alpha}(P)\neq \emptyset}}\frac{\sup_{y\in W}\sq{G}(y)}{l(W)}l(W)\Big(\fint_{W}|\delta H|^2dy\Big)^{1/2}\mathcal{N}_{8\alpha}(v)d\sigma(P)
    \\
    &\lesssim \frac{1}{\lambda^2}\int_{8\Delta\setminus D_{\gamma\lambda}} \mathcal{M}[\sq{k}\chi_{8\Delta}]\mathcal{T}_{8\alpha,2}(\delta H)\mathcal{N}_{8\alpha}(v) d\sigma\lesssim \gamma^2\int_{8\Delta} \mathcal{M}[\sq{k}\chi_{8\Delta}] d\sigma.
\end{align*}

\textbf{Bound for $T_3$:} For this term, we can insert the mean value of the antisymmetric matrix by \refprop{prop:ContantAntiSymmMatrix} and obtain
\begin{align*}
\frac{1}{\lambda^2}\left|\int_\Omega T_3 \, dy \right| &= \Big| \frac{1}{2\lambda^2}\sum_{W\in \sq{\mathcal{D}}_{1,4}^{(l,\alpha)}}\int_{\Omega}(A^a-(A^a)_{\frac{3}{2}W})\nabla (v^2\sq{G}\phi_W) \cdot \nabla \psi\Big|
\\
&\lesssim\frac{1}{\lambda^2}\sum_{W\in \sq{\mathcal{D}}_{1,4}^{(l,\alpha)}}\int_{\Omega}|A^a-((A)^a)_{\frac{3}{2}W}|\left[|v\nabla v|\sq{G}\phi_W + v^2|\nabla \sq{G}|\phi_W + v^2\sq{G}|\nabla\phi_W|\right] |\nabla \psi|
\\
&\lesssim\frac{1}{\lambda^2}\sum_{W\in \sq{\mathcal{D}}_{1,4}^{(l,\alpha)}}l(W)^{n}\Big(\fint_{\frac{3}{2}W}|A^a-(A^a)_{\frac{3}{2}W}|^2\Big)^{1/2}
\\
&\qquad\cdot \left[\Big(\fint_{\frac{3}{2}W}|\nabla v|^2\Big)^{1/2}\sup_{\frac{3}{2}W}|v|\sup_{\frac{3}{2}W}\sq{G} + \sup_{\frac{3}{2}W} v^2\Big(\fint_{\frac{3}{2}W}|\nabla \sq{G}|^2\Big)^{1/2} + \sup_{\frac{3}{2}W}v^2\frac{\sup_{\frac{3}{2}W}\sq{G}}{\delta} \right]
\end{align*}

We can apply Caccioppoli and Harnack on the Green's function expressions and bound all remaining expressions by either the nontangential maximal function of \(v\) or the mean-valued nontangential maximal function of \(\delta|\nabla v|\) which are both controlled by \(\gamma\lambda\) in the good set \(E\). Thus, we obtain by \reflemma{lemma:IntOver1/delta}

\begin{align*}
    \frac{1}{\lambda^2}\left| \int_\Omega T_3 \, dy\right|&\lesssim \gamma^2\sum_{W\in \sq{\mathcal{D}}_{1,4}^{(l,\alpha)}} \int_{2W}\frac{\sq{G}}{\delta^2}\lesssim \gamma^2\Big(\int_{\sq{\Omega}_5^{(l,\alpha)}\setminus\sq{\Omega}_0^{(l,\alpha)}}\frac{\sq{G}}{\delta^2}\,dy\Big)
    \\
    &\lesssim \gamma^2\int_{8\Delta}\mathcal{M}[\sq{k}\chi_{8\Delta}]d\sigma.
\end{align*}
We postpone bounding this term for now.

\textbf{Bound for $T_4$:} For this term, we use the fact that $A^s$ is bounded to obtain
\[\frac{1}{\lambda^2}\left|\int_\Omega T_4 \, dy\right| \lesssim \frac{1}{\lambda^2}\int_\Omega |\nabla \psi||v|^2|\nabla \sq{G}| + |\sq{G}||v||\nabla v||\nabla\psi| \, dy.\]
To handle the first term on the right, we decompose the region \(\supp(\nabla\psi)\) into the Whitney cubes $W$ and using that \(|v|\leq \sup_E\mathcal{N}(v)\leq \lambda \gamma\) on \(\supp(\nabla \psi)\) we obtain
\[\frac{1}{\lambda^2}\int_{\Omega}|\nabla \psi||v|^2|\nabla \sq{G}| \, dy \lesssim \gamma^2 \int_{\supp(\nabla\psi)} \frac{|\nabla \sq{G}|}{\delta} \, dy  \approx \gamma^2 \sum_{W\in \sq{\mathcal{D}}_{1,4}^{(l,\alpha)}}  \int_{W}\frac{|\nabla \sq{G}|}{\delta} \, dy\]
and
\begin{align*}
    \frac{1}{\lambda^2}\int_{\Omega}|\sq{G}||v||\nabla v||\nabla\psi| \, dy & \lesssim \frac{\gamma}{\lambda} \sum_{W\in \sq{\mathcal{D}}_{1,4}^{(l,\alpha)}} \int_{W}\frac{|\sq{G}||\nabla v|}{\delta} \, dy \\
    & \lesssim \frac{\gamma}{\lambda} \sum_{W\in \sq{\mathcal{D}}_{1,4}^{(l,\alpha)}} |W| \left(\fint_{W}\frac{|\sq{G}|^2}{\delta^4} \, dy\right)^{1/2}\left(\fint_{W}\delta^2|\nabla v|^2 \, dy\right)^{1/2}
\end{align*}
We can apply Caccioppoli and Harnack on the Green's function expressions and bound the \(L^2\) Whitney cube average of \(\delta|\nabla v|\) by mean-valued nontangential maximal function of \(\delta|\nabla v|\) which is also controlled by \(\gamma\lambda\) in the good set \(E\).

This gives the bound
\[\frac{1}{\lambda^2}\left|\int_\Omega T_4dy\right| \lesssim \gamma^2 \sum_{W\in \sq{\mathcal{D}}_{1,4}^{(l,\alpha)}} \int_{\frac{3}{2}W}\frac{|\sq{G}|}{\delta^2} \, dy\lesssim \gamma^2\Big(\int_{\sq{\Omega}_5^{(l,\alpha)}\setminus\sq{\Omega}_0^{(l,\alpha)}}\frac{\sq{G}}{\delta^2}\,dy\Big)\lesssim \gamma^2\int_{8\Delta} \mathcal{M}[\sq{k}\chi_{8\Delta}] d\sigma,\]
where we used \reflemma{lemma:IntOver1/delta} for the last inequality.

All together, we obtain the bound
\[\frac{1}{\lambda^2}\int_\Omega A\nabla v\cdot \nabla v \sq{G}\psi \,dy=\frac{1}{\lambda^2}\int_\Omega T_1+T_2+T_3+T_4\,dy\lesssim \gamma^2\int_{8\Delta} \mathcal{M}[\sq{k}\chi_{8\Delta}] d\sigma.\]

Since \(\sq{\omega}\in B_p(\sigma)\) for some \(p>1\), we have by boundedness of the maximal function and doubling property of $\sq{\omega}$ that
\begin{align}
    \gamma^2\int_{8\Delta} \mathcal{M}[\sq{k}\chi_{8\Delta}] d\sigma&\lesssim \gamma^2\sigma(\Delta)\Big(\fint_{8\Delta} \mathcal{M}[\sq{k}\chi_{8\Delta}]^p d\sigma\Big)^{1/p}\nonumber
    \\
    &\lesssim \gamma^2\sigma(\Delta)\Big(\fint_{8\Delta} \sq{k}^p d\sigma\Big)^{1/p}\nonumber
    \\
    &\lesssim \gamma^2\sigma(\Delta)\fint_{8\Delta} \sq{k} d\sigma = \gamma^2\sq{\omega}(8\Delta)\lesssim \gamma^2\sq{\omega}(\Delta),\label{eq:BoundOFM[k]byomega}
\end{align}
which finishes the proof of \eqref{eq:ClaimInGoodLambdaInequality}.

\end{proof}

\subsection{Proof of \refthm{thm:MainThm1}}\label{sec:ProofMainThm1}

By assumption, we have that $(\operatorname{D})_{p'}^{L^*}$ is solvable, which implies $\omega_{L^*}\in A_\infty(\sigma)$. By \reflemma{lemma:SquareFunctionEstimate} and \eqref{eq:CarlAreaTcomparable} we obtain the square function bound \eqref{eq:SquareFunctionEstimate}.

In our setting of matrices with BMO antisymmetric part, the proofs of \eqref{eq:PD2*inequalityElliptic} and \eqref{eq:PDinequalityElliptic} are verbatim the same as in \cite[Theorem 1.11]{mourgoglou_solvability_2025} for matrices with bounded antisymmetric part. Apart from real variable arguments, their proof uses properties of the Green's function from \refthm{thm:green-function-construction}, \reflemma{lem:green-elliptic-measure-comparison}, Caccioppoli inequality, and the boundary Hölder inequality for solutions with vanishing Dirichlet data, all of which tools have  been established for matrices with antisymmetric part in BMO. 

\section{Solvability of the regularity problem for \texorpdfstring{$t$}{t}-independent operators (Proof of Theorem~\ref{thm:MainThm2})}\label{sec:ProofRegProblem}

We will work on the upper half space \(\Omega=\mathbb{R}^{n+1}_+=\mathbb{R}^n\times (0,\infty)\), where we call the last component the \(t\)-direction, or \(t\)-component. The \(t-\)independence assumption asserts \(A(x,t)=A(x)\), and we set the different components of the matrix \(A\) as
    \[A(x)=\begin{pmatrix} \AP(x) & b(x)\\c(x) & d(x)\end{pmatrix}=\begin{pmatrix} \AP^s(x) & b^s(x)\\c^s(x) & d^s(x)\end{pmatrix} + \begin{pmatrix} \AP^a(x) & b^a(x)\\c^a(x) & 0\end{pmatrix},\]
    where \(\AP(x)=\AP^s(x)+\AP^a(x)\in\mathbb{R}^{n\times n}, c(x)=c^s(x)+c^a(x)\in\mathbb{R}^{1\times n}, b(x)=b^s(x)+b^a(x)\in\mathbb{R}^{ n\times 1}, d(x)=d^s(x)\in\mathbb{R}\). The matrix consisting of \(\AP^s,b^s,c^s,d^s\) is the symmetric part of \(A\), while the matrix consisting of \(\AP^a,b^a,c^a\) is the antisymmetric part of \(A\). We assume that the symmetric components are bounded above, whereas the antisymmetric ones lie in $BMO(\Omega)$. Note that since $A$ is $t$-independent, $A^a\in BMO(\mathbb{R}^n\times(0,\infty))$ implies $A^a\in BMO(\mathbb{R}^n)$. 
    We define the elliptic operators \(L_{||}=\mathrm{div}_{||}(\AP\nabla_{||} \cdot)\), where the subscript \(||\) means that we take the gradient or divergence only with respect to the first \(n\) components or the \(x\)-components only.

    As in \cite{hofmann_regularity_2015}, we will make use of the approximation operator with elliptic homogeneity:  
    \[\mathcal{P}_tf:=e^{-t^2L_{||}}f.\]
    In \cite{hofmann_lp_2022}, a related operator, the vertical area function, has already been studied for matrices with BMO antisymmetric part. For the conical area function, the following bounds are well known (see for instance $(4.3)$ in \cite{hofmann_regularity_2015}), if the antisymmetric part of the matrix is bounded. Similar arguments also apply when the antisymmetric part belongs to BMO, and the proof of the following proposition can be found in Appendix~\ref{sec:AreaFctBounds}.
    \begin{prop}\label{prop:AreaFctBoundsInLp}
    For all \(1<p<\infty, f\in \dot{W}^{1,p}(\mathbb{R}^n)\), it holds that
    \begin{enumerate}
        \item  \(\Vert \mathcal{A}(\partial_t\mathcal{P}_tf)\Vert_{L^p(\mathbb{R}^n)}\leq C \Vert \nabla_{||}f\Vert_{L^p(\mathbb{R}^n)},\)
        \item  \(\Vert \mathcal{A}(t\partial_t^2\mathcal{P}_tf)\Vert_{L^p(\mathbb{R}^n)}\leq C \Vert \nabla_{||}f\Vert_{L^p(\mathbb{R}^n)},\)
        \item \(\Vert \mathcal{A}(t\nabla_{||}\partial_t\mathcal{P}_tf)\Vert_{L^p(\mathbb{R}^n)}\leq C \Vert \nabla_{||}f\Vert_{L^p(\mathbb{R}^n)}\).
    \end{enumerate}
\end{prop}

\begin{rem}[Nontangential maximal function bounds]\label{rem:NontangnetialBoundsForSemigroup}
    Since $\sq{\mathcal{N}}_{\alpha,2}(F)\lesssim \mathcal{A}_{2\alpha}(F)$ holds pointwise for any function $F$ and $\alpha>1$, we immediately obtain bounds for the nontangential maximal function that are analogous to the area function bounds in
    \refprop{prop:AreaFctBoundsInLp}, i.e. for all $1<p<\infty$ and $f\in W^{1,p}(\mathbb{R}^n)$ 
    \begin{align*} \Vert\sq{\mathcal{N}}_2(\partial_t\mathcal{P}_tf)\Vert_{L^p} + \Vert\sq{\mathcal{N}}_2(t\partial_t^2\mathcal{P}_tf)\Vert_{L^p}+ \Vert\sq{\mathcal{N}}_2(t\nabla_{||}\partial_t\mathcal{P}_tf)\Vert_{L^p} + \Vert\sq{\mathcal{N}}_2(t^2\nabla_{||}\partial^2_t\mathcal{P}_tf)\Vert_{L^p}
    \\
    \lesssim \Vert\nabla_{||}f\Vert_{L^p}.\end{align*} 
    
    Furthermore, we also have the bounds
    \begin{align}\Vert\sq{\mathcal{N}}_2(\mathcal{P}_tf)\Vert_{L^p}\lesssim\Vert\mathcal{N}(\mathcal{P}_tf)\Vert_{L^p}\lesssim\Vert f\Vert_{L^p}\qquad \textrm{ for all }1<p<\infty,\label{eq:nontangnetialofPBoundedbyf}\end{align}
    and \begin{align}\Vert\sq{\mathcal{N}}_2(\nabla_{||}\mathcal{P}_tf)\Vert_{L^p}\lesssim\Vert\nabla_{||}f\Vert_{L^p}\qquad \textrm{ for all }1<p<\infty.\label{eq:nontangnetialofNablaPBoundedbyNablaf}\end{align}
    These results essentially follow from the Gaussian kernel bounds of the semigroup and we provide a proof of \eqref{eq:nontangnetialofPBoundedbyf} and \eqref{eq:nontangnetialofNablaPBoundedbyNablaf} in Appendix~\ref{sec:AppendixNontangnetialBound}.
\end{rem}

In addition, we have \(L^2-L^2\) boundedness of the following operators.

\begin{prop}[\cite{hofmann_lp_2022}]\label{prop:L2boundednessOfSemigroupOperators}
    Let \(f\in L^2(\mathbb{R}^n)\). Then for \(S_t\in\{e^{-tL_{||}}, t\partial_te^{-tL_{||}}, \sqrt{t}\nabla_{||} e^{-tL_{||}}\}\) we have
    \begin{align*}
        \Vert S_tf\Vert_{L^2}\leq C\Vert f\Vert_{L^2},
    \end{align*}
    where the implicit constant \(C\) only depends on the ellipticity constant of \(\AP\). Furthermore, we also have for \(f\in W^{1,2}(\mathbb{R}^n)\)
    \[\Vert \nabla_{||} e^{-tL_{||}}f\Vert_{L^2}+\Vert t\nabla_{||} L_{||}e^{-tL_{||}}f\Vert_{L^2}\leq C\Vert \nabla_{||} f\Vert_{L^2}.\]
    This statement immediately follows from the Kato estimate.
\end{prop}
    
    Furthermore, we have off-diagonal estimates:
    
\begin{prop}[Off-diagonal estimates, {\cite[Prop. 4.9]{hofmann_lp_2022}}]\label{prop:Off-diagonalEstimates}
    For \(S_t\in \{e^{-tL_{||}},\sqrt{t}\nabla e^{-tL_{||}},t\partial_te^{-tL_{||}}\}\) there exists \(C,\alpha>0\) such that
    \[\Vert S_th\Vert_{L^2(E)}\lesssim e^{-\frac{\alpha d(E,F)^2}{t}}\Vert h\Vert_{L^2(\mathbb{R}^n)}\] 
    for all \(h\in L^2(\mathbb{R}^n)\) with \(\mathrm{supp}(h)\subset F\).
\end{prop}

A corollary from the off-diagonal estimates in \refprop{prop:Off-diagonalEstimates} is the off-diagonal estimate for the operator \(\sqrt{t}e^{-tL}\mathrm{div}\)
\begin{cor}\label{cor:Off-diagonalEstimatesForPdiv}
    There exists \(C,\alpha>0\) such that
    \[\Vert \sqrt{t}e^{-tL}\mathrm{div} (h)\Vert_{L^2(E)}\lesssim e^{-\frac{\alpha d(E,F)^2}{t}}\Vert h\Vert_{L^2(\mathbb{R}^n)}\] 
    for all \(h\in (L^2(\mathbb{R}^n))^{n}\) with \(\mathrm{supp}(h)\subset F\).
\end{cor}
Although this result follows from using the adjoint off-diagonal estimates for \(\sqrt{t}\nabla e^{-tL}\) from \refprop{prop:Off-diagonalEstimates}, a short direct proof
is provided in the appendix~\ref{sec:AppendixOffDiago}.

We will require the important extension of the Kato estimates for operators with BMO-antisymmetric matrices. The Kato conjecture for bounded matrices was solved in \cite{Auscher_Kato}.
\begin{prop}[Kato, \cite{escauriaza_kato_2018} for $L^2$, \cite{hofmann_lp_2022} for $L^p$]\label{prop:Kato}
    For every $1<p<\infty$ and $f\in \dot{W}^{1,p}(\mathbb{R}^n)$
    \[\Vert L^{1/2} f\Vert_{L^p(\mathbb{R}^n)}\lesssim \Vert \nabla f\Vert_{L^p(\mathbb{R}^n)}.\]
    Furthermore, there exists $\EPS>0$, such that for every $1<p<2+\EPS$ and $f\in \dot{W}^{1,p}(\mathbb{R}^n)$
    \[ \Vert \nabla f\Vert_{L^p(\mathbb{R}^n)}\lesssim \Vert L^{1/2} f\Vert_{L^p(\mathbb{R}^n)}.\]
    
\end{prop}

\begin{lemma}\label{lemma:LocalDivCurl}
    Let \(1\leq i,j\leq n+1\) and \(W\subset \mathbb{R}^{n+1}_+\). Furthermore, assume that \(b\in BMO(\Omega)\). Then for every \(\EPS>0\) there exists \(C=C(\EPS)>0\) such that for all \(g\in W_0^{1,2}(W)\) and \(h\in W^{1,2+\EPS}(W)\).  
    \[\Big|\fint_W b(x)(\partial_ig(x)\partial_jh(x)-\partial_jg(x)\partial_i h(x))\,dx\Big|\leq C \Vert b\Vert_{BMO(\Omega)}\Big(\fint_W |\nabla g|^2dx\Big)^{\frac{1}{2}}\Big(\fint_W|\nabla h|^{2+\EPS}dx\Big)^{\frac{1}{2+\EPS}}\]
    and for all \(g\in W_0^{1,2+\EPS}(W)\) and \(h\in W^{1,2}(W)\)
    \[\Big|\fint_W b(x)(\partial_ig(x)\partial_jh(x)-\partial_jg(x)\partial_i h(x))\,dx\Big|\leq C \Vert b\Vert_{BMO(\Omega)}\Big(\fint_W |\nabla g|^{2+\EPS}dx\Big)^{\frac{1}{2+\EPS}}\Big(\fint_W|\nabla h|^{2}dx\Big)^{\frac{1}{2}}\]
\end{lemma}

\begin{proof}
First, assume that \(g\in C_c^\infty(W)\) and \(h\in W^{1,2+\EPS}(\Omega)\). We observe that
\begin{align}\int_W\partial_ig\partial_jh-\partial_j g\partial_i h dx&=\int_W \big(\partial_i\partial_jg -\partial_j\partial_ig\big)hdx=0\label{eq:ObMV=0}
\end{align}
Hence, the mean value of \(\partial_ig\partial_jh-\partial_j g\partial_i h\) is zero. Denoting the mean value as \((b)_{\frac{3}{2}W}=\fint_{\frac{3}{2}W}b(x) \,dx\), we obtain
\begin{align*}
&\Big|\fint_W b(x)(\partial_ig(x)\partial_jh(x)-\partial_jg(x)\partial_i h(x))\,dx\Big|
\\
&\Big|\fint_W (b(x)-(b)_{\frac{3}{2}W})(\partial_ig(x)\partial_jh(x)-\partial_jg(x)\partial_i h(x))\,dx\Big|
\\
&\leq \Big(\fint_W |b-(b)_{\frac{3}{2}W}|^{\frac{4+2\EPS}{\EPS}}dx\Big)^{\frac{\EPS}{4+2\EPS}}\Big(\fint_W |\nabla g|^2dx\Big)^{\frac{1}{2}}\Big(\fint_W|\nabla h|^{2+\EPS}dx\Big)^{\frac{1}{2+\EPS}}
\\
&\leq C \Vert b\Vert_{BMO(\Omega)}\Big(\fint_W |\nabla g|^2dx\Big)^{\frac{1}{2}}\Big(\fint_W|\nabla h|^{2+\EPS}dx\Big)^{\frac{1}{2+\EPS}}.
\end{align*}
By density of \(C_c^\infty\) in the Sobolev space \(W_0^{1,2}(\Omega)\), we obtain the first assertion of the lemma.
The argument for the second assertion works analogously.

\end{proof}

\subsection{Proof of \refthm{thm:MainThm2}}

\subsubsection{Overview over the Proof structure}

Since solvability of \((\operatorname{D})_{p'}^{L^*}\) for \(t-\)independent operators was already proven in \cite{hofmann_dirichlet_2022}, it is enough to prove that solvability of \((\operatorname{D})_{p'}^{L^*}\) implies solvability of \((\operatorname{R})_p^{L^*}\). 
The following proof of \((\operatorname{D})_{p'}^{L^*}\Rightarrow (\operatorname{R})_p^L\) contains a lot of technical steps, and we split it into three parts. First we present the main structure and idea of the proof. As we will point out, we need to justify the convergence of several integrals. These justifications will be postponed to the end as to not divert attention away from the main structure and arguments of the proof. Furthermore, several error terms arise, which we will call $E_1^R,...,E_4^R$. All of these will tend to zero, and as before we postpone dealing with them until the end, so that the main structure of the proof remains clear.

Furthermore, the proof makes use of two lemmas whose proofs will be provided in a later subsection.

\subsection{Main structure of the proof}

To begin with, let us recall that we assume \(p>1\) to be such that \((\operatorname{D})^{L^*}_{p'}\). Let \(g\in \dot{W}^{1,p}(\partial\Omega)\cap C^\infty_c(\partial\Omega)\) and let \(u\in Y^{1,2}(\Omega)\) be the solution to \(Lu=0\) with boundary data \(f=f_\EPS:=\mathcal{P}_\EPS g\in \dot{W}^{1,p}(\partial\Omega)\). Since \(f_\EPS\to g\) in \(\dot{W}^{1,p}(\partial\Omega)\) as \(\EPS\to 0\), it suffices to prove 
\[\Vert\sq{\mathcal{N}}_2(\nabla u)\Vert_{L^p(\partial\Omega)}\leq C \Vert \nabla_{||}f_\EPS\Vert_{L^p(\partial\Omega)}=C \Vert \nabla_{||}f\Vert_{L^p(\partial\Omega)}\] with constants independent of \(\EPS\) or \(g\). Let us fix \(1<q<\min\{p,(n+1)/n\}\). By entire real variable arguments (as for instance in \cite{shen_relationship_2007}), we have \(\Big(\fint_{B(y,s,s/2)}|\nabla u|^2 dzd\tau\Big)^{1/2}\lesssim \Big(\fint_{B(y,s,3s/4)}|\nabla u|^{q} dzd\tau\Big)^{1/q}\) for every point \((y,s)\in\Omega\).

Hence, by Carleson's theorem \refprop{prop:NontangnetialAndCarelsonAreDual}, we have 
\[\Vert\sq{\mathcal{N}}_2(\nabla u)\Vert_{L^p}\lesssim \Vert\sq{\mathcal{N}}_q(\nabla u)\Vert_{L^p}\lesssim \int_\Omega \nabla u \cdot h dxdt\]
for a dualising function \(h\in L^{\infty}_c(\Omega)\) with \(\Vert \sq{\mathcal{C}}_{q'}(h) \Vert_{L^{p'}(\Omega)}=1\).

Now, let \(v\in W^{1,2}(\Omega)\) be the solution to the Poisson-Dirichlet problem of the adjoint operator
\begin{align}\begin{cases} L^*v=\mathrm{div}(h) &\textrm{in }\Omega,\\ v=0 &\textrm{on }\partial\Omega.\end{cases}\label{PDEForv}\end{align}
By \refthm{thm:MainThm1} and \eqref{eq:CarlAreaTcomparable}, we have the bounds 
\begin{align}
    \Vert \mathcal{N}(v)\Vert_{L^{p'}}+\Vert \mathcal{S}(v)\Vert_{L^{p'}}\lesssim \Vert \sq{\mathcal{C}}_{q'}(h)\Vert_{L^{p'}}=1, \quad \textrm{and}\quad \Vert \mathcal{T}_{q'}(h)\Vert_{L^{p'}}\approx \Vert \sq{\mathcal{C}}_{q'}(h)\Vert_{L^{p'}}=1.\label{eq:EstimatesFromPDforR}
\end{align}

Let us now fix a large $R$, and set $\eta^R\in C^\infty((0,\infty))$ with $0\leq\eta^R\leq 1$, $\eta^R\equiv 1$ on $(0,R/2)$, $\eta^R\equiv 0$ on $(R,\infty)$, and $|\partial_s\eta^R|R+|\partial_s^2\eta^R|R^2\lesssim 1$, and define
\[ \mathcal{P}_sf(x):=e^{-s^2L_{||}}f(x),\quad\mathcal{P}_s^Rf(x):=\mathcal{P}_sf(x)\eta^R(s),\quad\textrm{and}\qquad \Omega^R:=\mathbb{R}^n\times(0,R).\]

Without loss of generality, we can assume that $\mathrm{supp}(f)\subset\Delta(0,R/4)$ and $\mathrm{supp}(h)\subset B(0,R/4)$. Then, using integration by parts, \(u-\mathcal{P}_t^Rf\in W^{1,2}_0(\Omega)\) and the PDEs for \(v\) and \(u\) yield
\begin{align}
    \int_\Omega \nabla u \cdot h dxdt&=\int_\Omega \nabla (u-\mathcal{P}_t^Rf) \cdot h dxdt+\int_\Omega \nabla \mathcal{P}_t^Rf \cdot h dxdt\nonumber
    \\
    &=\int_\Omega A\nabla u\cdot \nabla v dxdt - \int_\Omega A\nabla \mathcal{P}_t^Rf\cdot \nabla v dxdt+\int_\Omega \nabla \mathcal{P}_t^Rf \cdot h dxdt\nonumber
    \\
    &=- \int_\Omega A\nabla \mathcal{P}_t^Rf\cdot \nabla v dxdt+\int_\Omega \nabla \mathcal{P}_t^Rf \cdot h dxdt\label{eq:FirstTermInDImpliesR}
\end{align}

By Carleson's Theorem \refprop{prop:NontangnetialAndCarelsonAreDual}, the last term is bounded by
\begin{align*}
    \Big|\int_\Omega \nabla \mathcal{P}_t^Rf \cdot h dxdt\Big|\lesssim \Vert \sq{\mathcal{N}}(\nabla_{||} \mathcal{P}^R_tf)\Vert_{L^p}+\Vert \sq{\mathcal{N}}(\partial_t \mathcal{P}^R_tf)\Vert_{L^p}\lesssim \Vert \nabla_{||} f\Vert_{L^p}, 
\end{align*}
where the last inequality follows from \refrem{rem:NontangnetialBoundsForSemigroup}.

Thus, it remains to bound the first term in \eqref{eq:FirstTermInDImpliesR}. We have
\begin{align*}
     &\int_\Omega A(x)\nabla \mathcal{P}_t^Rf(x)\cdot \nabla v(x)dx
     \\
     &=\int_{\Omega^R}\AP(x)\nabla_{||}\mathcal{P}^R_tf(x)\cdot\nabla_{||}v(x,t) dxdt + \int_{\Omega^R}d(x)\partial_{t}\mathcal{P}^R_tf(x)\partial_t v(x,t) dxdt
     \\
    &\qquad\int_{\Omega^R}c(x)\cdot\nabla_{||}\mathcal{P}^R_tf(x)\partial_t v(x,t)+\partial_{t}\mathcal{P}^R_tf(x)b(x)\cdot\nabla_{||} v(x,t) dxdt 
    \\
    &=:I_1+I_2+I_3.
\end{align*}

First, for \(I_2\), we can observe that $d=d^s$ is bounded, and hence Hölder's inequality, \eqref{eq:EstimatesFromPDforR}, and \refprop{prop:AreaFctBoundsInLp} give
\begin{align}
I_2\lesssim \Vert\mathcal{A}(\partial_t\mathcal{P}_tf)\Vert_{L^p(\partial\Omega)} \Vert\mathcal{S}(v)\Vert_{L^{p'}(\partial\Omega)}
\lesssim\Vert\nabla_{||}f\Vert_{L^p(\partial\Omega)}.\label{eq:BoundI2}
\end{align}

Next, for \(I_1\), we use integration by parts and $\partial_t t=1$ to get
\begin{align*}
    I_1&=\int_\Omega A_{||}\nabla_{||}\mathcal{P}^R_tf(x) \cdot\nabla_{||} v(x,t) \partial_tt\,dxdt
    \\
    &=-\int_\Omega A_{||}\nabla_{||}\partial_t\mathcal{P}^R_tf(x) \cdot\nabla_{||} v(x,t) t\,dxdt - \int_\Omega A_{||}\nabla_{||}\mathcal{P}_tf(x)\eta^R(t) \cdot\nabla_{||} \partial_tv(x,t) t\,dxdt
    \\
    &=:I_{11}+I_{12}.
\end{align*}
For the second term \(I_{12}\), again, we integrate by parts and get, by \eqref{eq:EstimatesFromPDforR}, and \refprop{prop:AreaFctBoundsInLp},
\begin{align*}I_{12}&\lesssim \int_\Omega L_{||}\mathcal{P}_tf(x)\eta^R(t)\partial_tv(x,t) t\,dxdt
\lesssim \Vert\mathcal{A}(\partial_t\mathcal{P}_tf)\Vert_{L^p(\partial\Omega)} \Vert\mathcal{S}(v)\Vert_{L^{p'}(\partial\Omega)}
\lesssim\Vert\nabla_{||}f\Vert_{L^p(\partial\Omega)}.\end{align*}

For the first term, we split the matrix into symmetric and antisymmetric part
\[I_{11}=-\int_\Omega A^s_{||}\nabla_{||}\partial_t\mathcal{P}^R_tf(x) \cdot\nabla_{||} v(x,t) t\,dxdt-\int_\Omega A_{||}^a\nabla_{||}\partial_t\mathcal{P}^R_tf(x) \cdot\nabla_{||} v(x,t) t\,dxdt=:I_{11}^s+I_{11}^a.\]
For the symmetric part, Hölder inequality, \eqref{eq:EstimatesFromPDforR}, and \refprop{prop:AreaFctBoundsInLp} gives

\begin{align*}I_{11}^s&\lesssim \Vert\mathcal{A}(t\nabla_{||}\partial_t\mathcal{P}^R_tf)\Vert_{L^p(\partial\Omega)} \Vert\mathcal{S}(v)\Vert_{L^{p'}(\partial\Omega)}
\\
&\leq \Vert\mathcal{A}(t\nabla_{||}\partial_t\mathcal{P}_tf)\Vert_{L^p(\partial\Omega)} \Vert\mathcal{S}(v)\Vert_{L^{p'}(\partial\Omega)}+\Vert\mathcal{A}(t\partial_t\eta^R\nabla_{||}\mathcal{P}_tf)\Vert_{L^p(\partial\Omega)} \Vert\mathcal{S}(v)\Vert_{L^{p'}(\partial\Omega)}
\\
&\lesssim\Vert\nabla_{||}f\Vert_{L^p(\partial\Omega)} + E_1^R,\end{align*}
if we set \(E_1^R:=\Vert\mathcal{A}(t\partial_t\eta^R\nabla_{||}\mathcal{P}_tf)\Vert_{L^p(\partial\Omega)}\). The term \(E_1^R\) is an error term and we postpone its bound until the end. For the other term \(I_{11}^a\) we use the Whitney decomposition and partition of unity introduced in Section~\ref{Section:Defphi_W} and \reflemma{lemma:LocalDivCurl} to get
\begin{align*}
    I_{11}^a&=-\sum_{W\in\mathcal{D}(\Omega)}|W|\fint_{\frac{3}{2}W} A^a_{||}\nabla_{||}\big(\phi_W\partial_t\mathcal{P}^R_tf\big)(x) \cdot\nabla_{||} v(x,t) t\,dxdt
    \\
    &\lesssim \Vert \AP^a\Vert_{BMO(\Omega)}\sum_{W\in\mathcal{D}(\Omega)}l(W)|W|\Big(\fint_{\frac{3}{2}W}\frac{1}{l(W)}|\partial_t\mathcal{P}_t^Rf|^2dxdt\Big)^\frac{1}{2}\Big(\fint_{\frac{3}{2}W}|\nabla v|^{2+\EPS}dxdt\Big)^\frac{1}{2+\EPS}
    \\
    &\quad + \Vert \AP^a\Vert_{BMO(\Omega)}\sum_{W\in\mathcal{D}(\Omega)}l(W)|W|\Big(\fint_{\frac{3}{2}W}|\nabla_{||}\partial_t\mathcal{P}_t^Rf|^2dxdt\Big)^\frac{1}{2}\Big(\fint_{\frac{3}{2}W}|\nabla v|^{2+\EPS}dxdt\Big)^\frac{1}{2+\EPS}.
\end{align*}
Here, we chose $\EPS$ sufficiently small so that \reflemma{lemma:HigherIntegrabilityForGradient} can be applied, and combined with \reflemma{lemma:dyadicToNondyadic} this yields
\begin{align*}
    I_{11}^a&\lesssim \Vert \AP^a\Vert_{BMO(\Omega)}\sum_{W\in\mathcal{D}(\Omega)}l(W)|W|\Bigg[\Big(\fint_{\frac{3}{2}W}\frac{1}{l(W)}|\partial_t\mathcal{P}_t^Rf|^2dxdt\Big)^\frac{1}{2} + \Big(\fint_{\frac{3}{2}W}|\nabla_{||}\partial_t\mathcal{P}_t^Rf|^2dxdt\Big)^\frac{1}{2}\Bigg]
    \\
    &\hspace{50mm}\cdot\Bigg[\Big(\fint_{\frac{3}{2}W}|\nabla v|^{2}dxdt\Big)^\frac{1}{2}+\Big(\fint_{\frac{3}{2}W}|h|^{2+\EPS}dxdt\Big)^\frac{1}{2+\EPS}\Bigg]
    \\
    &\lesssim \big(\Vert\mathcal{A}(\partial_t\mathcal{P}^R_tf)\Vert_{L^p(\partial\Omega)}+ \Vert\mathcal{A}(t\partial_t\nabla_{||}\mathcal{P}^R_tf)\Vert_{L^p(\partial\Omega)}\big)\big(\Vert\mathcal{S}(v)\Vert_{L^{p'}(\partial\Omega)}+\Vert\mathcal{T}_{2+\EPS}(v)\Vert_{L^{p'}(\partial\Omega)}\big).
\end{align*}
Setting a second error term as \(E_2^R:=\Vert\mathcal{A}(\partial_t\eta^R\mathcal{P}_tf)\Vert_{L^p(\partial\Omega)}\), we obtain by \eqref{eq:EstimatesFromPDforR} and \refprop{prop:AreaFctBoundsInLp}
\[I_{11}^a\lesssim \Vert\mathcal{A}(\partial_t\mathcal{P}_tf)\Vert_{L^p(\partial\Omega)}+\Vert\mathcal{A}(t\partial_t\nabla_{||}\mathcal{P}^R_tf)\Vert_{L^p(\partial\Omega)}+E_1^R+E_2^R\lesssim \Vert \nabla_{||}f\Vert_{L^p(\partial\Omega)}+E_1^R+E_2^R.\]
We postpone bounding the error terms until the end. This finishes the bound of \(I_1\).

Hence, it remains to bound \(I_3\). Before we continue, we would like to note that a justification of finiteness of \(I_3\) is required, and will be provided at the end of this proof. Now, assuming finiteness, for the antisymmetric components we have that \(b^a_j=-c_j^a\) for all \(j=1,...,n\) and for the symmetric ones \(b_j^s=c_j^s\). We can now split the integral \(I_3\) into sums over symmetric and antisymmetric components and get
\begin{align*}
    I_3&=\int_{\Omega^R} \partial_t\mathcal{P}_t^Rf(x) b^{a}(x)\cdot \nabla _{||}v(x,t)  -  \nabla_{||}\mathcal{P}_t^Rf(x)\cdot b^{a}(x) \partial_t v(x,t) dxdt
    \\
    &\qquad+\int_{\Omega^R} \partial_t\mathcal{P}_t^Rf(x) b^{s}(x)\cdot \nabla _{||}v(x,t) +  \nabla_{||}\mathcal{P}_t^Rf(x)\cdot b^{s}(x) \partial_t v(x,t) dxdt
    \\
    &=:J^a+J^s.
\end{align*}
Again, we need to justify the convergence of $J^a$ and $J^s$, and postpone the required arguments until the end. To bound $J^a$ more precisely, let \(\mathfrak{A}_t\) given by 
\begin{align}
\mathfrak{A}_tf(x):=\sum_{W\in \mathcal{D}(\Omega)}\phi_{W'}(x,t)\fint_{W_t}f(y,t)dy,\qquad\textrm{where }\quad W_s:=W\cap\{(x,s)\in\mathbb{R}^{n+1}_+|s=t\},\label{eq:DefSmoothAveragingOp}
\end{align}
be the averaging operator on slices on height \(t\). We can proceed by splitting \(J^a\) into
\begin{align*}
    J^a&=\sum_{1\leq j\leq n}\int_{\Omega^R} b^{a}_j \Big[\partial_{t}\mathcal{P}^R_tf\, \partial_j(I-\mathfrak{A}_t)v -  \partial_j\mathcal{P}^R_tf\,\partial_t (I-\mathfrak{A}_t)v\Big] dxdt
    \\
    &\qquad +\sum_{1\leq j\leq n}\int_{\Omega^R} b^{a}_j \Big[\partial_{t}\mathcal{P}^R_tf \,\partial_j\mathfrak{A}_tv -  \partial_j\mathcal{P}^R_tf\, \partial_t \mathfrak{A}_tv\Big] dxdt
    \\
    &=:J_1^a+J_2^a.
\end{align*}
Note that convergence of $J_1^a$ and $J_2^a$ will also be justified in the next subsection.
For \(J_2^a\) we integrate the second term by parts in $t$, to obtain
\begin{align*}
    J_2^a&=\sum_{1\leq j\leq n}\int_{\Omega^R} b^{a}_j \Big[\partial_{t}\mathcal{P}_t^Rf \,\partial_j\mathfrak{A}_tv +  \partial_j\partial_{t}\mathcal{P}_t^Rf\, \mathfrak{A}_tv\Big] dxdt
    \\
    &=\sum_{1\leq j\leq n}\int_{\Omega^R} b^{a}_j \partial_j\big(\partial_{t}\mathcal{P}^R_tf \mathfrak{A}_tv\big) dxdt.
\end{align*}

A similar term will also appear in the symmetric part \(J^s\), so we postpone bounding this term more precisely to later (see the bound of $\sq{J}$).

For \(J_1^a\) we insert the Whitney decomposition and partition of unity (see Section \ref{Section:Defphi_W}) and get

\begin{align*}
    J_1^a&=\sum_{\substack{1\leq j\leq n,\\ W\in\mathcal{D}(\Omega)}}\int_{\Omega^R} b^{a}_j \Big[\partial_{t}\mathcal{P}^R_tf\, \partial_j(\phi_W(I-\mathfrak{A}_t)v) -  \partial_j\mathcal{P}^R_tf\,\partial_t (\phi_W(I-\mathfrak{A}_t)v)\Big] dxdt
\end{align*}

Due to an observation \eqref{eq:ObMV=0} from the proof of \reflemma{lemma:LocalDivCurl} with \(g=\mathcal{P}_tf(x)\) and \(h=\phi_{W}(I-\mathfrak{A}_t)v\), we note that the mean value of \(\partial_tg\partial_jh-\partial_jg \partial_th\) integrated over \(\frac{3}{2}W\) is zero. Hence, \refprop{prop:ContantAntiSymmMatrix} applies, and we obtain
\begin{align*}
    J^a_1=\sum_{j, W\in \mathcal{D}}\int_{\frac{3}{2}W} (b^{a}_j-(b^{a}_j)_{\frac{3}{2}W}) \Big[\partial_{t}\mathcal{P}^R_tf\, \partial_j(\phi_W(I-\mathfrak{A}_t)v)
      -  \partial_j\mathcal{P}^R_tf\,\partial_t (\phi_W(I-\mathfrak{A}_t)v)\Big] dxdt.
\end{align*}
We simplify and use integration by parts in $s$ on the second summand to continue with
\begin{align*}
    &=\sum_{j,  W\in \mathcal{D}}\int_{\frac{3}{2}W} (b^{a}_j-(b^{a}_j)_{\frac{3}{2}W}) \Big[\partial_{t}\mathcal{P}^R_tf \,\partial_j(\phi_W(I-\mathfrak{A}_t)v)
      +  t\partial_j\partial_{t}\mathcal{P}_t^Rf\,\phi_W(I-\mathfrak{A}_t)v\Big] dxdt
\end{align*}
Now, splitting \(\partial_j(\phi_W(I-\mathfrak{A}_t)v)=\phi_W\partial_j(I-\mathfrak{A}_t)v+\partial_j\phi_W(I-\mathfrak{A}_t)v\) gives 
\begin{align*}
    &=\sum_{j,  W\in \mathcal{D}}\int_{\frac{3}{2}W} (b^{a}_j-(b^{a}_j)_{\frac{3}{2}W}) \Big[\partial_{t}\mathcal{P}_t^Rf\, \phi_{W} \partial_j(I-\mathfrak{A}_t)v + \partial_{t}\mathcal{P}_t^Rf\, \partial_j\phi_{W} (I-\mathfrak{A}_t)v
    \\
    &\hspace{8cm}+ \partial_j\partial_{t}\mathcal{P}_t^Rf\, \phi_{W} (I-\mathfrak{A}_t)v\Big] dxdt
    \\
    &\lesssim\sum_{j,  W\in \mathcal{D}}|W|\Big(\fint_{\frac{3}{2}W} |b^{a}_j-(b^{a}_j)_{\frac{3}{2}W}|^r \Big)^{\frac{1}{r}}\Big[\Big(\fint_{\frac{3}{2}W}|\partial_{t}\mathcal{P}_t^Rf|^2\Big)^\frac{1}{2}
    \Big(\fint_{\frac{3}{2}W}|\nabla_{||} (I-\mathfrak{A}_t)v|^{2+\EPS}\Big)^{\frac{1}{2+\EPS}}
    \\
    &\qquad+ \frac{2}{l(W)}\Big(\fint_{\frac{3}{2}W}|\partial_{t}\mathcal{P}_t^Rf|^2\Big)^\frac{1}{2}\Big(\fint_{\frac{3}{2}W}|(I-\mathfrak{A}_t)v|^{2+\EPS}\Big)^{\frac{1}{2+\EPS}}
    \\
    &\qquad + \Big(\fint_{\frac{3}{2}W}|\partial_j\partial_{t}\mathcal{P}^R_tf|^2\Big)^\frac{1}{2}\Big(\fint_{\frac{3}{2}W}|(I-\mathfrak{A}_t)v|^{2+\EPS}\Big)^{\frac{1}{2+\EPS}}\Big].
\end{align*}
In the last line we used Hölder's inequality with a large \(r\) chosen such that \(\frac{1}{r}+\frac{1}{2+\EPS}=\frac{1}{2}\), where \(\EPS\) is small enough such that the reverse H\"{o}lder inequality \reflemma{lemma:HigherIntegrabilityForGradient} can be used.
We would like to discuss the terms involving \(\mathfrak{A}_tv\) next. Since \(\sum_{W'\in\mathcal{D}}\nabla \phi_{W'}\equiv 0\), we obtain
\begin{align}
     &\Big(\fint_{\frac{3}{2}W}|\nabla \mathfrak{A}_tv|^{2+\EPS}\Big)^{\frac{1}{2+\EPS}}\leq  \Big(\fint_{\frac{3}{2}W}\big|\sum_{W'}\nabla \phi_{W'}(v)_{W'}\big|^{2+\EPS}\Big)^{\frac{1}{2+\EPS}}\label{eq:NablaI-mathfrakA}
     \\
     &\leq \Big(\fint_{\frac{3}{2}W}\big|\sum_{W'}\nabla \phi_{W'}(v)_{W'}-\big(\sum_{W'}\nabla \phi_{W'}\big) v\big|^{2+\EPS}\Big)^{\frac{1}{2+\EPS}}\nonumber
    \\
     &\leq \sum_{W'\in\mathcal{D}; \frac{3}{2}W'\cap\frac{3}{2}W\neq\emptyset}\Big(\fint_{\frac{3}{2}W}|\nabla \phi_{W'}((v)_{W'}- v)|^{2+\EPS}\Big)^{\frac{1}{2+\EPS}}\nonumber
     \\
     &\lesssim \sum_{W', \frac{3}{2}W'\cap\frac{3}{2}W\neq \emptyset}\frac{1}{l(W')}\Big(\fint_{\frac{3}{2}W'}|(v)_{W'}- v|^{2+\EPS}\Big)^{\frac{1}{2+\EPS}}.\nonumber
\end{align}
Since the overlap of the Whitney cubes is finite, and all of the Whitney cubes \(W'\) in the sum have length comparable to \(l(W)\), we can use Poincar\'{e}'s inequality and \reflemma{lemma:HigherIntegrabilityForGradient} for \(\nabla v\) to continue with

\begin{align}\Big(\fint_{\frac{3}{2}W}|\nabla \mathfrak{A}_tv|^{2+\EPS}\Big)^{\frac{1}{2+\EPS}}\lesssim\Big(\fint_{\frac{31}{16}W}|\nabla v|^{2+\EPS}\Big)^{\frac{1}{2+\EPS}}
\lesssim \Big(\fint_{2W}|\nabla v|^{2}\Big)^{\frac{1}{2}}+\Big(\fint_{2W}|h|^{2+\EPS}\Big)^{\frac{1}{2+\EPS}}.\nonumber\end{align}
Consequently, we obtain
\begin{align}
\Big(\fint_{\frac{3}{2}W}|\nabla (I-\mathfrak{A}_t)v|^{2+\EPS}\Big)^{\frac{1}{2+\EPS}}&\lesssim\Big(\fint_{\frac{3}{2}W}|\nabla v|^{2+\EPS}\Big)^{\frac{1}{2+\EPS}} + \Big(\fint_{\frac{3}{2}W}|\nabla \mathfrak{A}_tv|^{2+\EPS}\Big)^{\frac{1}{2+\EPS}}\nonumber
\\
&\lesssim \Big(\fint_{2W}|\nabla v|^{2}\Big)^{\frac{1}{2}}+\Big(\fint_{2W}|h|^{2+\EPS}\Big)^{\frac{1}{2+\EPS}}.
\nonumber
\end{align}
Similarly, but more directly, we also obtain that 
\begin{align}
    \Big(\fint_{\frac{3}{2}W}|(I-\mathfrak{A}_t)v|^{2+\EPS}\Big)^{\frac{1}{2+\EPS}} & \lesssim l(W)\Big(\fint_{\frac{3}{2}W}|\nabla v|^{2+\EPS}\Big)^{\frac{1}{2+\EPS}} \\
    &\lesssim l(W)\Big[\Big(\fint_{2W}|\nabla v|^{2}\Big)^{\frac{1}{2}}+\Big(\fint_{2W}|h|^{2+\EPS}\Big)^{\frac{1}{2+\EPS}}\Big].\label{eq:I-mathfrakA}
\end{align}

By using the boundedness of \(b_j^a\) in BMO we get \(\Big(\fint_{\frac{3}{2}W} |b^{a}_j-(b^{a}_j)_{\frac{3}{2}W}|^r \Big)^{\frac{1}{r}}\leq C\), and hence we can get in total for \(J_1^a\) with \reflemma{lemma:dyadicToNondyadic}
\begin{align*}
    J_1^a&\lesssim\sum_{j,W\in\mathcal{D}}|W|l(W)\Big[\Big(\fint_{\frac{3}{2}W}|\partial_{t}\mathcal{P}^R_tf|^2\Big)^\frac{1}{2}
    + \Big(\fint_{\frac{3}{2}W}|t\partial_j\partial_{t}\mathcal{P}^R_tf|^2\Big)^\frac{1}{2}\Big]
     \\ &\qquad{}\cdot\Big[\Big(\fint_{2W}|\nabla v|^{2}\Big)^{\frac{1}{2}}+\Big(\fint_{2W}|h|^{2+\EPS}\Big)^{\frac{1}{2+\EPS}}\Big]
    \\
    &\lesssim(\Vert \mathcal{S}(v)\Vert_{L^{p'}(\partial\Omega)}+\Vert \mathcal{T}_{2+\EPS}(h)\Vert_{L^{p'}})\big(\Vert \mathcal{A}(\partial_t\mathcal{P}^R_tf)\Vert_{L^p(\partial\Omega)}+\Vert \mathcal{A}(t\nabla_{||}\partial_t\mathcal{P}^R_tf)\Vert_{L^p(\partial\Omega)}\big)
    \\
    &\lesssim (\Vert \mathcal{S}(v)\Vert_{L^{p'}(\partial\Omega)}+\Vert \mathcal{T}_{q'}(h)\Vert_{L^{p'}})\Big[\Vert \mathcal{A}(\partial_t\mathcal{P}_tf)\Vert_{L^p(\partial\Omega)}+\Vert \mathcal{A}(t\nabla_{||}\partial_t\mathcal{P}_tf)\Vert_{L^p(\partial\Omega)}
    \\
    &\qquad+ \Vert \mathcal{A}(\partial_t\eta^R\mathcal{P}_tf)\Vert_{L^p(\partial\Omega)} + \Vert \mathcal{A}(t\partial_t\eta^R\nabla_{||}\mathcal{P}_tf)\Vert_{L^p(\partial\Omega)}\Big].
\end{align*}
By \refprop{prop:AreaFctBoundsInLp}, \eqref{eq:CarlAreaTcomparable}, and \eqref{eq:EstimatesFromPDforR} we can obtain the bound
\begin{align*}
J_1^a&\lesssim C\Vert\nabla_{||}f\Vert_{L^p}+\Vert \mathcal{A}(\partial_t\eta^R\mathcal{P}_tf)\Vert_{L^p(\partial\Omega)}+\Vert \mathcal{A}(t\partial_t\eta^R\nabla_{||}\mathcal{P}_tf)\Vert_{L^p(\partial\Omega)}.
\end{align*}
If we set the first error term as $E^R_2:=\Vert \mathcal{A}(\partial_t\eta^R\mathcal{P}_sf)\Vert_{L^p(\partial\Omega)}$, then \(J_1^a\lesssim C\Vert\nabla_{||}f\Vert_{L^p}+E_2^R+E_1^R\), and it remains to bound the error terms at the end.

Now, let us turn to \(J^s\).
We see that by the same splitting argument as for \(J^a\)
\begin{align*}
    J^s&=\sum_{1\leq j\leq n}\int_{\Omega^R} b^{s}_j \Big[\partial_t\mathcal{P}^R_t\, \partial_j(I-\mathfrak{A}_t)v +  \partial_j\mathcal{P}^R_tf\,\partial_t(I-\mathfrak{A}_t)v\Big] dxdt
    \\
    &\qquad +\sum_{1\leq j\leq n}\int_{\Omega^R} b^{s}_j \Big[\partial_t\mathcal{P}^R_tf\, \partial_j\mathfrak{A}_tv +  \partial_j\mathcal{P}^R_tf\,\partial_t \mathfrak{A}_tv\Big] dxdt
    \\
    &=:J_1^s+J_2^s.
\end{align*}
Justification of finiteness of both integrals will be provided in the following subsection.
Now, for \(J_1^s\) we use integration by parts in $t$ for the second term and the boundedness of \(b^s\) to get
\begin{equation}\label{eq:ConvJustJ_s^1}
\begin{split}
    J^s_1&=\sum_{1\leq j\leq n}\int_\Omega b^{s}_j \Big[\partial_{t}\mathcal{P}^R_tf\partial_j(I-\mathfrak{A}_t)v +  \partial_j\partial_t\mathcal{P}^R_tf(I-\mathfrak{A}_t)v\Big] dxdt
    \\
    &\lesssim \Vert \mathcal{A}(\partial_t\mathcal{P}^R_tf)\Vert_{L^p(\partial\Omega)}\Vert \mathcal{A}(s\nabla_{||}(I-\mathfrak{A}_s)v)\Vert_{L^{p'}(\partial\Omega)}+\Vert \mathcal{A}(t\nabla_{||}\partial_t\mathcal{P}^R_tf)\Vert_{L^p(\partial\Omega)}\Vert \mathcal{A}((I-\mathfrak{A}_t)v)\Vert_{L^{p'}(\partial\Omega)}
    \\
    &\lesssim \Big[\Vert \mathcal{A}(\partial_t\mathcal{P}_tf)\Vert_{L^p(\partial\Omega)}+E_2^R+\Vert \mathcal{A}(t\nabla_{||}\partial_t\mathcal{P}_sf)\Vert_{L^p(\partial\Omega)}+E_1^R\Big]\Vert \mathcal{S}(v)\Vert_{L^{p'}(\partial\Omega)},
\end{split}
\end{equation}
where we used arguments as in \eqref{eq:I-mathfrakA} for the area function of $(I-\mathfrak{A}_t)v$. As for \(J_1^a\), this is bounded by \(\Vert\nabla_{||}f\Vert_{L^p}+E_1^R+E_2^R\) according to \refprop{prop:AreaFctBoundsInLp} and \eqref{eq:EstimatesFromPDforR}.

Next, for \(J_2^s\) we have by integration by parts, first in $t$ for the second term, then in $j$, that
\begin{align*}
    J_2^s&=\sum_{1\leq j\leq n}\int_{\Omega^R} b^{s}_j \Big[\partial_{t}\mathcal{P}^R_tf\, \partial_j\mathfrak{A}_tv +  \partial_j\mathcal{P}^R_tf\, \partial_t \mathfrak{A}_tv\Big] dxdt
    \\
    &=\sum_{1\leq j\leq n}\int_{\Omega^R} b^{s}_j \Big(\partial_{t}\mathcal{P}^R_tf\, \partial_j\mathfrak{A}_tv\Big) dxdt - \int_{\Omega^R} b^{s}_j\Big( \partial_j\partial_{t}\mathcal{P}^R_tf \,\mathfrak{A}_tv\Big) dxdt
    \\
    &=\sum_{1\leq j\leq n}2\int_{\Omega^R} b^{s}_j \Big(\partial_{t}\mathcal{P}^R_tf \, \partial_j\mathfrak{A}_tv\Big) dxdt - \int_{\Omega^R} b^{s}_j \partial_j\big(\partial_{t}\mathcal{P}^R_tf \,\mathfrak{A}_tv\big) dxdt=J^s_{21}+J^s_{22}
\end{align*}
Due to the boundedness of \(b^s\), the first term \(J_s^{21}\) is bounded by
\(\Vert \mathcal{A}(\partial_t\mathcal{P}^R_tf)\Vert_{L^p}\Vert \mathcal{S}(v)\Vert_{L^{p'}}\), which is bounded by \(C\Vert\nabla_{||}f\Vert_{L^p}+CE_2^R\) as before.
Justification of finiteness of $J_{22}^s$ is postponed to the next subsection, and it can be combined with \(J^a_2\) to give \(\sq{J}:=J^s_{22}+J_2^a\) and by defining \(\sq{b}_j:=b^s_j-b^a_j\) we get
\begin{align*}
    \sq{J}&:=\int_{\Omega^R} b^s \nabla \big(\partial_{t}\mathcal{P}^R_tf \,\mathfrak{A}_tv\big) dxdt-\int_{\Omega^R} b^a \nabla \big(\partial_{t}^2\mathcal{P}^R_tf \,\mathfrak{A}_tv\big) dxdt
    \\
    &=\int_{\Omega^R} \sq{b} \nabla\big(\partial_{t}\mathcal{P}^R_tf\, \mathfrak{A}_tv\big) dxdt.
\end{align*}
From here onwards, the only property used for the vector \(\sq{b}(x)\) is that \(\sq{b}\in BMO(\mathbb{R}^n)\) which follows since \(b^s\in L^\infty\) and \(b^a\in BMO(\Omega)\). We introduce a smooth cut-off function $\psi^R\in C^\infty(\Delta(0,2R))$ in the spatial variable $x$ only, with $0\leq\psi^R\leq 1$, $\psi^R\equiv 1$ on $\Delta(0,R)$, and $|\nabla_{||}\psi^R|\lesssim \frac{1}{R}$. 
We can continue to split the integral and get
\begin{align*}
    \sq{J}&=\int_{\Omega^R} \sq{b} \nabla_{||}\big(\psi^R\partial_{t}\mathcal{P}^R_tf \mathfrak{A}_tv\big) dxdt+\int_{\Omega^R} \sq{b} \nabla_{||}\big((1-\psi^R)\partial_{t}\mathcal{P}^R_tf \mathfrak{A}_tv\big) dxdt
    \\
    &=:\sq{J}_1+E^R_3.
\end{align*}
Finiteness of $\sq{J}_1$ will be proved in the following subsection.
We postpone the error term $E_3^R$ until the end and proceed with $\sq{J}_1$. Since $\sq{b}\in BMO$, we have $\sq{b}\in L^2(\Delta(0,2R))$

Thus, 
\begin{align*}
    \sq{J}_1&=\int_{\Omega^R} \sq{b} \nabla_{||}\big(\psi^R\partial_{t}\mathcal{P}^R_tf \mathfrak{A}_tv\big) dxdt
    \\
    &=\int_{\Omega^R} \sq{b} \nabla_{||}\Big[\psi^R \big(-2\eta^RtL_{||}\mathcal{P}_tf+\partial_t\eta^R\mathcal{P}_tf\big)\mathfrak{A}_tv\Big] dxdt
    \\
    &=\int_{\Omega^R} \eta^R\mathcal{P}_{\sqrt{3}t/2}^*\big(\mathrm{div}(\sq{b})\psi^R \mathfrak{A}_sv\big)(-2tL_{||}\mathcal{P}_{t/2}f)\,dxdt
     +\int_{\Omega^R}\sq{b}\nabla_{||}\Big[\psi^R \mathfrak{A}_tv\big(\partial_t\eta^R\mathcal{P}_tf\big)\Big] dxdt
     \\
    &=:\sq{J}_{11}+E^R_4.
\end{align*}
We will deal with the error term $E^R_4$ in the end. We note here that we understand $\mathcal{P}_{\sqrt{3}t/2}^*\big(\mathrm{div}(\sq{b})\psi^R \mathfrak{A}_tv\big)$ as 
$ \mathcal{P}_{\sqrt{3}t/2}^*\mathrm{div}\big(\sq{b}\psi^R \mathfrak{A}_tv\big)- \mathcal{P}_{\sqrt{3}t/2}^*\big(\sq{b}\cdot\nabla(\psi^R\eta^R \mathfrak{A}_tv)\big)$, and since $\psi^R$ is smooth, $b\in L^2(\mathrm{supp}(\psi^R))$, and $\mathfrak{A}_tv, \nabla\mathfrak{A}_t(v)\in L^\infty(\mathrm{supp}(\psi^R))$ for each fixed $t$, both these operators act on $L^2$ functions, and can be understood in the classical sense. Convergence of the integral $\sq{J}_{11}$ also needs justification.

Now to be precise, we introduce a rough averaging operator
\begin{align}\overline{\mathfrak{A}}_tf(x):=\sum_{W\in\mathcal{D}(\Omega)}\chi_W(x,t) \fint_{W_t}f(y,t)dy,\qquad\textrm{where }\quad W_t:=W\cap\{(x,s)\in\mathbb{R}^{n+1}_+|s=t\}, \label{eq:DefRoughAveragingOp}\end{align}
and can split again
\begin{align*}
    \sq{J}_{11}&=\int_{\Omega^R} \eta^R\mathcal{P}_{\sqrt{3}t/2}^*\big(\mathrm{div}(\sq{b})\psi^R \mathfrak{A}_tv\big)(-2tL_{||}\mathcal{P}_{t/2}f)\,dxdt
    \\
    &=\int_{\Omega^R} \eta^R\Big[\mathcal{P}_{\sqrt{3}t/2}^*(\mathrm{div}(\sq{b})\psi^R\mathfrak{A}_tv)-\mathcal{P}^*_{\sqrt{3}t/2}(\mathrm{div}(\sq{b})\psi^R)\overline{\mathfrak{A}}_tv\Big](-2tL_{||}\mathcal{P}_{t/2}f)\, dxdt 
    \\
    &\qquad + \int_{\Omega^R} \eta^R\mathcal{P}^*_{\sqrt{3}t/2}(\mathrm{div}(\sq{b})\psi^R)\, \overline{\mathfrak{A}}_tv \,(-2tL_{||}\mathcal{P}_{t/2}f)\,dxdt
    \\
    &=:\sq{J}_{111}+\sq{J}_{112},
\end{align*}
Now, since $-2tL_{||}\mathcal{P}_{t/2}f=4\partial_{t}\mathcal{P}_{t/2}f$, a stopping time argument yields for \(\sq{J}_{112}\)
\[\sq{J}_{112}\lesssim \Vert\mathcal{A}(\partial_{t}\mathcal{P}_{t/2}f)\Vert_{L^p}\Vert\mathcal{C}(t\eta^R\mathcal{P}_{\sqrt{3}t/2}^*(\mathrm{div}(\sq{b})\psi^R))\Vert_{L^\infty}\Vert \mathcal{N}(\overline{\mathfrak{A}}_tv)\Vert_{L^{p'}}.\]
Since \(\Vert \mathcal{N}_\alpha(\overline{\mathfrak{A}}_tv)\Vert_{L^{p'}}\lesssim \Vert \mathcal{N}_{\alpha+2}(v)\Vert_{L^{p'}}\), \refprop{prop:AreaFctBoundsInLp} and \eqref{eq:EstimatesFromPDforR} imply that 
\[\sq{J}_{112}\lesssim \Vert\mathcal{C}(t\eta^R\mathcal{P}_{\sqrt{3}t/2}^*(\mathrm{div}(\sq{b})\psi^R))\Vert_{L^\infty}\Vert\nabla_{||}f\Vert_{L^p}.\]
Note that the bounds in \refprop{prop:AreaFctBoundsInLp} are also true for $t/2$ instead of $t$ by change of variables. We claim that \(\Vert\mathcal{C}(t\eta^R\mathcal{P}_{\sqrt{3}t/2}^*(\mathrm{div}(\sq{b})\psi^R))\Vert_{L^\infty}\lesssim \Vert \sq{b}\Vert_{BMO}<\infty\), whose proof will be postponed to \reflemma{lemma:CarlesonFctBoundedBMO} in the following subsection.

Lastly, for \(\sq{J}_{111}\) we claim that
\[\sq{J}_{111}\lesssim \Vert \mathcal{A}(\partial_t\mathcal{P}_{t/2}f)\Vert_{L^p}(\Vert \mathcal{S}(v)\Vert_{L^{p'}}+\Vert \mathcal{T}_{2+\EPS}(h)\Vert_{L^{p'}}),\]
which will be proved in \reflemma{lemma:DifferenceOperatorRisbounded} in the following subsection. With this claim, \eqref{eq:CarlAreaTcomparable}, \eqref{eq:EstimatesFromPDforR} and \refprop{prop:AreaFctBoundsInLp}, we obtain \(J_{111}\lesssim\Vert\nabla_{||}f\Vert_{L^p}\).

\medskip

To summarize, up to this point we have proved the following: Assuming convergence of all integrals and assuming \reflemma{lemma:DifferenceOperatorRisbounded} and \reflemma{lemma:CarlesonFctBoundedBMO}, if $f\in \dot{W}^{1,p}(\partial\Omega)\cap C_c^\infty(\partial\Omega)$, then there exists a sufficiently large $R_0>0$, such that $\mathrm{supp}(f)\subset B(0,R_0/4)$, and for all $R>R_0$
\[\Vert\sq{\mathcal{N}}_2(\nabla u)\Vert_{L^p(\partial\Omega)}\lesssim \Vert\nabla_{||}f\Vert_{L^p}+\sum_{i=1}^4 E_i^R.\]

Hence, it remains to prove, first, the convergence of all integrals, second, that $E_i^R\to 0$ as $R$ tends to infinity, and lastly \reflemma{lemma:DifferenceOperatorRisbounded} and \reflemma{lemma:CarlesonFctBoundedBMO}.

\subsubsection{Convergence of integrals}
Here, we justify the finiteness of all integrals appearing the above proof. Recall that $f=f_\EPS=\mathcal{P}_\EPS g$ for $g\in \dot{W}^{1,2}(\partial\Omega)\cap\dot{W}^{1,p}(\partial\Omega)$. This implies that $\Vert \nabla_{||}L_{||}f\Vert_{L^2(\mathbb{R}^n)}$, $\Vert L_{||}f\Vert_{L^2(\mathbb{R}^n)}$, $\Vert \nabla_{||}f\Vert_{L^2(\mathbb{R}^n)}$, and $\Vert f\Vert_{L^2(\mathbb{R}^n)}$ are finite. Furthermore, we note that to prove finiteness of the integrals it is enough to prove boundedness by a constant depending on these terms and $R$. Moreover, the following observation for two functions \(\sq{g}\) and \(\sq{h}\) will be helpful: By $\mathcal{H}^1-BMO$ duality and \refprop{prop:DivCurl}
\begin{align}
    \int_{\Omega} b^a_j(\partial_j \sq{g}\partial_t \sq{h}-\partial_t\sq{g}\partial_j\sq{h}) dx\lesssim \Vert A^a\Vert_{BMO(\Omega)}\Vert\nabla \sq{g}\Vert_{L^2(\Omega)}\Vert\nabla \sq{h}\Vert_{L^2(\Omega)}\lesssim \Vert\nabla \sq{g}\Vert_{L^2(\Omega)}\Vert\nabla \sq{h}\Vert_{L^2(\Omega)}\label{eq:JustificationConvIntegralsb}.
\end{align}
Now let us show convergence of each integral appearing in the proof of \refthm{thm:MainThm2} in turn.

\begin{itemize}
    \item[$\mathbf{I_1:}$] We can split the integral \(I_1\) into
    \[\int_\Omega \AP\nabla_{||}\mathcal{P}^R_tf\cdot \nabla_{||}v \,dx=\int_\Omega \AP^s\nabla_{||}\mathcal{P}^R_tf\cdot \nabla_{||}v \,dx+\int_\Omega \AP^a\nabla_{||}\mathcal{P}^R_tf\cdot \nabla_{||}v \,dx,\]
and use Cauchy-Schwarz, $\mathcal{H}^1-BMO$ duality, and \refprop{prop:DivCurl}, respectively, to obtain
\begin{align*}
    \int_{\Omega}\AP^s\nabla_{||}\mathcal{P}^R_tf\cdot \nabla_{||}v \, dx\lesssim \Vert A^s\Vert_{L^\infty}\Vert\nabla_{||}\mathcal{P}_tf\Vert_{L^2(\Omega)}\Vert\nabla_{||}v\Vert_{L^2(\Omega)}\lesssim \Vert\nabla_{||}\mathcal{P}_tf\Vert_{L^2(\Omega)}\Vert\nabla_{||}v\Vert_{L^2(\Omega)},
\end{align*}
and 
\begin{align*}
    \int_{\Omega}\AP^a\nabla_{||}\mathcal{P}^R_tf\cdot \nabla_{||}v \, dx\lesssim \Vert A^a\Vert_{BMO(\Omega)}\Vert\nabla_{||}\mathcal{P}_tf\Vert_{L^2(\Omega)}\Vert\nabla_{||}v\Vert_{L^2(\Omega)}\lesssim \Vert\nabla_{||}\mathcal{P}_tf\Vert_{L^2(\Omega)}\Vert\nabla_{||}v\Vert_{L^2(\Omega)}.
\end{align*}
Hence, by \refprop{prop:L2boundednessOfSemigroupOperators} and \(v\in Y^{1,2}_0(\Omega\), 
\[I_1\lesssim \Vert\nabla_{||}\mathcal{P}_tf\Vert_{L^2(\Omega)}\Vert\nabla_{||}v\Vert_{L^2(\Omega)}\lesssim R^{1/2}\Vert\nabla_{||}f\Vert_{L^2(\Omega)}\Vert\nabla v\Vert_{L^2(\Omega)}<\infty.\]
\item[$\mathbf{I_2}$:] follows immediately for the bound in \eqref{eq:BoundI2}.

\item[$\mathbf{I_3}$:]  The finiteness of $I_3$ follows from the finiteness of $J^a$ and $J^s$.
\item[$\mathbf{I_{11}}$:]  The finiteness of $I_{11}$ follows from finiteness of $I^s_{11}$ and $I^a_{11}$
\item[$\mathbf{I^s_{11}}$:] The bound of $I_{11}^s$ also implies finiteness.
\item[$\mathbf{I^a_{11}}$:] This bound is similar to $I_1$. By \(\mathcal{H}^1-BMO\) duality, \refprop{prop:DivCurl}, and \refprop{prop:L2boundednessOfSemigroupOperators}, we get
\begin{align*}
    \int_{\Omega}\AP^a\nabla_{||}\partial_t\mathcal{P}^R_tf\cdot \nabla_{||}v t \, dx\lesssim \Vert A^a\Vert_{BMO(\Omega)}\Vert t\nabla_{||}\partial_t\mathcal{P}_tf\Vert_{L^2(\Omega)}\Vert\nabla_{||}v\Vert_{L^2(\Omega)}\lesssim \Vert\nabla_{||}f\Vert_{L^2(\Omega)}\Vert\nabla_{||}v\Vert_{L^2(\Omega)},
\end{align*}
which is finite, since \(\nabla f\in L^2(\partial\Omega)\) and \(v\in Y^{1,2}_0(\Omega)\).
\item[$\mathbf{I_{12}}$:]  The finiteness of $I_{12}$ follows from the equality \(I_{1}=I_{11}+I_{12}\) and finiteness of both $I_1$ and $I_{11}$.
\item[$\mathbf{J^a}$:] By \eqref{eq:JustificationConvIntegralsb} for \(\sq{g}=\mathcal{P}^R_tf\) and \(\sq{h}=v\), we have
\begin{align*}
    J^a&\lesssim\Vert \nabla_{x,t}\mathcal{P}^R_tf\Vert_{L^2(\Omega)}\Vert\nabla v\Vert_{L^2(\Omega)}
    \\
    &\lesssim \big(\Vert \nabla_{||}\mathcal{P}^R_tf\Vert_{L^2(\Omega)}+\Vert 2t\eta^R\mathcal{P}_t(L_{||}f)\Vert_{L^2(\Omega)}+\Vert \partial_t\eta^R\mathcal{P}_tf\Vert_{L^2(\Omega)}\big)\Vert\nabla v\Vert_{L^2(\Omega)}
\end{align*}
Since \(v\in Y^{1,2}\) and $A\in BMO(\Omega)$, we get by \refprop{prop:L2boundednessOfSemigroupOperators}
\begin{align*}J^a&\lesssim \Vert \eta^R\nabla_{||}\mathcal{P}_tf\Vert_{L^2(\Omega)}+\Vert 2t\eta^R\mathcal{P}_t(L_{||}f)\Vert_{L^2(\Omega)}+\Vert \partial_t\eta^R\mathcal{P}_tf\Vert_{L^2(\Omega)}
\\
&\lesssim \Big(\int_0^R\Vert \nabla_{||} f\Vert_{L^2(\mathbb{R}^n)}^2 + R\Vert L_{||} f\Vert_{L^2(\mathbb{R}^n)}^2 + \frac{1}{R}\Vert f\Vert_{L^2(\mathbb{R}^n)}^2 dt\Big)^{1/2}=C(R)<\infty.\end{align*} 

\item[$\mathbf{J^s}$:]  For $J^s$ we can apply Cauchy-Schwarz and get
\[J^s\leq \Vert A^s\Vert_{L^\infty}\Vert \nabla_{x,t}\mathcal{P}^R_tf\Vert_{L^2(\Omega)}\Vert\nabla v\Vert_{L^2(\Omega)}\lesssim \Vert \nabla_{x,t}\mathcal{P}^R_tf\Vert_{L^2(\Omega)}\Vert\nabla v\Vert_{L^2(\Omega)}.\]

From here, we can proceed verbatim as for $J^a$.
\item[$\mathbf{J^a_1}$:] By \eqref{eq:JustificationConvIntegralsb} for \(\sq{g}=\mathcal{P}^R_tf\) and \(\sq{h}=(I-\mathfrak{A}_t)v\), we get
\[|J_1^a|\lesssim\Vert \nabla_{x,t}\mathcal{P}^R_tf\Vert_{L^2(\Omega)}\Vert \nabla(I-\mathfrak{A}_t)v \Vert_{L^2(\Omega)}.\]
Since \(\Vert \nabla(I-\mathfrak{A}_t)v \Vert_{L^2(\Omega)}\lesssim \Vert \nabla v \Vert_{L^2(\Omega)}\), we can proceed verbatim to \(J^a\) from here.
\item[$\mathbf{J^a_2}$:] By \eqref{eq:JustificationConvIntegralsb} for \(\sq{g}=\mathcal{P}^R_tf\) and \(\sq{h}=(I-\mathfrak{A}_t)v\), we get
\[|J_2^a|\lesssim\Vert \nabla_{x,t}\mathcal{P}^R_tf\Vert_{L^2(\Omega)}\Vert \nabla \mathfrak{A}_tv \Vert_{L^2(\Omega)}.\]
Since \(\Vert \nabla \mathfrak{A}_tv \Vert_{L^2(\Omega)}\lesssim \Vert \nabla v \Vert_{L^2(\Omega)}\), we can proceed as for \(J^a\).
\item[$\mathbf{J^s_1}$:]  The bound of $J^s_1$ in \eqref{eq:ConvJustJ_s^1} also proves finiteness.
\item[$\mathbf{J^s_2}$:]  The finiteness of $J^s_2$ follows from the finiteness of $J^s_{21}$ and $J^s_{22}$.
\item[$\mathbf{J^s_{21}}$:]  The bound of $J^s_1$ also proves finiteness.
\item[$\mathbf{J^s_{22}}$:] We can apply Cauchy-Schwarz to get 
\[|J^s_{22}|\lesssim \Vert \nabla_{x,t}\partial_t\mathcal{P}^R_tf\Vert_{L^2(\Omega^R)}\Vert \mathfrak{A}_tv \Vert_{L^2(\Omega^R)}+\Vert \partial_t\mathcal{P}^R_tf\Vert_{L^2(\Omega)}\Vert \nabla\mathfrak{A}_tv \Vert_{L^2(\Omega)}.\]
From here on, we can proceed verbatim as for $J_2^a$.
\item[$\mathbf{\sq{J}}$:]  The finiteness of $\sq{J}$ follows from the finiteness of $\sq{J}_1$ and the boundedness of the error term $E_5^R$.
\item[$\mathbf{\sq{J}_{1}}$:] By $\mathcal{H}^1-BMO$ duality and \refprop{prop:H1normofproduct}
\[|\sq{J}_1|\lesssim \Vert \nabla_{x,t}\partial_t\mathcal{P}^R_tf\Vert_{L^2(\Omega^R)}\Vert \mathfrak{A}_tv \Vert_{L^2(\Omega^R)}+\Vert \partial_t\mathcal{P}^R_tf\Vert_{L^2(\Omega)}\Vert \nabla\mathfrak{A}_tv \Vert_{L^2(\Omega)}.\]
Since \(\Vert \partial_t\mathcal{P}^R_tf\Vert_{L^2(\Omega)}\lesssim \Vert \nabla_{x,t}\mathcal{P}^R_tf\Vert_{L^2(\Omega)}\), we have finiteness of \(\Vert \partial_t\mathcal{P}^R_tf\Vert_{L^2(\Omega)}\Vert \nabla\mathfrak{A}_tv \Vert_{L^2(\Omega)}\) as in \(J^a\).

For \(\Vert \nabla_{x,t}\partial_t\mathcal{P}^R_tf\Vert_{L^2(\Omega^R)}\Vert \mathfrak{A}_tv \Vert_{L^2(\Omega^R)}\), we observe that Poincar\'{e}'s inequality applies, since \( v\) vanishes at the boundary of \(\Omega\), and implies \(\Vert \mathfrak{A}_tv \Vert_{L^2(\Omega^R)}\lesssim\Vert v \Vert_{L^2(\Omega^R)}\lesssim R\Vert \nabla v\Vert_{L^2(\Omega^R)}\), which is finite. Moreover,
\begin{align*}
    \Vert \nabla_{x,t}\partial_t\mathcal{P}^R_tf\Vert_{L^2(\Omega)}&\lesssim \Vert (-2L_{||}+4t^2L_{||}L_{||})\mathcal{P}_tf\Vert_{L^2(\Omega)}+ \Vert 2t\partial_t\eta^RL_{||}\mathcal{P}_tf\Vert_{L^2(\Omega)}
    \\
    &\qquad + \Vert \partial_t^2\eta^R \mathcal{P}_tf\Vert_{L^2(\Omega)}+\Vert 2\nabla_{||}L_{||}\mathcal{P}_tf\Vert_{L^2(\Omega)}+ \Vert \partial_t\eta^R\nabla_{||}\mathcal{P}_tf\Vert_{L^2(\Omega)}
    \\
    &\lesssim R^2\Vert L_{||}L_{||}f\Vert_{L^2(\mathbb{R}^n)}+ \Vert L_{||}f\Vert_{L^2(\mathbb{R}^n)}
    \\
    &\qquad+ \frac{1}{R^2}\Vert f\Vert_{L^2(\mathbb{R}^n)}+\Vert \nabla_{||}L_{||}f\Vert_{L^2(\mathbb{R}^n)}+\frac{1}{R}\Vert \nabla_{||}f\Vert_{L^2(\mathbb{R}^n)}. 
\end{align*}
Hence \(\Vert \nabla_{x,t}\partial_t\mathcal{P}^R_tf\Vert_{L^2(\Omega^R)}\Vert \mathfrak{A}_tv \Vert_{L^2(\Omega^R)}\) is also finite and \(\sq{J}_1\) converges.
\item[$\mathbf{\sq{J}_{11}}$:] By Cauchy-Schwarz and \refprop{prop:L2boundednessOfSemigroupOperators}, we have
\begin{align*}
    &\int_{\Omega^R} \eta^R\mathcal{P}_{\sqrt{3}t/2}^*\big(\mathrm{div}(\sq{b})\psi^R \mathfrak{A}_tv\big)(-2tL_{||}\mathcal{P}_{t/2}f)dxdt
    \\
    &\lesssim \int_0^R \big(\Vert\mathcal{P}_{\sqrt{3}t/2}^*\big(\sq{b}\nabla(\psi^R \mathfrak{A}_tv)\big)\Vert_{L^2(\mathbb{R}^n)}+\Vert\mathcal{P}_{\sqrt{3}t/2}^*\mathrm{div}\big(\sq{b}\psi^R \mathfrak{A}_tv\big)\Vert_{L^2(\mathbb{R}^n)}\big)\Vert tL_{||}\mathcal{P}_{t/2}f\Vert_{L^2(\mathbb{R}^n)}dt
    \\
    &\lesssim \int_0^R \big(\Vert\sq{b}\nabla(\psi^R \mathfrak{A}_tv)\Vert_{L^2(\mathbb{R}^n)}+\frac{1}{t}\Vert\sq{b}\psi^R \mathfrak{A}_tv\Vert_{L^2(\mathbb{R}^n)}\big)\Vert tL_{||}f\Vert_{L^2(\mathbb{R}^n)}dt
    \\
    &\lesssim (\frac{1}{R}+1+R)\int_0^R \big(\Vert\sq{b}\Vert_{L^{2(2+\EPS)/\EPS}(\Delta(0,2R))}\Vert  \nabla v\Vert_{L^{2+\EPS}(\Delta(0,2R))}+\Vert\sq{b}\Vert_{L^{n}(\Delta(0,2R))}\Vert v\Vert_{L^{2^*}(\mathbb{R}^n)}\big)\Vert L_{||}f\Vert_{L^2(\mathbb{R}^n)}ds,
\end{align*}
where \(\EPS>0\) is chosen small enough such the reverse H\"{o}lder inequality \reflemma{lemma:HigherIntegrabilityForGradient} applies. Then
\begin{align*}
    \sq{J}_{11}&\lesssim (\frac{1}{R}+1+R)(R^\frac{1+\EPS}{2+\EPS}\Vert\nabla v\Vert_{L^{2+\EPS}(T(\Delta(0,2R)))}+R^\frac{2^*-1}{2^*}\Vert v\Vert_{L^{2^*}(\mathbb{R}_+^{n+1})}\big)\Vert L_{||}f\Vert_{L^2(\mathbb{R}^n)},
\end{align*}
Since, by \reflemma{lemma:HigherIntegrabilityForGradient},
\[\Vert\nabla v\Vert_{L^{2+\EPS}(T(\Delta(0,2R)))}\lesssim R^{-\frac{\EPS(n+1)}{2(2+\EPS)}}\Vert\nabla v\Vert_{L^{2}(T(\Delta(0,4R)))}+\Vert h\Vert_{L^{2+\EPS}(T(\Delta(0,2R)))},\]
and $h\in L_c^\infty(\Omega)$, we obtain finiteness of \(\Vert\nabla v\Vert_{L^{2+\EPS}(T(\Delta(0,2R)))}\) and hence convergence of the integral $\sq{J}_{11}$.
\item[$\mathbf{\sq{J}_{111}}$:]  The bound of $\sq{J}_{111}$ also proves finiteness.
\item[$\mathbf{\sq{J}_{112}}$:]  The bound of $\sq{J}_{112}$ also proves finiteness.
\end{itemize}

\subsubsection{Error term bounds}

We are going to show that for every $i$ the error term $E_i^R\to 0$ as $R\to \infty$.
\begin{itemize}
\item[$E_1^R$:] Recall that $E_1^R=\Vert\mathcal{A}(t\partial_t\eta^R\nabla_{||}\mathcal{P}_tf)\Vert_{L^p}$. For a fixed $x\in \mathbb{R}^n$, we can see that
\begin{align*}
\mathcal{A}(t\partial_t\eta^R\nabla_{||}\mathcal{P}_tf)^2(x)&=\int_{\Gamma(x)}\frac{|s\partial_s\eta^R\nabla_{||}\mathcal{P}_sf(x)|^2}{s^{n+1}}dxds
\\
&=\fint_{\Gamma(x)\cap(\mathbb{R}^n\times(R/2,R))}|s\partial_s\eta^R\nabla_{||}\mathcal{P}_sf(x)|^2dxds
\\
&\lesssim \fint_{\Delta(x,{2R})\times(R/2,R)}|\nabla_{||}\mathcal{P}_sf(x)|^2dxds
\end{align*}
By the Caccioppoli inequality \reflemma{lemma:CaccTypeInequalityForPSimple} for \(\mathcal{P}_t\), we get
\begin{align*}
    \fint_{\Delta(x,{2R})\times(R/2,R))}|\nabla_{||}\mathcal{P}_sf(x)|^2dxds&\lesssim \frac{1}{R^2}\fint_{\Delta(x,{2R})\times(R/2,R)}|\mathcal{P}_sf(x)|^2dxds
    \\
    &\lesssim \frac{1}{R}\mathcal{N}(\mathcal{P}_tf)(x)^2.
\end{align*}

Hence, by \refrem{rem:NontangnetialBoundsForSemigroup}, 
\[|E_1^R|\lesssim \frac{1}{R}\Vert\mathcal{N}(\mathcal{P}_tf)\Vert_{L^p}\lesssim \frac{1}{R}\Vert f\Vert_{L^p}\to 0\qquad\textrm{ as }R\to \infty.\] 

\item[$E_2^R$:] Recall that $E_2^R=\Vert\mathcal{A}(\partial_t\eta^R\mathcal{P}_tf)\Vert_{L^p}$. Similar to \(E_1^R\), we fix a $x\in \mathbb{R}^n$, and have
\begin{align*}
\mathcal{A}(\partial_t\eta^R\mathcal{P}_tf)^2(x)&=\int_{\Gamma(x)}\frac{|\partial_s\eta^R\mathcal{P}_sf(x)|^2}{s^{n+1}}dxds=\fint_{\Gamma(x)\cap(\mathbb{R}^n\times(R/2,R))}|\partial_s\eta^R\mathcal{P}_sf(x)|^2dxds
\\
&\lesssim \frac{1}{R^2}\sq{\mathcal{N}}_2(\mathcal{P}_sf)^2(x)
\end{align*}
Hence, by \refrem{rem:NontangnetialBoundsForSemigroup}, 
\[|E_1^R|\lesssim \frac{1}{R}\Vert\sq{\mathcal{N}}_2(\mathcal{P}_sf)\Vert_{L^p}\lesssim \frac{1}{R}\Vert f\Vert_{L^p}\to 0\qquad\textrm{ as }R\to \infty.\] 

\item[$E_3^R$:]  For $E^R_3$ we use an annular decomposition combined with off-diagonal estimates of \(\mathcal{P}_s\). Recall that $\psi^R\in C_0^\infty(\Delta(0,2R))$ with $0\leq\psi^R\leq 1$, $\psi^R\equiv 1$ on $\Delta(0,R)$, and $|\nabla_{||}\psi^R|\lesssim \frac{1}{R}$, and let $\sq{\psi}_k\in C_0^\infty(\Delta(0,2^{k+1}R)\cup\Delta(0,2^{k-1}R))$ be nonnegative smooth cut-off functions such that \(\sum_{k\geq 1}\sq{\psi}^2_k=(1-\psi^R)\) with the uniform bound \(|\nabla\sq{\psi}_k|\lesssim C\). By \refprop{prop:H1normofproduct}, we have
\begin{align*}
    |E^R_3|&=\Big|\int_{\Omega^R} \sq{b} \nabla_{||}\big((1-\psi^R)\partial_{t}\mathcal{P}^R_tf \mathfrak{A}_tv\big) dxdt\Big|
    \\
    &\lesssim\sum_{k\geq 1}\int_0^R\Big|\int_{\mathbb{R}^n} \sq{b} \nabla_{||}\big(\sq{\psi}_k^2\partial_{t}\mathcal{P}^R_tf \mathfrak{A}_tv\big) dx\Big| dt
    \\
    &\lesssim\sum_{k\geq 1}\int_0^R \Vert\sq{b}\Vert_{BMO}\Big(\Vert \nabla_{||}\big(\sq{\psi}_k\partial_{t}\mathcal{P}^R_tf\big)\Vert_{L^2(\mathbb{R}^n)} \Vert\sq{\psi}_k\mathfrak{A}_tv\Vert_{L^2(\mathbb{R}^n)}
    \\
    &\qquad\qquad +\Vert \sq{\psi}_k\partial_{t}\mathcal{P}^R_tf\Vert_{L^2(\mathbb{R}^n)} \Vert\nabla_{||}(\sq{\psi}_k\mathfrak{A}_tv)\Vert_{L^2(\mathbb{R}^n)}\Big)dt.
\end{align*}
Since $\mathrm{supp}f\subset\Delta(0,R/4)$, off-diagonal estimates \refprop{prop:Off-diagonalEstimates}, and \refprop{prop:L2boundednessOfSemigroupOperators} imply
\begin{align*}
    &\Vert \nabla_{||}\big(\sq{\psi}_k\partial_{t}\mathcal{P}^R_tf\big)\Vert_{L^2(\mathbb{R}^n)}
    \\
    &\lesssim \Vert\nabla_{||}\partial_{t}\mathcal{P}^R_tf\Vert_{L^2(\Delta(0,2^{k+1}R)\setminus \Delta(0,2^{k-1}R))} + \frac{1}{R}\Vert \partial_{t}\mathcal{P}^R_tf\Vert_{L^2(\Delta(0,2^{k+1}R)\setminus \Delta(0,2^{k-1}R))}
    \\
    &\lesssim \Vert \nabla_{||}\partial_{t}\mathcal{P}_tf\Vert_{L^2(\Delta(0,2^{k+1}R)\setminus \Delta(0,2^{k-1}R))}+\Vert \partial_t\eta^R\nabla_{||}\mathcal{P}_tf\Vert_{L^2(\Delta(0,2^{k+1}R)\setminus \Delta(0,2^{k-1}R))} 
    \\
    &\qquad
    +  \frac{1}{R}\Vert \partial_{t}\mathcal{P}_tf\Vert_{L^2(\Delta(0,2^{k+1}R)\setminus \Delta(0,2^{k-1}R))}+\frac{1}{R}\Vert \partial_t\eta^R\mathcal{P}_tf\Vert_{L^2(\Delta(0,2^{k+1}R)\setminus \Delta(0,2^{k-1}R))}
    \\
    &\lesssim \big(\frac{1}{t^2} + \frac{1}{tR}+\frac{1}{R^2}\big) e^{-c\frac{2^{2k}R^2}{t^2}}\Vert f\Vert_{L^2(\mathbb{R}^n)}.  
\end{align*}
Since $t\mapsto \big(\frac{1}{t^2} + \frac{1}{tR}+\frac{1}{R^2}\big) e^{-c\frac{2^{2k}R^2}{t^2}}$ is monotonically increasing on $(0,R)$, we can bound 
\[\sup_{t\in (0,R)}\left[\big(\frac{1}{t^2} + \frac{1}{tR}++\frac{1}{R^2}\big) e^{-c\frac{2^{2k}R^2}{t^2}}\Vert f\Vert_{L^2(\mathbb{R}^n)}\right]\lesssim \frac{1}{R^2}e^{-2^{2k}}\Vert f\Vert_{L^2(\mathbb{R}^n)}.\]
Analogously, 
\[\Vert \sq{\psi}_k\partial_{t}\mathcal{P}^R_tf\Vert_{L^2(\mathbb{R}^n)}\lesssim \frac{1}{R}\Vert f\Vert_{L^2(\mathbb{R}^n)}\]
Moreover,
\begin{align*}
    &\Vert\nabla_{||}(\sq{\psi}_k\mathfrak{A}_tv)\Vert_{L^2(\mathbb{R}^n)}
    \lesssim \Vert \nabla v(\cdot,t)\Vert_{L^2(\Delta(0,2^{k+1})\setminus\Delta(0,2^{k-1}))}+\frac{1}{R}\Vert v(\cdot,t)\Vert_{L^2(\Delta(0,2^{k+1})\setminus\Delta(0,2^{k-1}))}
    \\
    &\lesssim \Vert \nabla v(\cdot,t)\Vert_{L^2(\Delta(0,2^{k+1})\setminus\Delta(0,2^{k-1}))}+(2^k)^{\frac{n}{n+1}}R^{-\frac{1}{n+1}}\Vert v(\cdot,t)\Vert_{L^{2^*}(\Delta(0,2^{k+1})\setminus\Delta(0,2^{k-1}))}.
\end{align*}
and
\begin{align*}
    &\Vert(\sq{\psi}_k\mathfrak{A}_tv)\Vert_{L^2(\mathbb{R}^n)}
    \lesssim\Vert v(\cdot,t)\Vert_{L^2(\Delta(0,2^{k+1})\setminus\Delta(0,2^{k-1}))}\lesssim(2^kR)^{\frac{n}{n+1}}\Vert v(\cdot,t)\Vert_{L^{2^*}(\Delta(0,2^{k+1})\setminus\Delta(0,2^{k-1}))}.
\end{align*}

Consequently,
\begin{align*}
|E^R_3|&\leq\Vert f\Vert_{L^2(\mathbb{R}^n)}\sum_{k\geq 1}\int_0^R  \frac{\Vert \nabla v(\cdot,t)\Vert_{L^2(\Delta(0,2^{k+1})\setminus\Delta(0,2^{k-1}))}}{R}e^{-2^{2k}}dt
\\
&\qquad +\Vert f\Vert_{L^2(\mathbb{R}^n)}\sum_{k\geq 1}\int_0^R  \frac{(2^kR)^{\frac{n}{n+1}}\Vert v(\cdot,t)\Vert_{L^{2^*}(\Delta(0,2^{k+1})\setminus\Delta(0,2^{k-1}))}}{R^2}e^{-2^{2k}}dt
\\
&\lesssim \|f\|_2\sum_{k\ge1}e^{-2^{2k}}\big(
R^{-3/2}
\Vert \nabla v\Vert_{L^2(\Omega)}+R^{-1/2}2^{\frac{kn}{n+1}}\Vert v\Vert_{L^{2^*}(\Omega)}
\big).
\\
&\lesssim \frac{\Vert f\Vert_{L^2(\mathbb{R}^n)}}{\sqrt{R}}\Vert v\Vert_{Y^{1,2}(\Omega^R)}\to 0
\end{align*}
as $R\to \infty$.

\item[$E_4^R$:] Lastly, $E_4^R$ can be handled with ideas from $E_3^R$. We have by \refprop{prop:H1normofproduct}
\begin{align*}
|E_4^R|&\leq\int_{R/2}^R\Big|\int_{\mathbb{R}^n}\sq{b}\nabla_{||}\Big[\psi^R \mathfrak{A}_tv\big(\partial_t\eta^R\mathcal{P}_tf\big)\Big] dx\Big|dt
\\
&\lesssim\int_{R/2}^R\Vert\sq{b}\Vert_{BMO}\Big(\Vert\nabla_{||}(\psi^R \mathfrak{A}_tv)\Vert_{L^2(\Delta(0,2R))}\Vert \partial_t\eta^R\mathcal{P}_tf\Vert_{L^2(\mathbb{R}^n)}
\\
&\qquad\qquad+\Vert\psi^R \mathfrak{A}_tv\Vert_{L^2(\Delta(0,2R))}\Vert \partial_t\eta^R\nabla_{||}\mathcal{P}_tf\Vert_{L^2(\mathbb{R}^n)}\Big)dt.
\end{align*}
The $L^2$ boundedness \refprop{prop:L2boundednessOfSemigroupOperators} gives
\begin{align*}
|E_4^R|&\lesssim\int_{R/2}^R\Vert\nabla_{||}v\Vert_{L^2(\mathbb{R}^n)}\frac{\Vert\mathcal{P}_tf\Vert_{L^2(\mathbb{R}^n)}}{R}
+\Vert v\Vert_{L^2(\Delta(0,3R))}\frac{\Vert \nabla_{||}\mathcal{P}_tf\Vert_{L^2(\mathbb{R}^n)}}{R}dt
\\
&\lesssim\int_{R/2}^R\frac{1}{R}\Vert\nabla_{||}v\Vert_{L^{2}(\mathbb{R}^n)}\Vert f\Vert_{L^2(\mathbb{R}^n)} + R^{-\frac{n+2}{n+1}}\Vert v\Vert_{L^{2^*}(\Delta(0,3R))}\Vert f\Vert_{L^2(\mathbb{R}^n)}dt
\\
&\lesssim \frac{1}{\sqrt{R}}\Vert f\Vert_{L^2(\mathbb{R}^n)}\Vert v\Vert_{Y^{1,2}(\mathbb{R}^n\times(R/2,R))}.
\end{align*}
The right hand side converges to $0$ as $R$ approaches $\infty$.

\end{itemize}

This finishes the proof that all error terms $E_i^R$ converge to zero as $R$ goes to infinity, and hence the proof of \refthm{thm:MainThm2} is complete. It remains to prove \reflemma{lemma:DifferenceOperatorRisbounded} and \reflemma{lemma:CarlesonFctBoundedBMO} in the following subsection.

\subsection{Difference of averaging operator and Carleson function bound}

 We recall the definition of the smooth averaging operator $\mathfrak{A}_t$ in \eqref{eq:DefSmoothAveragingOp} and the rough averaging operator $\bar{\mathfrak{A}}_t$ in \eqref{eq:DefRoughAveragingOp}. Furthermore, we continue to use the notation \(\Omega^R:=\mathbb{R}^n\times(0,R)\) with cut-off functions \(\psi^R\) in $x$ and \(\eta^R\) in $t$ from the proof of \refthm{thm:MainThm2} in the previous subsection.
\begin{lemma}\label{lemma:DifferenceOperatorRisbounded}
    Let \(R>0\) and \(c>0\). Let \(\sq{b}\in BMO(\mathbb{R}^n)\) and $\mathcal{S}_t$ an operator. Then
    \begin{align*}\sq{J}_{111}&=\int_{\Omega t^R}\mathcal{S}_tf \cdot\Big[\mathcal{P}_{ct}^*(\mathrm{div}(\sq{b})\psi^R\mathfrak{A}_tv)-\mathcal{P}^*_{ct}(\mathrm{div}(\sq{b})\psi^R)\overline{\mathfrak{A}}_tv\Big] dxdt
    \\
    &\lesssim \Vert \mathcal{A}(t\mathcal{S}_tf)\Vert_{L^p}( \Vert \mathcal{S}(v)\Vert_{L^{p'}} + \Vert \mathcal{T}_{2+\EPS}(h)\Vert_{L^{p'}}).\end{align*}
The implicit constant is independent of $R>0$.
\end{lemma}

\begin{proof}
    First, let us fix \(t\) and consider a Whitney cube \(W\in\mathcal{D}(\Omega)\) such that \((x,t)\in W\). To clarify the dependence on \(t\) and \(x\), recall \(W_t:=\{(y,s)\in W|s=t\}\) and \(\phi^{W_t}(x)=\phi_{W}(x,t)\).
    
    Then
    \begin{align*}
        \Big(\mathcal{P}_{ct}^*(\mathrm{div}(\sq{b})\psi^R\mathfrak{A}_tv)-\mathcal{P}^*_{ct}(\mathrm{div}(\sq{b})\psi^R)\overline{\mathfrak{A}}_tv\Big)(x)&=\mathcal{P}_{ct}^*\Big(\mathrm{div}(\sq{b})\psi^R\big[\sum_{W'\in\mathcal{D}(\Omega)}(v)_{W'_t}\phi^{W'_t}-(v)_{W_t}\big]\Big)(x)
        \\
        &=\sum_{W'\in\mathcal{D}(\Omega)}\mathcal{P}_{ct}^*\Big(\mathrm{div}(\sq{b})\psi^R\big[(v)_{W'_t}-(v)_{W_t}\big]\phi^{W'_t}\Big)(x),
    \end{align*}
    where we used that \(\sum_{W'}\phi^{W'_t}\equiv \chi_{\mathbb{R}^n\times\{t\}}\).

    We further rewrite
    \begin{align*}
        &\sum_{W'\in\mathcal{D}(\Omega)}\mathcal{P}_{ct}^*\Big(\mathrm{div}(\sq{b})\psi^R\big[(v)_{W'_t}-(v)_{W_t}\big]\phi^{W'_t}\Big)(x)
        \\
        &=\sum_{W'\in\mathcal{D}(\Omega)}\mathcal{P}_{ct}^*\Big(\mathrm{div}(\sq{b}-(\sq{b})_{W'})\psi^R\big[(v)_{W'_t}-(v)_{W_t}\big]\phi^{W'_t}\Big)(x)
        \\
        &=\sum_{W'\in\mathcal{D}(\Omega)}\mathcal{P}_{ct}^*\Big(\mathrm{div}\Big((\sq{b}-(\sq{b})_{W'})\psi^R\big[(v)_{W'_t}-(v)_{W_t}\big]\phi^{W'_t}\Big)\Big)(x)
        \\
        &\qquad-\sum_{W'\in\mathcal{D}(\Omega)}\mathcal{P}_{ct}^*\Big((\sq{b}-(\sq{b})_{W'})\nabla\big[((v)_{W'_t}-(v)_{W_t})\phi^{W'_t}\big]\Big)(x)
        \\
        &\qquad-\sum_{W'\in\mathcal{D}(\Omega)}\mathcal{P}_{ct}^*\Big((\sq{b}-(\sq{b})_{W'})\nabla\psi^R\big[((v)_{W'_t}-(v)_{W_t})\phi^{W'_t}\big]\Big)(x)
        \\
        &=:I_t+II_t+III_t.
    \end{align*}

    With the off-diagonal estimates, \refcor{cor:Off-diagonalEstimatesForPdiv}, we now get
    \begin{align*}
        \Vert I_t\Vert_{L^2(W)}&\lesssim \frac{1}{l(W)}\sum_{k\geq 0}e^{-c2^{k}}\sum_{\substack{W'\in \mathcal{D}(\Omega),\\ \mathrm{dist}(W,W')\approx 2^kt}}\Big(\int_{l(W)/4}^{2l(W)}\Vert (\sq{b}-(\sq{b})_{W'})\big[(v)_{W'_t}-(v)_{W_t}\big]\phi^{W'_t}\Vert_{L^2(\frac{3}{2}W'_t)}^2dt\Big)^{1/2}
        \\
        &\lesssim \frac{1}{l(W)}\sum_{k\geq 0}e^{-c2^{k}t}|W|^{1/2}\sum_{\substack{W'\in \mathcal{D}(\Omega),\\ \mathrm{dist}(W,W')\approx 2^k}}\Big(\fint_{l(W)/4}^{2l(W)}\big(\fint_{\frac{3}{2}W'_t} |\sq{b}-(\sq{b})_{W'}|^rdx\big)^{2/r}
        \\
        &\hspace{50mm}\cdot\big(\fint_{\frac{3}{2}W'_t} |(v)_{W'_t}-(v)_{W_t}|^{2+\EPS}\big)^{2/(2+\EPS)}dt\Big)^{1/2}.
    \end{align*}
    Here in the last line, $\EPS>0$ is chosen sufficiently small such that \reflemma{lemma:HigherIntegrabilityForGradient} can be applied and $r$ is chosen such that $\frac{1}{2}=\frac{1}{r}+\frac{1}{2+\EPS}$.
    Now we use Poincar\'{e}'s and Jensen's inequality to obtain
    \begin{align*}
        \Vert I_t\Vert_{L^2(W)}&\lesssim \sum_{k\geq 0}ke^{-c2^{k}}|W|^{1/2}\sum_{W'\in \mathcal{R}_{2^k}(W), }\big(\fint_{\frac{3}{2}W'} |\nabla v|^{2+\EPS}dxdt\big)^{1/(2+\EPS)},
    \end{align*}
    where the set \(\mathcal{R}_{2^k}(W)\) denotes a covering of a fat annulus centered at \(W\) given by
    \begin{align*}
        &\mathcal{R}_{2^k}(W):=\{W'\in\mathcal{D}(\Omega)|\frac{1}{4}l(W')\leq l(W)\leq 4l(W'), \mathrm{dist}(W,W')\approx 2^k\}
        \\
        \textrm{and}\qquad &\mathsf{R}_{2^k(W)}=\bigcup_{W'\in\mathcal{R}_{2^k}}\frac{7}{4}W'.
    \end{align*}
   Next, we can apply \reflemma{lemma:HigherIntegrabilityForGradient} to get
    \begin{align*}
        &\lesssim \sum_{k\geq 0}ke^{-c2^{k}}|W|^{1/2}\sum_{W'\in\mathcal{R}_{2^k}(W)}\Big[\big(\fint_{\frac{7}{4}W'} |\nabla v|^{2}\big)^{1/2} + \big(\fint_{\frac{7}{4}W'} |h|^{2+\EPS}\big)^{1/(2+\EPS)}\Big]
        \\
        &\lesssim \sum_{k\geq 0}k2^{nk}e^{-c2^{k}}|W|^{1/2}\Big[\big(\int_{\mathsf{R}_{2^k}(W)} \frac{|\nabla v|^{2}}{s^{n+1}}dyds\big)^{1/2}+\sum_{W'\in\mathcal{R}_{2^k}(W)}\big(\fint_{\frac{7}{4}W'} |h|^{2+\EPS}\big)^{1/(2+\EPS)}\Big].
    \end{align*}
    With a similar argument we obtain for \(II_t\)
    \[\Vert II_t\Vert_{L^2(W)}^2\lesssim \sum_{k\geq 0}k2^{nk}e^{-c2^{k}}|W|^{1/2}\Big[\big(\int_{\mathsf{R}_{2^k}(W)} \frac{|\nabla v|^{2}}{s^{n+1}}dyds\big)^{1/2}+\sum_{W'\in\mathcal{R}_{2^k}(W)}\big(\fint_{\frac{7}{4}W'} |h|^{2+\EPS}\big)^{1/(2+\EPS)}\Big],\]
    and for \(III_t\)
    \[\Vert III_t\Vert_{L^2(W)}^2\lesssim \sum_{k\geq 0}\frac{l(W)}{R}k2^{nk}e^{-c2^{k}}|W|^{1/2}\Big[\big(\int_{\mathsf{R}_{2^k}(W)} \frac{|\nabla v|^{2}}{s^{n+1}}dyds\big)^{1/2}+\sum_{W'\in\mathcal{R}_{2^k}(W)}\big(\fint_{\frac{7}{4}W'} |h|^{2+\EPS}\big)^{1/(2+\EPS)}\Big].\]
    Since $t\in (0,R)$, the additional factor of $\frac{l(W)}{R}$ is bounded by $1$, all the right hand sides of the bounds for \(I_t,II_t\) and \(III_t\) agree.

    Together, we have now
    \begin{align*}
        \sq{J}_{111}&=\int_\Omega t\mathcal{S}_tf \cdot\Big[\mathcal{P}_{ct}^*(\mathrm{div}(\sq{b})\mathfrak{A}v)-\mathcal{P}^*_{ct}(\mathrm{div}(\sq{b}))\overline{\mathfrak{A}}v\Big] dxdt
        \\
        &\lesssim \sum_{W\in\mathcal{D}}|W|^{1/2}l(W)\big(\fint_W|\partial_t^2\mathcal{P}_tf|^2dxdt\big)^{1/2}\big(\int_W|I+II|^2dxdt\big)^{1/2}
        \\
        &\lesssim \sum_{k\geq 0}k2^{kn}e^{-c2^{k}}\sum_{W\in\mathcal{D}}|W|l(W)\big(\fint_W|\partial_t^2\mathcal{P}_tf|^2dxdt\big)^{1/2}
        \\
        &\qquad\cdot\Big[\big(\int_{\mathsf{R}_{2^k}(W)} \frac{|\nabla v|^{2}}{s^{n+1}}dyds\big)^{1/2}+\sum_{W'\in\mathcal{R}_{2^k}(W)}\big(\fint_{\frac{7}{4}W'} |h|^{2+\EPS}\big)^{1/(2+\EPS)}\Big]
    \end{align*}

    If we fix \(k\), we can follows the same proof idea of \reflemma{lemma:dyadicToNondyadic}, and obtain like in \eqref{eq:EqInProofDyadicToNondyadic} that
    \begin{align*}
        &\sum_{W\in\mathcal{D}}|W|l(W)\big(\fint_W|\mathcal{S}_tf|^2dxdt\big)^{1/2}\Big[\big(\int_{\mathsf{R}_{2^k}(W)} \frac{|\nabla v|^{2}}{s^{n+1}}dyds\big)^{1/2}+\sum_{W'\in\mathcal{R}_{2^k}(W)}\big(\fint_{\frac{7}{4}W'} |h|^{2+\EPS}\big)^{1/(2+\EPS)}\Big]
        \\
        &\lesssim \int_{\partial\Omega}\Big(\sum_{W:W\cap\Gamma_{2^{k+2}}(x)\neq \emptyset}\int_{W}\frac{\Big(\fint_{W}|\partial_t^2\mathcal{P}_tf|^2\Big)^{1/2}\Big[\big(\fint_{W} |\nabla v|^{2}dyds\big)^{1/2}+\big(\fint_{\frac{7}{4}W} |h|^{2+\EPS}\big)^{1/(2+\EPS)}\Big]}{\delta(y)^{n-1}}dy\Big)d\sigma(x),
    \end{align*}
    where the only difference lies in the enlarged aperture. From here on, we can continue as in the proof of \reflemma{lemma:dyadicToNondyadic}. Due to this enlarged aperture, we obtain the additional factor of \(2^{kn}\), and hence in total
    \begin{align*}
        &\sum_{k\geq 0}k2^{kn}e^{-c2^{k}}\sum_{W\in\mathcal{D}}|W|l(W)\big(\fint_W|\mathcal{S}_tf|^2dxdt\big)^{1/2}
        \\
        &\qquad\qquad\cdot\Big[\big(\int_{\mathsf{R}_{2^k}(W)} \frac{|\nabla v|^{2}}{s^{n+1}}dyds\big)^{1/2}+\sum_{W'\in\mathcal{R}_{2^k}(W)}\big(\fint_{\frac{7}{4}W'} |h|^{2+\EPS}\big)^{1/(2+\EPS)}\Big]
        \\
        &\lesssim \sum_{k\geq 0}k2^{2kn}e^{-c2^{k}}\Vert\mathcal{A}_{2^{k+2}}(t\mathcal{S}_tf)\Vert_{L^p}\big(\Vert\mathcal{A}_{2^{k+2}}(t|\nabla v|)\Vert_{L^{p'}}+\Vert \mathcal{T}_{2^{k+2}, 2+\EPS}(h)\Vert_{L^{p'}}\big)
    \end{align*}
    From \cite{Auscher_change_2011} it follows that there is \(c=c(n,p)\) such \(\Vert\mathcal{A}_{2^{k+2}}(f)\Vert_{L^{p'}}\lesssim 2^{ck}\Vert\mathcal{A}(t|\nabla v|)\Vert_{L^{p'}}\). Similarly, \(\Vert\mathcal{A}_{2^{k+2}}(f)\Vert_{L^{p'}}\lesssim 2^{ck}\Vert\mathcal{A}(t|\nabla v|)\Vert_{L^{p'}}\) holds (see \refcor{cor:appendixDiffApCompT}). Performing this change of angle on all involved area functions and using \eqref{eq:CarlAreaTcomparable} yields
    \begin{align*}
        &\sum_{k\geq 0}k2^{2kn}e^{-c2^{k}}\Vert\mathcal{A}_{2^{k+2}}(t\mathcal{S}_tf)\Vert_{L^p}\big(\Vert\mathcal{A}_{2^{k+2}}(t|\nabla v|)\Vert_{L^{p'}}+\Vert \mathcal{T}_{2^{k+2}, 2+\EPS}(h)\Vert_{L^{p'}}\big)
        \\
        &\qquad\lesssim \sum_{k\geq 0}k2^{k(2n+c)}e^{-c2^{k}}\Vert\mathcal{A}(t\mathcal{S}_tf)\Vert_{L^p}\big(\Vert\mathcal{A}(t|\nabla v|)\Vert_{L^{p'}}+\Vert \mathcal{T}_{2+\EPS}(h)\Vert_{L^{p'}}\big)
        \\
        &\qquad \lesssim\Vert\mathcal{A}(t\mathcal{S}_tf)\Vert_{L^p}\big(\Vert\mathcal{A}(t|\nabla v|)\Vert_{L^{p'}}+\Vert \mathcal{T}_{2+\EPS}(h)\Vert_{L^{p'}}\big).
    \end{align*}

\end{proof}

\begin{lemma}\label{lemma:CarlesonFctBoundedBMO}
    Let \(\sq{b}\in BMO(\mathbb{R}^n)\) and $c>0$. Then
    \[\Vert\mathcal{C}(t\eta^R\mathcal{P}_{ct}(\mathrm{div}(\sq{b})\psi^R))\Vert_{L^\infty(\mathbb{R}^n)}^2=\sup_{\Delta\subset \partial\Omega}\frac{1}{\sigma(\Delta)}\int_{T(\Delta)}|\eta^R\mathcal{P}_{ct}(\mathrm{div}(\sq{b})\psi^R)|^2t dxdt\lesssim \Vert \sq{b}\Vert_{BMO}^2.\]
\end{lemma}

\begin{proof}
To bound this expression, we split into a local and an away part by fixing a boundary ball \(\Delta\subset \partial\Omega\), and discussing
\begin{align*}
    &\int_{T(\Delta)}|\eta^R\mathcal{P}_{ct}(\mathrm{div}(\sq{b})\psi^R)|^2t dxdt=\int_{T(\Delta)}|\eta^R\mathcal{P}_{ct}(\mathrm{div}[(\sq{b}-(\sq{b})_{3\Delta})]\psi^R)|^2t dxdt
    \\
    &\qquad\leq\int_{T(\Delta)}|\eta^R\mathcal{P}_{ct}\mathrm{div}[(\sq{b}-(\sq{b})_{3\Delta})\psi^R]|^2t dxdt + \int_{T(\Delta)}|\eta^R\mathcal{P}_{ct}[(\sq{b}-(\sq{b})_{3\Delta})\nabla\psi_R]|^2t dxdt
    \\
    &\qquad\leq\int_{T(\Delta)}|\mathcal{P}_{ct}\mathrm{div}[(\sq{b}-(\sq{b})_{3\Delta})\psi^R\chi_{3\Delta}]|^2t dxdt + \int_{T(\Delta)}|\mathcal{P}_{ct}\mathrm{div}[(\sq{b}-(\sq{b})_{3\Delta})\psi^R\chi_{\mathbb{R}^n\setminus 3\Delta}]|^2t dxdt
    \\
    &\qquad\qquad+\int_{T(\Delta)}|\eta^R\mathcal{P}_{ct}[(\sq{b}-(\sq{b})_{3\Delta})\nabla \psi^R\chi_{3\Delta}]|^2t dxdt+ \int_{T(\Delta)}|\eta^R\mathcal{P}_{ct}[(\sq{b}-(\sq{b})_{3\Delta})\nabla \psi^R\chi_{\mathbb{R}^n\setminus 3\Delta}]|^2t dxdt
    \\
    &=I_1+I_2+I_3+I_4.
\end{align*}

For the local part \(I_1\) we use the Kato results, \refprop{prop:Kato}, and obtain
\begin{align*}
    I_1&=\int_0^{l(\Delta)}\Big(\int_{\mathbb{R}^n}|\mathcal{P}_{ct}\mathrm{div}[(\sq{b}-(\sq{b})_{3\Delta})\psi^R\chi_{3\Delta}]|^2 dx\Big) tdt
    \\
    &=\int_0^{l(\Delta)}\Big(\int_{\mathbb{R}^n}|L_{||}^{1/2}\mathcal{P}_{ct}L^{-1/2}_{||}\mathrm{div}[(\sq{b}-(\sq{b})_{3\Delta})\psi^R\chi_{3\Delta}]|^2 dx\Big) tdt
    \\
    &=\int_0^{l(\Delta)}\Big(\int_{\mathbb{R}^n}|\nabla_{||}\mathcal{P}_{ct}L^{-1/2}_{||}\mathrm{div}[(\sq{b}-(\sq{b})_{3\Delta})\psi^R\chi_{3\Delta}]|^2 dx\Big) tdt
    \\
    &=\int_0^{l(\Delta)}\Big(\int_{\mathbb{R}^n}A_{||}\nabla_{||}\mathcal{P}_{ct}L^{-1/2}_{||}\mathrm{div}[(\sq{b}-(\sq{b})_{3\Delta})\psi^R\chi_{3\Delta}] \cdot \nabla_{||}\mathcal{P}_{ct}L^{-1/2}_{||}\mathrm{div}[(\sq{b}-(\sq{b})_{3\Delta})\psi^R\chi_{3\Delta}] dx\Big) tdt
    \\
    &=\frac{1}{2c^2}\int_0^{l(\Delta)}\Big(\int_{\mathbb{R}^n}\partial_t\mathcal{P}_{ct}L^{-1/2}_{||}\mathrm{div}[(\sq{b}-(\sq{b})_{3\Delta})\psi^R\chi_{3\Delta}] \cdot \mathcal{P}_{ct}L^{-1/2}_{||}\mathrm{div}[(\sq{b}-(\sq{b})_{3\Delta})\psi^R\chi_{3\Delta}] dx\Big) dt
    \\
    &=\frac{1}{4c^2}\int_0^{l(\Delta)}\partial_t\Vert\mathcal{P}_{ct}L^{-1/2}_{||}\mathrm{div}[(\sq{b}-(\sq{b})_{3\Delta})\psi^R\chi_{3\Delta}]\Vert_{L^2(\mathbb{R}^n)}^2 dt
    \\
    &=\Vert\mathcal{P}_{cl(\Delta)}L^{-1/2}_{||}\mathrm{div}[(\sq{b}-(\sq{b})_{3\Delta})\psi^R\chi_{3\Delta}]\Vert_{L^2(\mathbb{R}^n)}^2
    - \Vert L^{-1/2}_{||}\mathrm{div}[(\sq{b}-(\sq{b})_{3\Delta})\psi^R\chi_{3\Delta}]\Vert_{L^2(\mathbb{R}^n)}^2
    \\
    &\lesssim \Vert \sq{b}-(\sq{b})_{3\Delta}\Vert_{L^2(3\Delta)}^2\leq \sigma(\Delta)\Vert \sq{b}\Vert_{BMO}^2.
\end{align*}

For the away part \(I_2\) we have by the off-diagonal estimates \refcor{cor:Off-diagonalEstimatesForPdiv} that
\begin{align*}
    I_2&\leq\int_0^{l(\Delta)}\Big(\sum_{k\geq 1}\Vert\mathcal{P}_{ct}\mathrm{div}[(\sq{b}-(\sq{b})_{3\Delta})\psi^R\chi_{2^{k+1}\Delta\setminus2^k\Delta}]\Vert_{L^2(\Delta)}\Big)^2 tdt
    \\
    &\lesssim\int_0^{l(\Delta)}\Big(\sum_{k\geq 1}\frac{1}{t}e^{-\frac{2^{2k}l(\Delta)^2}{t^2}}\Vert(\sq{b}-(\sq{b})_{3\Delta})\Vert_{L^2(2^{k+1}\Delta\setminus2^k\Delta)}\Big)^2 tdt.
\end{align*}
We can estimate
\begin{align}
    &\Vert \sq{b}-(\sq{b})_{3\Delta} \Vert_{L^2(2^{k+1}\Delta\setminus2^k\Delta)}\nonumber
    \\
    &\leq\Vert \sq{b}-(\sq{b})_{2^{k+1}\Delta} \Vert_{L^2(2^{k+1}\Delta)} +\sum_{l=1}^{k-1} \Vert (\sq{b})_{2^{l+1}\Delta}-(\sq{b})_{2^l\Delta} \Vert_{L^2(2^{k+1}\Delta\setminus2^k\Delta)}\leq k2^{kn/2}\sigma(\Delta)^{1/2}\Vert \sq{b}\Vert_{BMO}.\label{eq:BoundBybBMO}
\end{align}
Furthermore, since the function \((0,l(\Delta)]\to\mathbb{R}, t\mapsto\frac{1}{t}e^{-\frac{2^{2k}l(\Delta)^2}{t^2}}\) takes on its maximum in \(t=l(\Delta)\), we obtain
\begin{align*}
    I_2\lesssim &\sigma(\Delta)\Vert \sq{b}\Vert_{BMO}^2\int_0^{l(\Delta)}\Big(\sum_{k\geq 1}\frac{k2^{kn/2}}{2^kl(\Delta)}\Big)^2 tdt\lesssim \sigma(\Delta)\Vert \sq{b}\Vert_{BMO}^2\int_0^{l(\Delta)}\frac{t}{l(\Delta)^2}dt\lesssim \sigma(\Delta)\Vert \sq{b}\Vert_{BMO}^2.
\end{align*}

For \(I_3\), we denote \(m:=\min\{l(\Delta),R\}\), and use \refprop{prop:L2boundednessOfSemigroupOperators} to get
\begin{align*}
    &\int_{T(\Delta)}|\eta^R\mathcal{P}_{ct}[(\sq{b}-(\sq{b})_{3\Delta})\nabla \psi^R\chi_{3\Delta}]|^2t dxdt\leq \int_0^m\Vert\mathcal{P}_{ct}[(\sq{b}-(\sq{b})_{3\Delta})\nabla \psi^R\chi_{3\Delta}]\Vert_{L^2(\mathbb{R}^n)}^2dt
    \\
    &\qquad\lesssim\int_0^m\Vert(\sq{b}-(\sq{b})_{3\Delta})\nabla \psi^R\chi_{3\Delta}\Vert_{L^2(\mathbb{R}^n)}^2dt \lesssim\int_0^m\frac{1}{R}\Vert\sq{b}-(\sq{b})_{3\Delta}\Vert_{L^2(3\Delta)}^2dt\lesssim \sigma(\Delta)\Vert \sq{b}\Vert_{BMO}^2\int_0^m\frac{1}{R}dt
    \\
    &\qquad\leq \sigma(\Delta)\Vert \sq{b}\Vert_{BMO}^2.
\end{align*}

The integral \(I_4\) can be dealt with similar to \(I_2\). First, by the off-diagonal estimates \refprop{prop:Off-diagonalEstimates},
\begin{align*}
    I_4&\leq\int_0^{l(\Delta)}\eta^R\Big(\sum_{k\geq 1}\Vert\mathcal{P}_{ct}[(\sq{b}-(\sq{b})_{3\Delta})\nabla\psi^R\chi_{2^{k+1}\Delta\setminus2^k\Delta}]\Vert_{L^2(\Delta)}\Big)^2 tdt
    \\
    &\lesssim\int_0^{\min\{l(\Delta),R\}}\Big(\sum_{k\geq 1}e^{-\frac{2^{2k}l(\Delta)^2}{t^2}}\Vert(\sq{b}-(\sq{b})_{3\Delta})\nabla\psi^R\Vert_{L^2(2^{k+1}\Delta\setminus2^k\Delta)}\Big)^2 tdt
    \\
    &\lesssim\int_0^{\min\{l(\Delta),R\}}\Big(\sum_{k\geq 1}\frac{1}{R}e^{-\frac{2^{2k}l(\Delta)^2}{t^2}}\Vert(\sq{b}-(\sq{b})_{3\Delta})\Vert_{L^2(2^{k+1}\Delta\setminus2^k\Delta)}\Big)^2 tdt.
\end{align*}
Recall that we denote \(m:=\min\{l(\Delta),R\}\), and use \eqref{eq:BoundBybBMO}, we get
\begin{align*}
    &\int_0^{m}\Big(\sum_{k\geq 1}\frac{1}{R}e^{-\frac{2^{2k}l(\Delta)^2}{t^2}}\Vert(\sq{b}-(\sq{b})_{3\Delta})\Vert_{L^2(2^{k+1}\Delta\setminus2^k\Delta)}\Big)^2 tdt
    \\
    &\lesssim \sigma(\Delta)\int_0^{m}\Big(\sum_{k\geq 1}\frac{k}{m}e^{-\frac{2^{2k}l(\Delta)^2}{m^2}}\Vert\sq{b}\Vert_{BMO}\Big)^2 tdt
\end{align*}
Since the function \(m\mapsto\frac{1}{m}e^{-\frac{2^{2k}l(\Delta)^2}{m^2}}\) takes on its maximum in \(m=2^kl(\Delta)\), we obtain
\begin{align*}
    I_4\lesssim &\sigma(\Delta)\Vert \sq{b}\Vert_{BMO}^2\int_0^{l(\Delta)}\Big(\sum_{k\geq 1}\frac{k}{2^kl(\Delta)}\Big)^2 tdt\lesssim \sigma(\Delta)\Vert \sq{b}\Vert_{BMO}^2\int_0^{l(\Delta)}\frac{t}{l(\Delta)^2}dt\lesssim \sigma(\Delta)\Vert \sq{b}\Vert_{BMO}^2.
\end{align*}

\end{proof}

\appendix

\section{Conical Area function bounds}\label{sec:AreaFctBounds}

We would like to prove the area function bounds in \refprop{prop:AreaFctBoundsInLp}. The strategy of proof is as in the case for operators with bounded matrices \(A\). For each operator \(T_t\in \{\partial_t\mathcal{P}_t,t\partial_{t}^2\mathcal{P}_t,\nabla_{||}\partial_t\mathcal{P}_t\}\), we first prove that the Carleson function is bounded by the \(L^\infty\)-norm of the gradient of \(f\) (see Section~\ref{appendix:SubSecCarlesonBounds}). Then, we show that the \(L^1\)-norm of the area function is bounded by the Hardy-Sobolev norm of \(f\) (see Section~\ref{appendix:SubsectionHardyBound}, before we apply real interpolation to obtain the \(L^p\) area function bounds (see Section~\ref{appendix:Subsection:Interpolation}).

Before we are able to proceed, we are going to need the following preliminary statements.
\begin{prop}[Prop 4.3 in \cite{hofmann_lp_2022}]\label{Prop:KernelBounds}
    Let us denote the kernel of the semigroup \(e^{tL_{||}}\) by \(K_t(x,y)\).
    Then, for any \(l\in\mathbb{N}\) there exists \(C=C(n,\lambda,l), \beta=\beta(n,\lambda)>0\) such that
    \begin{align}
        |\partial_t^l K_t(x,y)|\leq C_l t^{-\frac{n}{2}-l}e^{-\beta\frac{|x-y|^2}{t}}
    \end{align}
    for all \(x,y\in\mathbb{R}^n\).
\end{prop}

Furthermore, let us observe the following Caccioppoli type inequalities for the operator \(\mathcal{P}_t\).

\begin{lemma}\label{lemma:CaccTypeInequalityForPSimple}
Let $B\subset\Omega$ be a ball with diameter $R$ such that \(\mathrm{dist}(B,\partial\Omega)=3R\). Then
\[\int_{B}|\nabla_{||}\mathcal{P}_tf|^2dxdt\lesssim \frac{1}{R^2}\int_{2B}|\mathcal{P}_tf|^2dxdt.\]
\end{lemma}

\begin{lemma}\label{lemma:CaccTypeInequalityForP}
Let $B\subset\Omega$ be a ball with diameter $R$ such that \(\mathrm{dist}(B,\partial\Omega)=3R\). Then
\[\int_{B}|\nabla_{||} \partial_t\mathcal{P}_tf|^2dxdt=\int_{B}|2t\nabla_{||} L_{||}e^{-t^2L_{||}}f|^2dxdt\lesssim \frac{1}{R^2}\int_{2B}|2tL_{||}e^{-t^2L_{||}}f|^2dxdt=\frac{1}{R^2}\int_{2B}|\partial_t\mathcal{P}_tf|^2dxdt.\]
\end{lemma}

Since the proof of \reflemma{lemma:CaccTypeInequalityForPSimple} and \reflemma{lemma:CaccTypeInequalityForP} are analogous, we only provide the proof of \reflemma{lemma:CaccTypeInequalityForP}. 

\begin{proof}[Proof of \reflemma{lemma:CaccTypeInequalityForP}]
We observe that \(v(x,t):=\partial_t\mathcal{P}_tf(x)=\partial_t e^{-t^2L_{||}}f(x)\) satisfies the PDE
\begin{align*}
\partial_t \partial_t\mathcal{P}_tf&=\partial_t\big(-2tL_{||}e^{-t^2L_{||}}f\big)= -2L_{||} e^{-t^2L_{||}}f - 2tL_{||} \partial_te^{-t^2L_{||}}f=\frac{1}{t}\partial_t\mathcal{P}_t f - 2tL_{||} \partial_t\mathcal{P}_tf.
\end{align*}

Analogously to the standard proof for Caccioppoli's inequality we take a smooth cut-off function \(\psi\in C^\infty(\Omega)\) that satisfies 
\[\psi\equiv 1 \textrm{ on } B,\]
and
\[\psi\equiv 0 \textrm{ on } \Omega\setminus \frac{3}{2}B,\]
and
\[\Vert \nabla \psi \Vert \lesssim \frac{1}{l(W)}\approx \frac{1}{R}.\]
Thus, ellipticity of $A$ implies
\begin{align*}
    \int_{B}|\nabla_{||} \partial_t\mathcal{P}_tf|^2dxdt&\lesssim \int_{\frac{3}{2}B}|\nabla_{||} \partial_t\mathcal{P}_tf|^2\psi^2dxdt\lesssim \int_{\frac{3}{2}B} A\nabla_{||} \partial_t\mathcal{P}_tf \cdot \nabla_{||} \partial_t\mathcal{P}_tf \psi^2dxdt
    \\
    &= -\int_{\frac{3}{2}B} (L_{||} \partial_t\mathcal{P}_tf)(    \partial_t\mathcal{P}_tf\psi^2)dxdt -\int_{\frac{3}{2}B} A\nabla_{||}  \partial_t\mathcal{P}_tf \cdot \nabla_{||}\psi^2  \partial_t\mathcal{P}_tf dxdt
    \\
    &=:I_1+I_2.
\end{align*}
Since $\frac{1}{t}\partial_t\mathcal{P}_tf\psi^2\in W^{1,2}_0(\frac{3}{2}B)$, we can apply the PDE and obtain for the first term
\begin{align*}
    I_1&=\frac{1}{2}\int_{\frac{3}{2}B} \frac{\partial_t^2\mathcal{P}_tf}{t}\partial_t\mathcal{P}_tf\psi^2dxdt-\frac{1}{2}\int_{\frac{3}{2}B} \frac{\partial_t\mathcal{P}_tf}{t}\partial_t\mathcal{P}_tf\psi^2dxdt=:I_{11}+I_{12}.
\end{align*}
While \(I_{12}\) is one of the terms we expect on the right hand side, we get for the other term
\begin{align*}
    I_{11}=\frac{1}{2}\int_{\frac{3}{2}B} \partial_t(\partial_t\mathcal{P}_tf)^2\frac{\psi^2}{t}dxdt
    \\
    &=\frac{1}{2}\int_{\frac{3}{2}B} \partial_t\big((\partial_t\mathcal{P}_tf)^2\frac{\psi^2}{t}\big)dxdt-\frac{1}{2}\int_{\frac{3}{2}B} (\partial_t\mathcal{P}_tf)^2\big(\frac{\partial_t\psi^2}{t}-\frac{\psi^2}{t^2}\big)dxdt.
\end{align*}
The first term vanishes, and hence
\[|I_1|\leq |I_{11}|+|I_{12}|\lesssim \frac{1}{R^2}\int_{\frac{3}{2}B}|\partial_t\mathcal{P}_tf|^2dxdt.\]
For the second term, we use \refprop{prop:ContantAntiSymmMatrix} and Hölder's inequality to get
\begin{align*}
I_2&=\int_{\frac{3}{2}B} A\nabla_{||}  (\partial_t\mathcal{P}_tf)^2 \cdot \nabla_{||}\psi^2 dxdt=\frac{1}{2}\int_{\frac{3}{2}B} (A-(A^a)_{\frac{3}{2}B})\nabla_{||}  (\partial_t\mathcal{P}_tf)^2 \cdot \nabla_{||}\psi^2 dxdt
\\
&\lesssim R^{n+1}\Big(\fint_{\frac{3}{2}B}|A-(A^a)_{\frac{3}{2}B}|^4\Big)^{1/4}\Big(\fint_{\frac{3}{2}B}|\nabla_{||} \partial_t\mathcal{P}_tf|^2 \psi^2dxdt\Big)\Big(\fint_{\frac{3}{2}B}|\nabla_{||} \psi|^4 |\partial_t\mathcal{P}_tf|^{4}dxdt\Big)^{1/4}.
\end{align*}
Next, we apply \cite[Proposition 6]{hofmann_dirichlet_2022}, which states that
\[\sup_{\frac{3}{2}B}|\partial_t\mathcal{P}_tf|\lesssim\Big(\fint_{2B} |\partial_t\mathcal{P}_tf|^{2}dxdt\Big)^{1/2},\]
which yields
\begin{align*}
    I_2\lesssim R^{n+1}\Big(\fint_{\frac{3}{2}B}|A-(A^a)_{\frac{3}{2}B}|^4\Big)^{1/4}\Big(\fint_{\frac{3}{2}B}|\nabla \partial_t\mathcal{P}_tf|^2 \psi^2dxdt\Big)\Big(\fint_{2B} \frac{|\partial_t\mathcal{P}_tf|^{2}}{t^2}dxdt\Big)^{1/2}.
\end{align*}
Using that $A\in BMO(\Omega)$, and Young's inequality for a small $\sigma>0$, leads to
\begin{align*}
I_2\lesssim \sigma\int_{\frac{3}{2}B}|\nabla \partial_t\mathcal{P}_tf|^2 \psi^2dxdt+C_\sigma\int_{2B} \frac{|\partial_t\mathcal{P}_tf|^{2}}{t^2}dxdt.
\end{align*}
For a sufficiently small choice of $\sigma>0$, absorbing the first term in the left hand side completes the proof.

\end{proof}

We note that \cite{hofmann_lp_2022} establishes the following $L^2$ bounds for the vertical square function:

\begin{prop}\label{prop:VerticalAreaFunctionL2Bounds}
    Let $T_t\in \{\partial_t\mathcal{P}_t, t\nabla\partial_t\mathcal{P}_t, t\partial_t^2\mathcal{P}_t\}$. Then for all $f\in W^{1,2}(\mathbb{R}^n)$
    \[\Big\Vert\Big(\int_0^\infty |T_tf|^2\frac{dt}{t}\Big)^{1/2}\Big\Vert_{L^2(\mathbb{R}^n)}\leq C\Vert\nabla_{||}f\Vert_{L^2(\mathbb{R}^n)}.\]
\end{prop}
We note that we use the conical area function throughout the proofs of \refthm{thm:MainThm1} and \refthm{thm:MainThm2}. For the \(L^2\) norms, \cite[Lemma 3.14]{milakis_harmonic_2013} implies the relationship
\[\Big\Vert\Big(\int_0^\infty |T_tf|^2\frac{dt}{t}\Big)^{1/2}\Big\Vert_{L^2(\mathbb{R}^n)}\approx\Vert\mathcal{A}(T_tf)\Vert_{L^2(\mathbb{R}^n)},\]
whence their results can be used for our conical area functions.

\subsection{Carleson function bounds}\label{appendix:SubSecCarlesonBounds}

\begin{lemma}\label{lemma:CarlesonFunctionBound}
    Let $T_t\in \{\partial_t\mathcal{P}_t, t\partial_t^2\mathcal{P}_t\}$ and \(f\in W^{1,\infty}(\mathbb{R}^n)\), then
     \begin{align}\sup_{\Delta\subset\partial\Omega}\frac{1}{\sigma(\Delta)}\int_{T(\Delta)}\frac{|T_tf(x,t)|^2}{t}dxdt\lesssim \Vert\nabla_{||}f\Vert_{L^\infty(\mathbb{R}^n)}^2.\label{eq:CarlesonFctBound}
     \end{align}
\end{lemma}

\begin{proof}
    Let us fix a boundary ball $\Delta\subset\partial\Omega$. We introduce a smooth cut-off function \(\eta\in C_0^\infty(3\Delta)\) with \(\eta\equiv 1\) on \(2\Delta\) and \(|\nabla \eta|\lesssim \frac{1}{l(\Delta)}\) and split the data \(f\) into a local and far-away part. Since all our operators $T_t$ kill constants, we obtain

    \begin{align*}
        &\int_{T(\Delta)}\frac{|T_tf(x)|^2}{t}dxdt=\int_{T(\Delta)}\frac{|T_t(f-(f)_{3\Delta})(x)|^2}{t}dxdt
        \\
        &\qquad \leq \int_{T(\Delta)}\frac{|T_t((f-(f)_{3\Delta})\eta)(x)|^2}{t}dxdt
        \\
        &\qquad \qquad+ \int_{T(\Delta)}\frac{|T_t((f-(f)_{3\Delta})(1-\eta))(x)|^2}{t}dxdt.
    \end{align*}
    From the $L^2$ area function bounds in \refprop{prop:VerticalAreaFunctionL2Bounds}, we obtain boundedness of the first integral by 
    \[\Vert\nabla( (f-(f)_{3\Delta})\eta)\Vert_{L^2(\mathbb{R}^n)}\lesssim \frac{1}{l(\Delta)}\Vert (f-(f)_{3\Delta}) \Vert_{L^2(3\Delta)} + \Vert \nabla_{||} f\Vert_{L^2(3\Delta)}\lesssim \Vert \nabla_{||} f\Vert_{L^2(3\Delta)},\]
    where Poincar\'{e}'s inequality was applied. Thus, we only need to deal with the second integral, the far-away part.

    \smallskip
    We have for a fixed \(t\in (0,l(\Delta))\) with off-diagonal estimates \refprop{prop:Off-diagonalEstimates}
    \begin{align*}
        &\Vert T_t ((f-(f)_{3\Delta})(1-\eta))\Vert_{L^2(\Delta)}
        \\
        &\lesssim \Vert T_t (\chi_{4\Delta}(f-(f)_{3\Delta})(1-\eta))\Vert_{L^2(\Delta)}+\sum_{k\geq 2}\Vert T_t (\chi_{2^{k+1}\Delta\setminus 2^{k}\Delta}(f-(f)_{3\Delta}))\Vert_{L^2(\Delta)}
        \\
        &\lesssim \frac{1}{t}e^{-c\frac{l(\Delta)^2}{t^2}}\Vert f-(f)_{3\Delta} \Vert_{L^2(4\Delta\setminus 2\Delta)} +  \sum_{k\geq 2}\frac{1}{t}e^{-c\frac{2^{2k}l(\Delta)^2}{t^2}}\Vert f-(f)_{3\Delta}\Vert_{L^2(2^{k+1}\Delta\setminus 2^k\Delta)}
        \\
        &\lesssim \frac{1}{t}e^{-c\frac{l(\Delta)^2}{t^2}}\Vert f-(f)_{4\Delta} \Vert_{L^2(4\Delta\setminus 2\Delta)} + \frac{1}{t}e^{-c\frac{l(\Delta)^2}{t^2}}\Vert (f)_{4\Delta}-(f)_{3\Delta} \Vert_{L^2(4\Delta\setminus 2\Delta)}
        \\
        &\qquad +\sum_{k\geq 2}\frac{1}{t}e^{-c\frac{2^{2k}l(\Delta)^2}{t^2}}\Big(\Vert f-(f)_{2^{k+1}\Delta}\Vert_{L^2(2^{k+1}\Delta\setminus 2^k\Delta)} 
        \\
        &\qquad + \sum_{j=2}^{k-1} \Vert (f)_{2^{j+1}\Delta}-(f)_{2^j\Delta}\Vert_{L^2(2^{k+1}\Delta\setminus 2^k\Delta)} + \Vert (f)_{4\Delta}-(f)_{3\Delta}\Vert_{L^2(2^{k+1}\Delta\setminus 2^k\Delta)}\Big).
    \end{align*}

    For each term of the form \(\Vert (f)_{2^{j+1}\Delta}-(f)_{2^j\Delta}\Vert_{L^2(2^{k+1}\Delta\setminus 2^k\Delta)}\) or \(\Vert f-(f)_{4\Delta}\Vert_{L^2(4\Delta)}\) we apply Poincar\'{e}'s inequality and obtain
    \begin{align*}
        &\lesssim\frac{l(\Delta)}{t}e^{-c\frac{l(\Delta)^2}{t^2}}\Vert \nabla_{||}f \Vert_{L^2(4\Delta\setminus 2\Delta)} + \frac{l(\Delta)}{t}e^{-c\frac{l(\Delta)^2}{t^2}}\Vert\nabla_{||}f \Vert_{L^2(4\Delta\setminus 2\Delta)}
        \\
        &\qquad +\sum_{k\geq 2}\frac{1}{t}e^{-c\frac{2^{2k}l(\Delta)^2}{t^2}}\Big(2^{k+1}l(\Delta)\Vert \nabla_{||}f\Vert_{L^2(2^{k+1}\Delta)}
        \\
        &\qquad\qquad + \sum_{j=2}^{k-1} 2^{(k+1)\frac{n}{2}+j+1}l(\Delta)^{\frac{n}{2}+1}(|\nabla_{||} f|)_{2^{j+1}\Delta}
         + 2^{(k+1)\frac{n}{2}}l(\Delta)^{\frac{n}{2}+1} (|\nabla_{||} f|)_{4\Delta}\Big)
        \\
        &\lesssim \frac{l(\Delta)^{\frac{n}{2}+1}}{t}e^{-c\frac{l(\Delta)^2}{t^2}}\Vert \nabla_{||}f \Vert_{\infty} + \frac{l(\Delta)^{\frac{n}{2}+1}}{t}e^{-c\frac{l(\Delta)^2}{t^2}}\Vert\nabla_{||}f \Vert_{\infty}
        \\
        &\qquad +\sum_{k\geq 2}\frac{1}{t}e^{-c\frac{2^{2k}l(\Delta)^2}{t^2}}\Big(2^{(k+1)(\frac{n}{2}+1)}l(\Delta)^{\frac{n}{2}+1}\Vert \nabla_{||}f\Vert_{\infty}
        \\
        &\qquad \qquad+ \sum_{j=2}^{k-1} 2^{(k+1)(\frac{n}{2}) + j+1}l(\Delta)^{\frac{n}{2}+1}\Vert \nabla_{||} f\Vert_{\infty} + 2^{(k+1)\frac{n}{2}}l(\Delta)^{\frac{n}{2}+1} \Vert \nabla_{||} f\Vert_{\infty}\Big)
        \\
        &\lesssim \frac{l(\Delta)^{\frac{n}{2}+1}}{t}\Big(e^{-c\frac{l(\Delta)^2}{t^2}} + \sum_{k\geq 2} (k+1)2^{(k+1)(\frac{n}{2}+1)}e^{-c\frac{2^{2k}l(\Delta)^2}{t^2}}\Big)\Vert \nabla_{||} f\Vert_{\infty}
        \\
        &\lesssim \frac{l(\Delta)^{\frac{n}{2}+1}}{t}\Big(\sum_{k\geq 1} (k+1)2^{(k+1)(\frac{n}{2}+1)}e^{-c\frac{2^{2k}l(\Delta)^2}{t^2}}\Big)\Vert \nabla_{||} f\Vert_{\infty}
        \\
        &\lesssim \frac{l(\Delta)^{\frac{n}{2}+1}}{t}\Big(\sum_{k\geq 1} 2^{(k+1)(\frac{n}{2}+2)}e^{-c\frac{2^{2k}l(\Delta)^2}{t^2}}\Big)\Vert \nabla_{||} f\Vert_{\infty}.
    \end{align*}

    With this estimate in hand, we return now to  
    \begin{align*}
    &\int_{T(\Delta)}\frac{|T_t((f-(f)_{3\Delta})(1-\eta))(x,t)|^2}{t}dxdt
    \\
    &\qquad\lesssim \int_0^{l(\Delta)}\frac{l(\Delta)^{n+2}}{t^3}\Big(\sum_{k\geq 1} 2^{(k+1)(\frac{n}{2}+2)}e^{-c\frac{2^{2k}l(\Delta)^2}{t^2}}\Big)^2\Vert \nabla_{||} f\Vert_{\infty}^2dt
    \\
    &\qquad \lesssim l(\Delta)^{n+2}\Vert \nabla_{||} f\Vert_{\infty}^2\int_0^{l(\Delta)}\Big(\sum_{k\geq 1} \frac{1}{t^{3/2}}2^{(k+1)(\frac{n}{2}+2)}e^{-c\frac{2^{2k}l(\Delta)^2}{t^2}}\Big)^2dt.
    \end{align*}

    We can study the function \(t\mapsto \frac{1}{t^{3/2}}e^{-c\frac{a^2}{t^2}}\) for \(t\in (0,l(\Delta))\) and observe that it attains its maximum on this interval if \(t=l(\Delta)\), whence
    \[\sum_{k\geq 1} \frac{1}{t^{3/2}}2^{(k+1)(\frac{n}{2}+2)}e^{-c\frac{2^{2k}l(\Delta)^2}{t^2}}\lesssim \frac{1}{l(\Delta)^{3/2}}\sum_{k\geq 1} 2^{(k+1)(\frac{n}{2}+2)}e^{-c2^{2k}}.\]
    We consider the sum as Riemann sums of the corresponding integral
    \[\int_0^\infty x^{\frac{n}{2}+1}e^{-cx^2}dx,\]
    which converges. Hence we obtain in total
    
    \begin{align*}
    &\int_{T(\Delta)}\frac{|T_t((f-(f)_{3\Delta})(1-\eta))(x,t)|^2}{t}dxdt
    \\
    &\qquad \lesssim l(\Delta)^{n+2}\Vert \nabla_{||} f\Vert_{\infty}^2\int_0^{l(\Delta)}\frac{1}{l(\Delta)^3}dt=\sigma(\Delta) \Vert \nabla_{||} f\Vert_{\infty}^2.
    \end{align*}
    This finishes the proof of \eqref{eq:CarlesonFctBound}.
\end{proof}

\subsection{\texorpdfstring{$L^1$}{L¹} Area function bounds on Hardy-Sobolev atoms}\label{appendix:SubsectionHardyBound}

\begin{defin}
    Let \(1<\beta<\infty\). We call a function \(a:\mathbb{R}^n\to\mathbb{R}\) a \textit{homogeneous Hardy-Sobolev \(\beta\)-atom} associated to a boundary ball \(\Delta\subset\mathbb{R}^n\) if
    \begin{enumerate}[(i)]
        \item \(\mathrm{supp}(a)\subset \Delta\); and
        \item \(\Vert \nabla a\Vert_{L^\beta(\mathbb{R}^n)}\leq |\Delta|^{-\frac{1}{\beta'}} \).
    \end{enumerate}
    If \(f\) can be written as
    \begin{align}f=\sum_{j=1}^\infty\lambda_ja_j\label{eq:AtomicHardyDecomposition}\end{align}
    for \(\beta\)-atoms \(a_j\) and coefficients \(\lambda_j\in\mathbb{R}\) with \(\sum_{j=1}^\infty|\lambda_j|<\infty\), we say that \(f\in \dot{HS}^{1,\beta}_{atom}\), where \(\Vert f\Vert_{\dot{HS}^{1,\beta}_{atom}}:=\inf\sum_{j=1}^\infty|\lambda_j|\) with an infimum that is take over all choices of decompositions \eqref{eq:AtomicHardyDecomposition}. 
\end{defin}

\begin{lemma}\label{lemma:HardyNormBounds}
    Let $T_t\in \{\partial_t\mathcal{P}_t, t\partial_t^2\mathcal{P}_t\}$ and let \(f\in \dot{HS}_{atom}^{1,2}(\mathbb{R}^n)\). Then
    \begin{align}
    \Vert\mathcal{A}(T_tf)\Vert_{L^1(\mathbb{R}^n)}\lesssim \Vert f\Vert_{\dot{HS}^{1,2}(\mathbb{R}^n)}.\label{eq:HardyAreaBound}
    \end{align}
\end{lemma}

\begin{proof}
    First, we note that it is enough to show
    \[\Vert\mathcal{A}(T_tf)\Vert_{L^1(\mathbb{R}^n)}\leq C\]
    for all homogeneous Hardy-Sobolev \(2\)-atoms \(f\) associated with \(\Delta\), whence we assume that \(f\) is such an atom going forward. 
    
    We begin with showing \eqref{eq:HardyAreaBound}. We split the integral into a local and a far away part 
    \begin{align*}
        \Vert\mathcal{A}(T_tf)\Vert_{L^1(\mathbb{R}^n)}=\Vert\mathcal{A}(T_tf)\Vert_{L^1(\mathbb{R}^n\setminus 5\Delta)} + \Vert\mathcal{A}(T_tf)\Vert_{L^1(5\Delta)}.
    \end{align*}
    For the local part, we have by Hölder's inequality, the $L^2$ area function bounds \refprop{prop:VerticalAreaFunctionL2Bounds}, and the properties of the Hardy space that
    \begin{align*}
        \Vert\mathcal{A}(T_tf)\Vert_{L^1(5\Delta)}\lesssim \Vert\mathcal{A}(T_tf)\Vert_{L^2(5\Delta)}|\Delta|^{1/2}\lesssim \Vert\nabla_{||} f\Vert_{L^2(\mathbb{R}^n)}|\Delta|^{1/2}\leq C.
    \end{align*}
    For the away part, we use the kernel bounds \refprop{Prop:KernelBounds} to observe
    \begin{align}
        \int_{\Delta_{t/2}(y)}|T_tf(x)|^2dx\lesssim \begin{cases} \frac{1}{t^{2n+2}}e^{-c\frac{\mathrm{dist}(\Delta_{t/2}(y),\Delta)^2}{t^2}}\Vert f\Vert_{L^1(\Delta)}^2 &\textrm{if }\mathrm{dist}(\Delta_{t/2}(y),\Delta)>0, \\ \frac{1}{t^{2n+2}}\Vert f\Vert_{L^1(\Delta)}^2& \mathrm{else}.\end{cases}\label{eq:CasesForL^1BoundofTermInProof}
    \end{align}
    We make use of \eqref{eq:CasesForL^1BoundofTermInProof} and Poincar\'{e}, since $f\in W^{1,p}_0(\Delta)$. We observe that if \(\mathrm{dist}(y,\Delta)\approx 2^jl(\Delta)\) then
    \begin{align*}
        &\int_{\Gamma(y)}\frac{|T_t(f)|^2}{t^{n+1}}dxdt=\int_0^\infty \frac{1}{t^{n+1}}\int_{\Delta_{t/2}(y)}|T_tf(x,t)|^2dx dt
        \\
        &\qquad=\int_0^{2^{j-1}l(\Delta)} \frac{e^{-2\frac{|2^{j}l(\Delta)-\frac{t}{2}|^2}{t^2}}}{t^{2n+3}}\Vert f\Vert_{L^1(\mathbb{R}^n)}^2dx dt + \int_{2^{j-1}l(\Delta)}^\infty \frac{1}{t^{2n+3}}\Vert f\Vert_{L^1(\mathbb{R}^n)}^2 dt
        \\
        &\qquad\lesssim\int_0^{2^{j-1}l(\Delta)} \frac{l(\Delta)^{n+2}e^{-\frac{2^{2(j-1)}l(\Delta)^2}{t^2}}}{t^{2n+3}}\Vert \nabla_{||}f\Vert_{L^2(\mathbb{R}^n)}^2dx dt
        \\
        &\qquad\qquad + \int_{2^{j-1}l(\Delta)}^\infty \frac{l(\Delta)^{n+2}}{t^{2n+3}}\Vert \nabla_{||}f\Vert_{L^2(\mathbb{R}^n)}^2 dt
        \\
        &\qquad \lesssim \frac{2^{-j(2n+2)}}{l(\Delta)^{n}}\Vert \nabla_{||}f\Vert_{L^2(\mathbb{R}^n)}^2\lesssim \frac{2^{-j(2n+2)}}{l(\Delta)^{2n}}.
    \end{align*}
    Summing over annuli, we obtain
    \begin{align*}
        \Vert\mathcal{A}(T_tf)\Vert_{L^1(\mathbb{R}^n\setminus 5\Delta)}&\lesssim \sum_{j\geq 3}\int_{2^j\Delta\setminus 2^{j-1}\Delta }\Big(\int_{\Gamma(y)}\frac{|T_tf|^2}{t^{n+1}}dxdt\Big)^{1/2}dy \lesssim \sum_{j\geq 3}\int_{2^j\Delta\setminus 2^{j-1}\Delta }\frac{2^{-j(n+1)}}{l(\Delta)^{n}}dy 
        \\
        &\lesssim \sum_{j\geq 3} 2^{-j}\lesssim C.
    \end{align*}
    This finishes the proof of \eqref{eq:HardyAreaBound}.
\end{proof}

\subsection{Tent spaces and real interpolation}\label{appendix:Subsection:Interpolation}

Let \(A_0,A_1\) be two normed vector spaces of functions \(a:\mathbb{R}^n\to\mathbb{R}\). For each \(a\in A_0 + A_1\) we define the \(K-\)functional of real interpolation by
\[K(a,t,A_0,A_1)=\inf_{a=a_0+a_1}\Vert a_0\Vert_{A_0} + t\Vert a_1\Vert_{A_1}.\]
For \(0<\theta<1,1\leq q\leq \infty\), we denote by \((A_0,A_1)_{\theta,q}\) the real interpolation space between \(A_0\) and \(A_1\) defined as
\[(A_0,A_1)_{\theta,q}=\Big\{a\in A_1+A_0:\Vert a\Vert_{\theta,q}=\Big(\int_0^\infty (t^{-\theta}K(a,t,A_0,A_1))^q\frac{\mathrm{d}t}{t}\Big)^{1/q}<\infty \Big\}.\]
According to Theorem 3.1.2 in \cite{bergh_interpolation_1976}, \(K\) can be seen as an exact interpolation functor, which means that if an operator \(T\) is bounded from \(A_0\to B_0\) and from \(A_1\to B_1\) for linear normed vector spaces of functions, then
\(T:(A_0,A_1)_{\theta,q}\to (B_0,B_1)_{\theta,q}\) is a bounded linear operator with \(\Vert T\Vert\leq C\Vert T\Vert^{1-\theta}_{A_0\to B_0}\Vert T\Vert^{\theta}_{A_1\to B_1}\).  

Now, we have the following real interpolation result.
\begin{prop}[Thm 0.4 in \cite{badr_abstract_2010}]
    For every \(\beta\in(1,\infty]\) and \(1<p<\infty\) the real interpolation space is
    \[(\dot{HS}^{1,\beta}_{atom}, \dot{W}^{1,\infty})_{1-\frac{1}{p},p}=\dot{W}^{1,p}.\]
\end{prop}

On the other hand, recall the definition of the area function \eqref{recall:areafunctiondef} and of the Carleson function \eqref{recall:carelsonfunctiondef} and define the tent spaces over \(\mathbb{R}^n\) with parameter \(1\leq p<\infty\) as
\begin{align*}
T^{p,2}(\Omega):=\{F\in L^2_{loc};\Vert F\Vert_{T^{p,2}}:=\Vert \mathcal{A}(F)\Vert_{L^p(\partial\Omega)}<\infty\}
\end{align*}
and
\begin{align*}
T^{\infty,2}(\Omega):=\{F\in L^2_{loc};\Vert F\Vert_{T^{\infty,2}}:=\Vert \mathcal{C}(F)\Vert_{L^\infty(\partial\Omega)}<\infty\}.
\end{align*}

By \cite{coifman_new_1985} we have
\begin{prop}[Theorem 4' in \cite{coifman_new_1985}]
    For every \(1<p<\infty\) the real interpolation space is
    \[(T^{1,2},T^{\infty,2})_{1-1/p,p}=T^{p,2}.\]
\end{prop}

Hence, real interpolation of the bounds from \reflemma{lemma:CarlesonFunctionBound} and \reflemma{lemma:HardyNormBounds} yields the following corollary:

\begin{cor}\label{cor:AreaFctBoundsInLp}
    For every \(1<p<\infty\), and for every \(f\in \dot{W}^{1,p}(\partial\Omega)\) we have
    \[\Vert \mathcal{A}(\partial_t\mathcal{P}_tf)\Vert_{L^p(\partial\Omega)}+\Vert \mathcal{A}(t\partial_t^2\mathcal{P}_tf)\Vert_{L^p(\partial\Omega)}\leq C\Vert \nabla_{||} f\Vert_{L^p(\partial\Omega)}.\]
\end{cor}

Furthermore, we also obtain
\begin{cor}\label{cor:AreaFctBoundsInLpNablaOperators}
    For every \(1<p<\infty\), and for every \(f\in \dot{W}^{1,p}(\partial\Omega)\) we have
    \[\Vert \mathcal{A}(t\nabla_{||}\partial_t\mathcal{P}_tf)\Vert_{L^p(\partial\Omega)}\leq C\Vert \nabla_{||} f\Vert_{L^p(\partial\Omega)}.\]
\end{cor}

This result might look surprising because the vertical area function of $t\nabla_{||}\partial_t\mathcal{P}_t$ admits this bound only in the range $1<p<2+\EPS$ (see \cite{hofmann_lp_2022}), which cannot be improved over $2+\EPS$, but the conical area function is bounded in the full range $1<p<\infty$.

\begin{proof}[Proof of \refcor{cor:AreaFctBoundsInLpNablaOperators}]
    First, we demonstrate the key argument only for $t\nabla_{||}\partial_t\mathcal{P}_tf$. For each $x\in\mathbb{R}^n$ we have
    \begin{align*}\mathcal{A}_{1/2}(t\nabla_{||}\partial_t\mathcal{P}_tf)^2(x)=\int_{\Gamma_{1/2}(x)}\frac{|t\nabla_{||}\partial_t\mathcal{P}_tf(y)|^2}{t^{n+1}}dydt
    \\
    \lesssim \sum_{W\in\mathcal{D};W\cap\Gamma_{1/2}(x)\neq\emptyset}\int_{W}\frac{|t\nabla_{||}\partial_t\mathcal{P}_tf(y)|^2}{t^{n+1}}dydt
    \end{align*}
    On each Whitney cube $W$, we can apply Caccioppoli inequality \reflemma{lemma:CaccTypeInequalityForP} and obtain
    \begin{align*}
    \sum_{W\in\mathcal{D};W\cap\Gamma_{1/2}(x)\neq\emptyset}\int_{W}\frac{|t\nabla_{||}\partial_t\mathcal{P}_tf(y)|^2}{t^{n+1}}dydt\lesssim \sum_{W\in\mathcal{D};W\cap\Gamma_{1/2}(x)\neq\emptyset}\int_{\frac{3}{2}W}\frac{|\partial_t\mathcal{P}_tf(y)|^2}{t^{n+1}}dydt
    \\
    \lesssim \mathcal{A}_5(\partial_t\mathcal{P}_tf)^2(x).
    \end{align*}
    Applying the comparability of $L^p$-norms for area function with different apertures \cite[Prop. 4.5]{milakis_harmonic_2013} gives
    \begin{align*}
    \Vert \mathcal{A}_\alpha(t\nabla_{||}\partial_t\mathcal{P}_tf)\Vert_{L^p(\partial\Omega)}&\lesssim \Vert \mathcal{A}_{1/2}(t\nabla_{||}\partial_t\mathcal{P}_tf)\Vert_{L^p(\partial\Omega)}\lesssim \Vert \mathcal{A}_5(\partial_t\mathcal{P}_tf)\Vert_{L^p(\partial\Omega)}
    \\
    &\lesssim \Vert \mathcal{A}_\alpha(\partial_t\mathcal{P}_tf)\Vert_{L^p(\partial\Omega)}\lesssim \Vert\nabla_{||}f\Vert_{L^p(\partial\Omega)},
    \end{align*}
    where we used \refcor{cor:AreaFctBoundsInLp} in the last inequality.

\end{proof}

Now, \refprop{prop:AreaFctBoundsInLp} follows from \refcor{cor:AreaFctBoundsInLp} and \refcor{cor:AreaFctBoundsInLpNablaOperators}.

\section{Technical lemmas}\label{appendix:technicalLemmas}

We present several results of rather technical nature that were used in above work, in particular in the proof of \refthm{thm:MainThm2}. The results of this section might not be surprising for the expert reader, however, we include them for completeness and rigor.

\subsection{Off-diagonal estimate for \texorpdfstring{$e^{-tL}\mathrm{div}$}{exp(-tL) div} (Proof of \refcor{cor:Off-diagonalEstimatesForPdiv})}\label{sec:AppendixOffDiago}

Let $h$ be a vector-valued function with components in $L^2(\mathbb{R}^n)$ with $\mathrm{supp}(f)\subset F$ and $E\subset\mathbb{R}^n$. Then there exist $g\in L^2(\mathbb{R}^n)$ with $\mathrm{supp}(g)\subset E$ and $\Vert g\Vert_{L^2(E)}=1$ such that
\begin{align*}
    \Vert \sqrt{t}e^{-tL}\mathrm{div}(h)\Vert_{L^2(E)}\lesssim\int_E \sqrt{t}e^{-tL}\mathrm{div}(h) \,g \,dx=\int_{\mathbb{R}^n}\sqrt{t}e^{-tL}\mathrm{div}(h) \,g \,dx.
\end{align*}
Using the adjoint operator, we obtain
\begin{align*}
    \Big|\int_{\mathbb{R}^n}\sqrt{t}e^{-tL}\mathrm{div}(h) \,g \,dx\Big|&=\Big|\int_{\mathbb{R}^n}h\cdot(\sqrt{t}\nabla e^{-tL^*}g) \,dx\Big|=\Big|\int_{F}h\cdot(\sqrt{t}\nabla e^{-tL^*}g) \,dx\Big|
    \\
    &\leq \Vert h\Vert_{L^2(F)}\Vert\sqrt{t}\nabla e^{-tL^*}g\Vert_{L^2(F)}.
\end{align*}
Using the off-diagonal estimates for $\nabla e^{-tL^*}$ in \refprop{prop:Off-diagonalEstimates}, yields \(\Vert\sqrt{t}\nabla e^{-tL^*}g\Vert_{L^2(F)}\lesssim e^{-\alpha\frac{d(E,F)^2}{t}}\), which finishes the proof of \refcor{cor:Off-diagonalEstimatesForPdiv}.

\subsection{Reverse Hölder inequality for the gradient of a Poisson solution (Proof of \reflemma{lemma:HigherIntegrabilityForGradient})}\label{sec:AppendixHigherIntProof}
Let $u\in W^{1,2}_{loc}(\Omega)$ be a solution to $Lu=H+\mathrm{div}(F)$, and let $B(x,2r)\subset\Omega$. Furthermore, let $\phi\in C^\infty_0(B(x,2r))$ be such that $\phi\equiv 1$ on $B(x,r)$, $0\leq \phi\leq 1$, and $|\nabla \phi|\lesssim \frac{1}{r}$. Then, for any constant $c$
\begin{align}
    \int_{B(x,r)}|\nabla u|^2 dy&\leq \int_{B(x,2r)}|\nabla u|^2\phi^2 dy\leq \lambda_0^{-1} \int_{B(x,2r)}A\nabla u\cdot\nabla (u-c) \phi^2 dy.\label{eq:appendixCacc}
\end{align}
We use that $A^a\nabla u\cdot\nabla u=0$ to continue with
\begin{align}
    &\nonumber\int_{B(x,2r)}A\nabla u\cdot\nabla (u-c) \phi^2 dy=\int_{B(x,2r)}(A-(A^a)_{B(x,2r)})\nabla u\cdot\nabla (u-c) \phi^2 dy
    \\
    &\label{eq:ProofRevHölderPDE}\qquad =\int_{B(x,2r)}A\nabla u\cdot\nabla ((u-c) \phi^2) dy - 2\int_{B(x,2r)}(A-(A^a)_{B(x,2r)})\nabla u \cdot \nabla\phi ((u-c)\phi) dy
    \\
    &\nonumber\qquad=\int_{B(x,2r)}H((u-c) \phi^2) dy+\int_{B(x,2r)}F\cdot\nabla ((u-c) \phi^2) dy
    \\
    &\nonumber\qquad\qquad- 2\int_{B(x,2r)}(A-(A^a)_{B(x,2r)})\nabla u \cdot \nabla\phi ((u-c)\phi) dy,
\end{align}
where we used the PDE for $u$ in the last equality. By Young's inequality,
\begin{align}
    &\int_{B(x,2r)}F\cdot\nabla ((u-c) \phi^2) dy=\int_{B(x,2r)}F\cdot(\nabla u \phi^2+2(u-c)\phi\nabla \phi) dy\nonumber
    \\
    &\leq  C_\kappa\int_{B(x,2r)}|F|^2 dy + \kappa\int_{B(x,2r)}|\nabla u|^2\phi^2 dy + C\int_{B(x,2r)}|u-c|^2|\nabla\phi|^2 dy.\label{eq:appendixCacc2}
\end{align}
where $\kappa$ is a small constant to be chosen. Similarly,
\begin{align}
    \int_{B(x,2r)}H ((u-c) \phi^2) dy\leq \int_{B(x,2r)}|\delta H|^2 dy + \frac{1}{r^2}\int_{B(x,2r)}|u-c|^2 dy,\label{eq:appendixCacc3}
\end{align}
Next, by H\"{o}lder
\begin{align*}
    &\int_{B(x,2r)}(A-(A^a)_{B(x,2r)})\nabla u \cdot \nabla\phi (u-c)\phi dy
    \lesssim \frac{1}{r}\int_{B(x,2r)}|A-(A^a)_{B(x,2r)}||\nabla u||\phi(u-c)| dy
    \\
    &\qquad \lesssim \frac{|B(x,2r)|}{r}\Big(\fint_{B(x,2r)}|A-(A^a)_{B(x,2r)}|^{\frac{s}{s-1}}dy\Big)^{\frac{s-1}{s}}\Big(\fint_{B(x,2r)}|\nabla u|^2\phi^2 dy\Big)^{1/2}\Big(\fint_{B(x,2r)}|u-c|^\frac{2s}{2-s} dy\Big)^{\frac{2-s}{2s}}
\end{align*}
where $s\in (1,\frac{n+1}{n})$. Since $A^a\in BMO(\Omega)$ and $A^s\in L^\infty(\Omega)$, it follows that $\Big(\fint_{B(x,2r)}|A-(A^a)_{B(x,2r)}|^{\frac{s}{s-1}}dy\Big)^{\frac{s-1}{s}}\leq C_s$. Then, by Young's inequality, we obtain
\[\int_{B(x,2r)}(A-(A^a)_{B(x,2r)})\nabla u \cdot \nabla\phi (u\phi) dy
    \lesssim  \frac{C_\kappa}{r^{2-2(n+1)/s'}}\Big(\int_{B(x,2r)}|u-c|^{\frac{2s}{2-s}} dy\Big)^{\frac{2-s}{s}} + \kappa\int_{B(x,2r)}|\nabla u|^2\phi^2 dy.\]
Combining this with \eqref{eq:appendixCacc}, \eqref{eq:appendixCacc2}, and \eqref{eq:appendixCacc3}, we obtain 

\begin{align*}
\int_{B(x,2r)}|\nabla u|^2\phi^2 dy\lesssim\int_{B(x,2r)}|\delta H|^2 dy + \int_{B(x,2r)}|F|^2 dy + \kappa\int_{B(x,2r)}|\nabla u|^2\phi^2 dy
\\
+  \frac{C_\kappa}{r^{2-2(n+1)/s'}}\Big(\int_{B(x,2r)}|u-c|^{\frac{2s}{2-s}} dy\Big)^{\frac{2-s}{s}}.
\end{align*}
Choosing $\kappa>0$ sufficiently small and hiding the second term on the left hand side yields a classical Caccioppoli-type inequality
\begin{align*}
\int_{B(x,r)}|\nabla u|^2 dy&\lesssim \int_{B(x,2r)}|\delta H|^2 dy+\int_{B(x,2r)}|F|^2 dy +  \frac{C_\kappa}{r^{2-2(n+1)/s'}}\Big(\int_{B(x,2r)}|u-c|^{\frac{2s}{2-s}} dy\Big)^{\frac{2-s}{s}},
\end{align*}
which completes the proof of \eqref{eq:RevBoundaryHölder2}. Now for \eqref{eq:RevBoundaryHölder1}, we can choose $c=(u)_{B(x,2r)}$ and apply Poincar\'{e}-Sobolev inequality, and obtain 
\begin{align*}
\int_{B(x,r)}|\nabla u|^2 dy&\lesssim \int_{B(x,2r)}|\delta H|^2dy+\int_{B(x,2r)}|F|^2 dy + r^{n+1}\Big(\fint_{B(x,2r)}|\nabla u|^{q} dy\Big)^{\frac{2}{q}}
\end{align*}
for some fixed $q<2$ solving $\frac{2s}{2-s}=q^*=\frac{q(n+1)}{n+1-q}$. The proof of \eqref{eq:RevBoundaryHölder1} is completed by \cite[Prop. 1.1 in Chapter V]{giaquinta_multiple_1983}.

Lastly, we observe that the same arguments apply, if \(x\in\partial\Omega\), we replace \(B(x,r)\) and \(B(x,2r)\) by \(T(\Delta(x,r))\) and \(T(\Delta(x,2r))\), respectively, and we set $\phi\in C^\infty(\Omega)$ to be such that $\phi\equiv 1$ on $T(\Delta(x,r))$, $\phi\equiv 0$ on $\Omega\setminus T(\Delta(x,2r))$, $0\leq \phi\leq 1$, and $|\nabla \phi|\lesssim \frac{1}{r}$. In that case, we just have to choose \(c=0\) so that the vanishing Dirichlet data of \(u\) allows to apply the PDE in \eqref{eq:ProofRevHölderPDE}, and notice that Poincar\'{e}-Sobolev inequality still applies to \(\Big(\int_{B(x,2r)}|u-c|^{\frac{2s}{2-s}} dy\Big)^{\frac{2-s}{s}}\). Other than that, the proof of \eqref{eq:RevBoundaryHölder3} is verbatim the same to the above presented proof for \eqref{eq:RevBoundaryHölder1}.

\subsection{Averaged area function}
Let us define the \textit{averaged area function} by
\[\sq{\mathcal{A}}_{\alpha,q}(F):=\int_{\Gamma_\alpha(x)}\Big(\fint_{B(y,\delta(y)/2)}|F|^q dz\Big)^{1/q}\frac{dy}{\delta(y)^{n}}.\]
This operator plays an underlying role because $\mathcal{T}_q$ is a discretized version of the averaged area function:
\begin{lemma}\label{lemma:AppendixAcompT}
    Let \(\alpha>1\). For any function $F$ and any point $x\in\partial\Omega$
    \[\sq{\mathcal{A}}_{\alpha/2,q}(F)(x)\leq \mathcal{T}_{\alpha,q}(F)(x)\leq\mathcal{T}_{3\alpha,q}(F)(x)\leq \sq{\mathcal{A}}_{5\alpha,q}(F)(x).\]
\end{lemma}

\begin{proof}
We note that $B(y,\delta(y),2)$ can be covered by a finite number of Whitney cubes in $\mathcal{D}$, that is uniform in $y$. Hence
\begin{align*}
\sq{\mathcal{A}}_{\alpha/2,q}(F)(x)&=\int_{\Gamma_\frac{\alpha}{2}(x)}\Big(\fint_{B(y,\delta(y)/2)}|F|^q dz\Big)^{1/q}\frac{dy}{\delta(y)^{n}}
\\
&\leq \sum_{W\cap\Gamma_\frac{\alpha}{2}(x)\neq\emptyset}\int_{W}\Big(\fint_{B(y,\delta(y)/2)}|F|^q dz\Big)^{1/q}\frac{dy}{\delta(y)^{n}}
\\
&\leq \sum_{W\cap\Gamma_\frac{\alpha}{2}(x)\neq\emptyset}l(W)\fint_W\Big(\sum_{\sq{W}\in\mathcal{D}:\sq{W}\cap B(y,\delta(y)/2)\neq \emptyset}\fint_{\sq{W}}|F|^q dz\Big)^{1/q}dy
\\
&\lesssim \sum_{\sq{W}\cap\Gamma_\alpha(x)\neq\emptyset}l(W)\Big(\fint_{\sq{W}}|F|^q dz\Big)^{1/q}=\mathcal{T}_{\alpha,q}(F)(x).
\end{align*}
Conversely, for every $W\in \mathcal{D}$ there exists a finite covering by balls with radius comparable to \(l(W)\) and hence a uniform constant independent of $W$ such that
\[\Big(\fint_W |F|^qdy\Big)^{1/q}\leq \fint_W\big(\fint_{B(y,\delta(y)/2)}|F|^qdz\big)^{1/q}dy.\]
 Thus,
\begin{align*}
\mathcal{T}_{3\alpha,q}(F)(x)&=\sum_{W\cap\Gamma_{3\alpha}(x)\neq\emptyset}l(W)\Big(\fint_{W}|F|^q dz\Big)^{1/q}
\\
&\leq \sum_{W\cap\Gamma_{3\alpha}(x)\neq\emptyset}l(W)^{-n}\int_W\big(\fint_{B(y,\delta(y)/2)}|F|^qdz\big)^{1/q}dy
\\
&\leq \int_{\Gamma_{5\alpha}(x)}\big(\fint_{B(y,\delta(y)/2)}|F|^qdz\big)^{1/q}\frac{dy}{\delta(y)^{n}}=\sq{\mathcal{A}}_{5\alpha,q}(F)(x).
\end{align*}
\end{proof}

Furthermore, the averaged area function has several properties of the classical area function, with proof that are verbatim the same:

\begin{lemma}\label{lemma:appendixAcompC}
    Let $\Omega$ be a CAD. For any function $F$, it holds that
    \begin{enumerate}[(a)]
\item\label{item:AleqCTilde} for all $0<p<\infty$: $\Vert\sq{\mathcal{A}}_q(F)\Vert_{L^p}\leq C_p\Vert\sq{\mathcal C}_q(F)\Vert_{L^p}$; and
\item\label{item:CleqATilde} for all $1<p<\infty$: $\Vert\sq{\mathcal C}_q(F)\Vert_{L^p}\leq C_p\Vert\sq{\mathcal{A}}_q(F)\Vert_{L^p}$.
\end{enumerate}
\end{lemma}

\begin{proof} 
This statement essentially follows from the ideas of \cite{coifman_new_1985} (or also \cite{milakis_harmonic_2013} on chord arc domains specifically). In \cite[Theorem 3]{coifman_new_1985} or \cite[Theorem 6.1]{milakis_harmonic_2013}, it is proven that for any function $F$:
\begin{align}\Vert\mathcal{A}(F)\Vert_{L^p}\leq C_p\Vert\mathcal C(F)\Vert_{L^p}\qquad \textrm{ for all }0<p<\infty\textrm{; and}\label{item:AleqC}
\\
\Vert\mathcal C(F)\Vert_{L^p}\leq C_p\Vert\mathcal{A}(F)\Vert_{L^p}\qquad \textrm{ for all }2<p<\infty.\label{item:CleqA}
\end{align}

The proof of \eqref{item:AleqC} employs a good $\lambda$ inequality, and if we replace $\mathcal{C}$ by $\sq{\mathcal{C}}_{q}$ and $\mathcal{A}$ by $\sq{\mathcal{A}}_{q}$,
then the proof of the good $\lambda$ inequality works completely analogous without any other modifications. This proves \eqref{item:AleqCTilde}

Next, we would like to prove \eqref{item:CleqATilde} in analogy to \eqref{item:CleqA}.
The proof idea stems from \cite[Theorem 3]{coifman_new_1985}, however, we would like to demonstrate here why the range of $p$ changes. To prove \eqref{item:CleqATilde}, we fix a boundary ball $\Delta$, and obtain by \cite{milakis_harmonic_2013}
\begin{align*}
    &\int_{T(\Delta)}\Big(\fint_{B(y,\delta(y)/2)}|F|^q dz\Big)^{1/q} dy=\int_{\Omega}\Big(\fint_{B(y,\delta(y)/2)}|F|^q dz\Big)^{1/q}\chi_{T(\Delta)}(y) dy
    \\
    &\qquad\lesssim \int_{\partial\Omega}\Big(\int_{\Gamma(x)}\Big(\fint_{B(y,\delta(y)/2)}|F|^q dz\Big)^{1/q} \frac{\chi_{T(\Delta)}}{\delta(y)^{n+1}}dy\Big) \, d\sigma(x)
    \\
    &\qquad\lesssim \int_{3\Delta}\Big(\int_{\Gamma(x)}\Big(\fint_{B(y,\delta(y)/2)}|F|^q dz\Big)^{1/q} \frac{dy}{\delta(y)^{n+1}}\Big) \, d\sigma(x).
\end{align*}
Dividing both sides by $\sigma(\Delta)$ and using that the boundary is $n$-Ahlfors regular, implies
\[\frac{1}{\sigma(\Delta)}\int_{T(\Delta)}\Big(\fint_{B(y,\delta(y)/2)}|F|^q dz\Big)^{1/q} dy\leq \frac{1}{\sigma(3\Delta)}\int_{3\Delta}\sq{\mathcal{A}}_{q}(F) \, d\sigma(x).\]
Taking a supremum implies the pointwise bound $\sq{\mathcal{C}}_q(F)\lesssim \mathcal{M}[\sq{\mathcal{A}}_q(F)]$, and hence boundedness of the maximal function implies that \eqref{item:CleqATilde} holds for all $1<p<\infty$.

\end{proof}

\begin{lemma}\label{lemma:appendixDiffApComp}
Let $\Omega$ be a CAD and let $1<p<\infty$. For all apertures $\alpha>1$, we have
\[\Vert\sq{\mathcal{A}}_{\alpha,q}(F)\Vert_{L^p}\approx_\alpha\Vert\sq{\mathcal{A}}_q(F)\Vert_{L^p} \]

If $\Omega=\mathbb{R}^{n+1}$, then the implicit constants are such that
\begin{align*}
    \Vert\sq{\mathcal{A}}_{\alpha,q}(F)\Vert_{L^p}&\leq C\alpha^{n} \Vert\sq{\mathcal{A}}_q(F)\Vert_{L^p}\qquad \textrm{and}
    \\
    \Vert\sq{\mathcal{A}}_q(F)\Vert_{L^p}&\leq C\alpha^{-n/p} \Vert\sq{\mathcal{A}}_{\alpha,q}(F)\Vert_{L^p}.
\end{align*}
\end{lemma}

 This proof goes back to \cite{coifman_new_1985}, and was proven in this form on chord arc domains in \cite[Prop. 4.5]{milakis_harmonic_2013} for the classical area function. The precise constants for the classical area function are proven in \cite{Auscher_change_2011}. However, the proof for the averaged area function works completely analogously.

\begin{cor}\label{cor:appendixDiffApCompT}
    Let $\Omega$ be a CAD. Let $1<p<\infty$. For all apertures $\alpha>0$, we have
\[\Vert\mathcal{T}_{\alpha,q}(F)\Vert_{L^p}\approx_\alpha\Vert\mathcal{T}_q(F)\Vert_{L^p}. \]

If $\Omega=\mathbb{R}^{n+1}$, then the implicit constants are such that
\begin{align*}
    \Vert\mathcal{T}_{\alpha,q}(F)\Vert_{L^p}&\leq C\alpha^n \Vert\mathcal{T}_q(F)\Vert_{L^p}\qquad \textrm{and}
    \\
    \Vert\mathcal{T}_q(F)\Vert_{L^p}&\leq C\alpha^{-n/p} \Vert\mathcal{T}_{\alpha,q}(F)\Vert_{L^p}.
\end{align*}
\end{cor}
    The combination of \reflemma{lemma:appendixDiffApComp} and \reflemma{lemma:AppendixAcompT} proves \refcor{cor:appendixDiffApCompT}.

    Finally, we conclude \eqref{eq:CarlAreaTcomparable}, i.e. that the averaged Carleson function $\sq{\mathcal{C}}$ and the discretized area function $\mathcal{T}$ are comparable in $L^p$.
\begin{proof}[Proof of \eqref{eq:CarlAreaTcomparable}]
    Let $\alpha>0$ be any aperture, $p>1$, and $q\geq 1$. Then by \reflemma{lemma:appendixAcompC} and \reflemma{lemma:appendixDiffApComp}
    \[\Vert \sq{\mathcal{C}}_{\alpha,q}(F)\Vert_{L^p}\approx \Vert \sq{\mathcal{A}}_{\alpha,q}(F)\Vert_{L^p}\approx \Vert \sq{\mathcal{A}}_{1/2,q}(F)\Vert_{L^p}\approx \Vert \sq{\mathcal{A}}_{5,q}(F)\Vert_{L^p}.\]
    By \refcor{cor:appendixDiffApCompT}, we have
    \[\Vert\mathcal{T}_{\alpha,q}(F)\Vert_{L^p}\approx\Vert\mathcal{T}_{3,q}(F)\Vert_{L^p},\]
    and \reflemma{lemma:AppendixAcompT} yields
    \[\Vert \sq{\mathcal{A}}_{1/2,q}(F)\Vert_{L^p}\lesssim \Vert \mathcal{T}_{3,q}(F)\Vert_{L^p}\lesssim \Vert \sq{\mathcal{A}}_{5,q}(F)\Vert_{L^p}.\]
    Combining all of the above implies \eqref{eq:CarlAreaTcomparable}.
\end{proof}

\subsection{Auxiliary lemmas for the proof of \refthm{thm:MainThm2}}\label{sec:AppendixAux}

The following lemma is used in the proof of \refthm{thm:MainThm2}. It is well established since \cite{coifman_new_1985} and \cite{milakis_harmonic_2013} that 
\[\int_\Omega \delta |fg| \, dx\lesssim \min\{\Vert\mathcal{A}(\delta f)\Vert_{L^p(\partial\Omega)}\Vert\mathcal{A}(\delta g)\Vert_{L^{p'}(\partial\Omega)}, \Vert\mathcal{A}(\delta f)\Vert_{L^p(\partial\Omega)}\Vert\mathcal{C}(\delta g)\Vert_{L^{p'}(\partial\Omega)}\}.\]
The following lemma is a discretized version thereof.

\begin{lemma}\label{lemma:dyadicToNondyadic}
    Let $\Omega\subset\mathbb{R}^{n+1}$, and let $\hat{W}:=\frac{3}{2}W$ be an enlargement of a Whitney cube $W\in\mathcal{D}(\Omega)$.
    Let \(f,g\in L^2_{loc}(\Omega)\), and \(1<p,q<\infty\), then
    \[\sum_{W\in\mathcal{D}(\Omega)}l(W)^{n+2}\Big(\fint_{\hat{W}}|f|^2\Big)^{1/2}\Big(\fint_{\hat{W}}|g|^2\Big)^{1/2}\lesssim \Vert\mathcal{A}(\delta f)\Vert_{L^p}\Vert\mathcal{A}(\delta g)\Vert_{L^{p'}},\]
    or
    \[\sum_{W\in\mathcal{D}(\Omega)}l(W)^{n+2}\Big(\fint_{\hat{W}}|f|^2\Big)^{1/2}\Big(\fint_{\hat{W}}|g|^q\Big)^{1/q}\lesssim \Vert\mathcal{A}(\delta f)\Vert_{L^p}\Vert \mathcal{T}_q(g)\Vert_{L^{p'}},\]
\end{lemma}

\begin{proof}
    We are going to prove both assertions analogously. First, set \[\Phi(x):=\sum_{W}\chi_W(x)l(W)\Big(\fint_{\hat{W}}|f|^2\Big)^{1/2}\Big(\fint_{\hat{W}}|g|^q\Big)^{1/q}\qquad\textrm{ for }x\in\Omega.\]
    Then
    \begin{align}
    \nonumber&\sum_{W}l(W)^{n+2}\Big(\fint_{\hat{W}}|f|^2\Big)^{1/2}\Big(\fint_{\hat{W}}|g|^q\Big)^{1/q}
    \\
    \nonumber&\lesssim \int_\Omega \Phi(x) dx\lesssim \int_{\partial\Omega}\Big(\int_{\Gamma(x)}\frac{\Phi(y)}{\delta(y)^{n}}dy\Big)d\sigma(x)
    \\
    \nonumber&\lesssim \int_{\partial\Omega}\Big(\sum_{W:W\cap\Gamma(x)\neq \emptyset}\int_{W}\frac{\Phi(y)}{\delta(y)^{n}}dy\Big)d\sigma(x)
    \\
    &\lesssim \int_{\partial\Omega}\Big(\sum_{W:W\cap\Gamma(x)\neq \emptyset}\int_{W}\frac{\Big(\fint_{\hat{W}}|f|^2\Big)^{1/2}\Big(\fint_{\hat{W}}|g|^q\Big)^{1/q}}{\delta(y)^{n-1}}dy\Big)d\sigma(x)
    \label{eq:EqInProofDyadicToNondyadic}
\end{align}
At this point we would like to apply Hölder's inequality, and we need to distinguish two cases.
If \(q=2\), then
\begin{align*}
    \eqref{eq:EqInProofDyadicToNondyadic}&\lesssim \int_{\partial\Omega}\sum_{W:W\cap\Gamma(x)\neq \emptyset}\Big(\int_{\hat{W}}\frac{|f|^2}{\delta^{n-1}}\Big)^{1/2}\Big(\int_{\hat{W}}\frac{|g|^2}{\delta^{n-1}}\Big)^{1/2}d\sigma(x)
    \\&\lesssim \int_{\partial\Omega}\Big(\int_{\Gamma_3(x)}\frac{|f|^2}{\delta^{n-1}}dy\Big)^{1/2}\Big(\int_{\Gamma_3(x)}\frac{|g|^2}{\delta^{n-1}}dy\Big)^{1/2}d\sigma(x)
    \lesssim \Vert\mathcal{A}_3(\delta f)\Vert_{L^p}\Vert\mathcal{A}_3(\delta g)\Vert_{L^{p'}}.
\end{align*}

If however \(q\neq 2\), then we additionally use that \(\ell^2\subset\ell^1\) and get
\begin{align*}
    \eqref{eq:EqInProofDyadicToNondyadic}
    &\lesssim \int_{\partial\Omega}\sum_{W:W\cap\Gamma(x)\neq \emptyset}\Big(\int_{\hat{W}}\frac{|f|^2}{\delta^{n-1}}\Big)^{1/2}\Big[l(W)\Big(\fint_{\hat{W}}|g|^q\Big)^{1/q}\Big]d\sigma(x)
    \\
    &\lesssim \int_{\partial\Omega}\Big(\int_{\Gamma_3(x)}\frac{|f|^2}{\delta^{n-1}}\Big)^{1/2}\Big(\sum_{W:W\cap\Gamma(x)\neq \emptyset}\Big[l(W)\Big(\fint_{\hat{W}}|g|^q\Big)^{1/q}\Big]^2\Big)^{1/2}d\sigma(x)
    \\
    &\lesssim \int_{\partial\Omega}\Big(\int_{\Gamma_3(x)}\frac{|f|^2}{\delta^{n-1}}\Big)^{1/2}\Big(\sum_{W:W\cap\Gamma(x)\neq \emptyset}l(W)\Big(\fint_{\hat{W}}|g|^q\Big)^{1/q}\Big)d\sigma(x)
    \\
    &\lesssim \Vert\mathcal{A}_3(\delta f)\Vert_{L^p}\Vert \mathcal{T}_{3,q}(g)\Vert_{L^{p'}}.
\end{align*}

Lastly, let us justify that the apertures on the right hand sides can be reduced to 1. For that we use that $\Vert\mathcal{A}_\alpha(F)\Vert_{L^p}\approx\Vert\mathcal{A}(F)\Vert_{L^p}$ which goes back to \cite{coifman_new_1985}, and was proven in this form on chord arc domains in \cite[Prop. 4.5]{milakis_harmonic_2013}.
In addition to that, for $\mathcal{T}_{3,q}$, we use \reflemma{lemma:AppendixAcompT} and \reflemma{lemma:appendixDiffApComp} to conclude
\[\Vert \mathcal{T}_{3,q}(g)\Vert_{L^{p'}}\leq \Vert \sq{\mathcal{A}}_{5,q}(g)\Vert_{L^{p'}}\lesssim \Vert \sq{\mathcal{A}}_q(g)\Vert_{L^{p'}}\leq \Vert \mathcal{T}_{q}(g)\Vert_{L^{p'}}.\]
    
\end{proof}

\subsection{\texorpdfstring
    {Bound of $\sq{\mathcal{N}}_2(\nabla_{||}\mathcal{P}_tf)$ and $\sq{\mathcal{N}}_2(\mathcal{P}_tf)$ (Proof of \eqref{eq:nontangnetialofNablaPBoundedbyNablaf}/\eqref{eq:nontangnetialofPBoundedbyf})}
    {Nontangential maximal bounds for the semigroup and its tangential gradient (Proof of (\ref{eq:nontangnetialofNablaPBoundedbyNablaf})/(\ref{eq:nontangnetialofPBoundedbyf}))}}\label{sec:AppendixNontangnetialBound}

To start with, let us recall a useful result $(3.35)$ in \cite[Proof of Prop. 11]{hofmann_dirichlet_2022}. There the authors show that for every \((y,t)\in \Gamma_{\alpha}(x)\)
\begin{align}
    |2tL_{||}e^{-t^2L_{||}}f(y)|\leq C_\alpha \mathcal{M}(\nabla_{||} f)(x).\label{eq:semigroupinConeBoundedNablaf}
\end{align}
Their proof uses the Gaussian kernel estimates as the only property of \(tL_{||}e^{-t^2L_{||}}\), and hence an analogous proof also gives
\begin{align}
    |e^{-t^2L_{||}}(f-(f)_{\Delta(x,2(1+\alpha)t)})(y)|\leq C_\alpha \mathcal{M}(\nabla_{||} f)(x).\label{eq:semigroup-constantinConeBoundedNablaf}
\end{align}
Moreover, a simplified version of the same proof also gives 
\begin{align}
    |e^{-t^2L_{||}}f(y)|\leq C_\alpha \mathcal{M}(|f|)(x).\label{eq:semigroupinConeBoundedf}
\end{align}
Indeed, the kernel bounds \refprop{Prop:KernelBounds} yield
\begin{align*}
\bigl|e^{-t^2L_{||}}f(y)\bigr|
&\lesssim t^{-n}\sum_{k=1}^{\infty}
\exp\left(-\frac{(2^k-1)^2\alpha^2}{c}\right)
\int_{\Delta(x,2^{k+1}(1+\alpha) t)\setminus \Delta(x,2^{k}(1+\alpha) t)}|f|\,d\sigma + t^{-n}\int_{\Delta(x,2(1+\alpha) t)}|f|\,d\sigma
\\
&\lesssim(1+\alpha)^{n+1}
\left[
1+\sum_{k=1}^{\infty}
2^{k(n+1)}
\exp\left(-\frac{(2^k-1)^2\alpha^2}{c}\right)
\right]\mathcal{M}(| f|)(x)\lesssim_\alpha \mathcal{M}(|f|)(x).
\end{align*}

Now, we turn to the proofs of \eqref{eq:nontangnetialofPBoundedbyf} and \eqref{eq:nontangnetialofNablaPBoundedbyNablaf}. For \eqref{eq:nontangnetialofPBoundedbyf}, the observation \eqref{eq:semigroupinConeBoundedf} immediately implies with \(L^p\) boundedness of the maximal function 
\[\Vert\mathcal{N}(\mathcal{P}_tf)\Vert_{L^p}\lesssim \Vert\mathcal{M}(|f|)\Vert_{L^p}\lesssim \Vert f\Vert_{L^p}.\]

Similarly, \eqref{eq:nontangnetialofNablaPBoundedbyNablaf} can be reduced to proving the claim
\begin{align}
    \sq{\mathcal{N}}_2(\nabla_{||}\mathcal{P}_tf)(x)\lesssim \mathcal{M}(|\nabla_{||}f|)(x)\qquad\textrm{ for all }x\in\partial\Omega,\label{eq:PointwiseN(Nablaf)lesssimMaximalFct}
\end{align}

To prove this claim, we fix a point $x\in \partial\Omega$, and let $(y,t)\in\Gamma_\alpha(x)$ be such that $\sq{\mathcal{N}}_2(\nabla_{||}\mathcal{P}_tf)(x)\lesssim \Big(\fint_{B(y,t,t/2)}|\nabla_{||}\mathcal{P}_s f |^2dzds\Big)^{1/2}$. We can introduce a smooth cut-off function $\psi\in C_0^\infty(B(y,t,3t/4))$ with $\psi\equiv 1$ on $B(y,t,t/2)$ and $|\nabla\psi|\lesssim\frac{1}{t}$. Then, we obtain by integration by parts for any constant $c\in\mathbb{R}$
\begin{align*}
    &\fint_{B(y,t,t/2)}|\nabla_{||}\mathcal{P}_s f |^2dzds\leq \fint_{B(y,t,3t/4)}|\nabla_{||}\mathcal{P}_s (f-c) |^2\psi^2dzds
    \\
    &\lesssim\fint_{B(y,t,3t/4)}A \nabla_{||}\mathcal{P}_s (f-c) \cdot\nabla_{||}\mathcal{P}_s (f-c) \psi^2 dzds
    \\
    &\lesssim-\fint_{B(y,t,3t/4)}L_{||}\mathcal{P}_s (f-c) \mathcal{P}_s (f-c) \psi^2 dzds- \fint_{B(y,t,3t/4)}A\nabla_{||}\mathcal{P}_s (f-c) \cdot\nabla_{||}\psi^2 \mathcal{P}_s (f-c)  dzds.
    \\
    &=:I_1+I_2.
\end{align*}
For $I_2$ we use \refprop{prop:ContantAntiSymmMatrix} and H\"{o}lder to obtain
\begin{align*}
    I_2&=\frac{1}{2}\fint_{B(y,t,3t/4)}A\nabla_{||}(\mathcal{P}_s (f-c))^2 \cdot\nabla_{||}\psi^2  dzds
    \\
    &=\frac{1}{2}\fint_{B(y,t,3t/4)}(A-(A^a)_{B(y,t,3t/4)})\nabla_{||}(\mathcal{P}_s (f-c))^2 \cdot\nabla_{||}\psi^2  dzds
    \\
    &\lesssim \Big(\fint_{B(y,t,3t/4)}|A-(A^a)_{B(y,t,3t/4)}|^4dzds\Big)^{1/4}\Big(\fint_{B(y,t,3t/4)}|\nabla_{||}\mathcal{P}_s (f-c)|^2\psi^2dzds\Big)^{1/2}
    \\
    &\qquad\cdot\Big(\fint_{B(y,t,3t/4)}|\mathcal{P}_s (f-c)|^4|\nabla_{||}\psi|^4dzds\Big)^{1/4}
\end{align*}
Since $A\in BMO(\Omega)$, the first factor is bounded by a constant only depending on $\Lambda_0$. We can apply Young's inequality for a small $\sigma>0$ and get
\begin{align*}
I_2\lesssim \sigma \fint_{B(y,t,3t/4)}|\nabla_{||}\mathcal{P}_s(f-c)|^2\psi^2 dzds + C_\sigma \Big(\fint_{B(y,t,3t/4)}\frac{|\mathcal{P}_s(f-c)|^4}{s^4}dzds\Big)^{1/2}.
\end{align*}
For a sufficiently small choice of $\sigma$, the first term can be hidden on the left hand side. For $I_1$ we obtain
\begin{align*}
    I_1\lesssim\fint_{B(y,t,3t/4)}\partial_s\mathcal{P}_s(f-c)\frac{\mathcal{P}_s(f-c)}{s}\psi^2 dzds,
\end{align*}
and hence together with the bound for $I_2$
\begin{align*}
    I_1+I_2\lesssim \fint_{B(y,t,3t/4)}\partial_s\mathcal{P}_s(f-c)\frac{\mathcal{P}_s(f-c)}{s}\psi^2 dzds+ \Big(\fint_{B(y,t,3t/4)}\frac{|\mathcal{P}_s(f-c)|^4}{s^4}dzds\Big)^{1/2}
\end{align*}
By \eqref{eq:semigroup-constantinConeBoundedNablaf} and \eqref{eq:semigroupinConeBoundedNablaf}, we get for all \((z,s)\in B(y,t,3t/4)\)
\begin{align}|\partial_s\mathcal{P}_s(f-c)(z)|+ \frac{|\mathcal{P}_s(f-c)(z)|}{s}\lesssim \mathcal{M}[|\nabla_{||}f|](x),\label{eq:PointWiseBoundOnP/2}\end{align}
which finishes the proof of the claim \eqref{eq:PointwiseN(Nablaf)lesssimMaximalFct}.

\bibliographystyle{alpha}
\raggedbottom
\bibliography{references} 

@phdthesis{li_lp_2019,
	title = {Lp {Dirichlet} problem for second order elliptic operators having a {BMO} anti-symmetric part},
	url = {https://repository.library.brown.edu/studio/item/bdr:918947/},
	language = {en},
	urldate = {2021-07-01},
	school = {Brown University},
	author = {Li, Linhan},
	year = {2019},
}

@article{christ_tb_1990,
	title = {A {T}(b) theorem with remarks on analytic capacity and the {Cauchy} integral},
	volume = {60-61},
	issn = {0010-1354},
	url = {https://eudml.org/doc/265825},
	language = {eng},
	number = {2},
	urldate = {2021-02-03},
	journal = {Colloquium Mathematicae},
	author = {Christ, Michael},
	year = {1990},
	pages = {601--628},
}

@article{Lewis_1988,
	title = {Uniformly Fat Sets},
	volume = {308},
	issn = {0001-8708},
	url = {https://doi.org/10.2307/2000957},
	language = {en},
	journal = {Transactions AMS},
    number = {1},
    month = {July}, 
	author = {Lewis, J.},
	year = {1988},
	pages = {177-196},
}

@article{dahlberg_area_1984,
	title = {Area integral estimates for elliptic differential operators with non-smooth coefficients},
	volume = {22},
	issn = {0004-2080, 1871-2487},
	url = {https://projecteuclid.org/journals/arkiv-for-matematik/volume-22/issue-1-2/Area-integral-estimates-for-elliptic-differential-operators-with-non-smooth/10.1007/BF02384374.full},
	doi = {10.1007/BF02384374},
	number = {1-2},
	urldate = {2021-06-24},
	journal = {Arkiv för Matematik},
	author = {Dahlberg, Björn E. J. and Jerison, David S. and Kenig, Carlos E.},
	month = dec,
	year = {1984},
	note = {Publisher: Institut Mittag-Leffler},
	pages = {97--108},
}

@article{Ataei2025,
  title={The Kato square root problem for parabolic operators with an anti-symmetric part in BMO},
  author={Alireza Ataei and Kaj Nystr{\"o}m},
  journal={Nonlinear Analysis},
  year={2025},
  url={https://api.semanticscholar.org/CorpusID:279639037}
}

@article{DNU2024, 
title = {Perturbation Theory for Second Order Elliptic Operators with BMO Antisymmetric Part},
url = {https://doi.org/10.1007/s10013-023-00653-z},
journal = {Vietnam J. Math., 52},
author= {Dindoš, M. and Nyström, E. and Ulmer, M.},
year = {2024},
pages = {519–566}
}

@article{hofmann_dirichlet_2022,
	title = {The {Dirichlet} problem for elliptic operators having a {BMO} anti-symmetric part},
	volume = {382},
	issn = {1432-1807},
	url = {https://doi.org/10.1007/s00208-021-02219-1},
	doi = {10.1007/s00208-021-02219-1},
	language = {en},
	number = {1},
	urldate = {2022-07-26},
	journal = {Mathematische Annalen},
	author = {Hofmann, Steve and Li, Linhan and Mayboroda, Svitlana and Pipher, Jill},
	month = feb,
	year = {2022},
	pages = {103--168},
}

@article{hofmann_regularity_2015,
	title = {The regularity problem for second order elliptic operators with complex-valued bounded measurable coefficients},
	volume = {361},
	issn = {1432-1807},
	url = {https://doi.org/10.1007/s00208-014-1087-6},
	doi = {10.1007/s00208-014-1087-6},
	language = {en},
	number = {3},
	urldate = {2024-04-03},
	journal = {Mathematische Annalen},
	author = {Hofmann, Steve and Kenig, Carlos and Mayboroda, Svitlana and Pipher, Jill},
	month = apr,
	year = {2015},
	pages = {863--907},
}

@article{hofmann_lp_2022,
	title = {\${L}{\textasciicircum}p\$ theory for the square roots and square functions of elliptic operators having a {BMO} anti-symmetric part},
	volume = {301},
	issn = {0025-5874, 1432-1823},
	url = {http://arxiv.org/abs/1908.01030},
	doi = {10.1007/s00209-021-02938-w},
	number = {1},
	urldate = {2022-12-12},
	journal = {Mathematische Zeitschrift},
	author = {Hofmann, Steve and Li, Linhan and Mayboroda, Svitlana and Pipher, Jill},
	month = may,
	year = {2022},
	note = {arXiv:1908.01030 [math]},
	pages = {935--976},
}

@article{jerison_kenig_1981,
 ISSN = {0003486X},
 URL = {http://www.jstor.org/stable/2006988},
 author = {David S. Jerison and Carlos E. Kenig},
 journal = {Annals of Mathematics},
 number = {2},
 pages = {367--382},
 publisher = {Annals of Mathematics},
 title = {The Dirichlet Problem in Non-Smooth Domains},
 urldate = {2023-07-17},
 volume = {113},
 year = {1981}
}

@article{Auscher_Kato,
 ISSN = {0003486X},
 URL = {http://www.jstor.org/stable/3597201},
 author = {Pascal Auscher and Steve Hofmann and Michael Lacey and Alan McIntosh and Ph. Tchamitchian},
 journal = {Annals of Mathematics},
 number = {2},
 pages = {633--654},
 publisher = {Annals of Mathematics},
 title = {The Solution of the Kato Square Root Problem for Second Order Elliptic Operators on  Rn},
 urldate = {2023-10-27},
 volume = {156},
 year = {2002}
}

@article{kenig_new_2000,
	title = {A {New} {Approach} to {Absolute} {Continuity} of {Elliptic} {Measure}, with {Applications} to {Non}-symmetric {Equations}},
	volume = {153},
	issn = {00018708},
	url = {https://linkinghub.elsevier.com/retrieve/pii/S0001870899918996},
	doi = {10.1006/aima.1999.1899},
	language = {en},
	number = {2},
	urldate = {2023-11-08},
	journal = {Advances in Mathematics},
	author = {Kenig, C. and Koch, H. and Pipher, J. and Toro, T.},
	month = aug,
	year = {2000},
	pages = {231--298},
}

@article{hofmann_square_2015,
	title = {Square {Function}/{Non}-{Tangential} {Maximal} {Function} {Estimates} and the {Dirichlet} {Problem} for {Non}-{Symmetric} {Elliptic} {Operators}},
	volume = {28},
	issn = {0894-0347},
	url = {https://www.jstor.org/stable/43302884},
	number = {2},
	urldate = {2024-11-05},
	journal = {Journal of the American Mathematical Society},
	author = {Hofmann, Steve and Kenig, Carlos and Mayboroda, Svitlana and Pipher, Jill},
	year = {2015},
	note = {Publisher: American Mathematical Society},
	pages = {483--529},
}

@article{badr_abstract_2010,
	title = {Abstract {Hardy}–{Sobolev} spaces and interpolation},
	volume = {259},
	issn = {0022-1236},
	url = {https://www.sciencedirect.com/science/article/pii/S0022123610000807},
	doi = {10.1016/j.jfa.2010.02.012},
	number = {5},
	urldate = {2024-06-26},
	journal = {Journal of Functional Analysis},
	author = {Badr, Nadine and Bernicot, Frédéric},
	month = sep,
	year = {2010},
	pages = {1169--1208},
}

@book{bergh_interpolation_1976,
	address = {Berlin, Heidelberg},
	series = {Grundlehren der mathematischen {Wissenschaften}},
	title = {Interpolation {Spaces}: {An} {Introduction}},
	volume = {223},
	copyright = {http://www.springer.com/tdm},
	isbn = {978-3-642-66453-3 978-3-642-66451-9},
	shorttitle = {Interpolation {Spaces}},
	url = {http://link.springer.com/10.1007/978-3-642-66451-9},
	language = {en},
	urldate = {2024-06-26},
	publisher = {Springer Berlin Heidelberg},
	author = {Bergh, Jöran and Löfström, Jörgen},
	year = {1976},
	doi = {10.1007/978-3-642-66451-9},
}

@article{coifman_new_1985,
	title = {Some new function spaces and their applications to {Harmonic} analysis},
	volume = {62},
	issn = {0022-1236},
	url = {https://www.sciencedirect.com/science/article/pii/0022123685900072},
	doi = {10.1016/0022-1236(85)90007-2},
	number = {2},
	urldate = {2024-08-13},
	journal = {Journal of Functional Analysis},
	author = {Coifman, R. R and Meyer, Y and Stein, E. M},
	month = jun,
	year = {1985},
	pages = {304--335},
}

@article{kenig_neumann_1995,
	title = {The {Neumann} problem for elliptic equations with nonsmooth coefficients: {Part} {II}},
	volume = {81},
	issn = {0012-7094},
	shorttitle = {The {Neumann} problem for elliptic equations with nonsmooth coefficients},
	url = {https://projecteuclid.org/journals/duke-mathematical-journal/volume-81/issue-1/The-Neumann-problem-for-elliptic-equations-with-nonsmooth-coefficients/10.1215/S0012-7094-95-08112-5.full},
	doi = {10.1215/S0012-7094-95-08112-5},
	language = {en},
	number = {1},
	urldate = {2023-01-09},
	journal = {Duke Mathematical Journal},
	author = {Kenig, Carlos E. and Pipher, Jill},
	month = jan,
	year = {1995},
}

@misc{mourgoglou_regularity_2023,
	title = {The regularity problem for the {Laplace} equation in rough domains},
	url = {http://arxiv.org/abs/2110.02205},
	doi = {10.48550/arXiv.2110.02205},
	urldate = {2025-06-23},
	publisher = {arXiv},
	author = {Mourgoglou, Mihalis and Tolsa, Xavier},
	month = aug,
	year = {2023},
	note = {arXiv:2110.02205 [math]},
}

@article{shen_relationship_2007,
	title = {A relationship between the dirichlet and regularity problems for elliptic equations},
	volume = {14},
	copyright = {Copyright 2018 Elsevier B.V., All rights reserved.},
	issn = {1073-2780},
	language = {eng},
	number = {2-3},
	journal = {Mathematical research letters},
	author = {Shen, Zhongwei},
	year = {2007},
	note = {Place: SOMERVILLE
Publisher: Int Press Boston, Inc},
	pages = {205--213},
}

@article{milakis_harmonic_2013,
	title = {Harmonic {Analysis} on {Chord} {Arc} {Domains}},
	volume = {23},
	issn = {1050-6926, 1559-002X},
	url = {http://link.springer.com/10.1007/s12220-012-9322-4},
	doi = {10.1007/s12220-012-9322-4},
	language = {en},
	number = {4},
	urldate = {2021-07-22},
	journal = {Journal of Geometric Analysis},
	author = {Milakis, Emmanouil and Pipher, Jill and Toro, Tatiana},
	month = oct,
	year = {2013},
	pages = {2091--2157},
}

@article{hofmann_greens_function_2007,
   title={The Green function estimates for strongly elliptic systems of second order},
   volume={124},
   ISSN={1432-1785},
   url={http://dx.doi.org/10.1007/s00229-007-0107-1},
   DOI={10.1007/s00229-007-0107-1},
   number={2},
   journal={Manuscripta mathematica},
   publisher={Springer Science and Business Media LLC},
   author={Hofmann, Steve and Kim, Seick},
   year={2007},
   pages={139–172} }

@misc{mourgoglou_greens_function_2023,
      title={Regularity theory and Green's function for elliptic equations with lower order terms in unbounded domains}, 
      author={Mihalis Mourgoglou},
      year={2023},
      eprint={1904.04722},
      archivePrefix={arXiv},
      primaryClass={math.AP},
      url={https://arxiv.org/abs/1904.04722}, 
}

@article{mourgoglou_solvability_2025,
	title = {Solvability of the {Poisson}–{Dirichlet} problem with interior data in ${L}^{p'}$-{Carleson} spaces and its applications to the ${L}^p$-regularity problem},
	issn = {1435-9855},
	url = {https://ems.press/journals/jems/articles/14298929},
	doi = {10.4171/jems/1660},
	language = {en},
	urldate = {2025-08-15},
	journal = {Journal of the European Mathematical Society},
	author = {Mourgoglou, Mihalis and Poggi, Bruno and Tolsa, Xavier},
	month = jul,
	year = {2025},
}

@article{escauriaza_kato_2018,
	title = {Kato square root problem with unbounded leading coefficients},
	volume = {146},
	issn = {0002-9939, 1088-6826},
	url = {https://www.ams.org/proc/2018-146-12/S0002-9939-2018-14224-0/},
	doi = {10.1090/proc/14224},
	language = {en},
	number = {12},
	urldate = {2024-04-02},
	journal = {Proceedings of the American Mathematical Society},
	author = {Escauriaza, Luis and Hofmann, Steve},
	month = dec,
	year = {2018},
	pages = {5295--5310},
}

@article{Auscher_change_2011,
title = {Change of angle in tent spaces},
journal = {Comptes Rendus Mathematique},
volume = {349},
number = {5},
pages = {297-301},
year = {2011},
issn = {1631-073X},
doi = {https://doi.org/10.1016/j.crma.2011.01.023},
url = {https://www.sciencedirect.com/science/article/pii/S1631073X11000409},
author = {Pascal Auscher}
}

@article{meyer_gradient_1963,
     author = {Meyers, Norman G.},
     title = {An $L^p$-estimate for the gradient of solutions of second order elliptic divergence equations},
     journal = {Annali della Scuola Normale Superiore di Pisa - Scienze Fisiche e Matematiche},
     pages = {189--206},
     year = {1963},
     publisher = {Scuola normale superiore},
     volume = {Ser. 3, 17},
     number = {3},
     mrnumber = {159110},
     zbl = {0127.31904},
     language = {en},
     url = {https://www.numdam.org/item/ASNSP_1963_3_17_3_189_0/}
}

@book{giaquinta_multiple_1983,
	title = {Multiple {Integrals} in the {Calculus} of {Variations} and {Nonlinear} {Elliptic} {Systems}. ({AM}-105)},
	isbn = {978-0-691-08331-5},
	url = {http://www.jstor.org/stable/j.ctt1b9s07q},
	urldate = {2026-06-17},
	publisher = {Princeton University Press},
	author = {GIAQUINTA, MARIANO},
	year = {1983},
}

@article{jones_extension_1980,
	title = {Extension {Theorems} for {BMO}},
	volume = {29},
	issn = {0022-2518},
	url = {https://www.jstor.org/stable/24892564},
	number = {1},
	urldate = {2026-07-01},
	journal = {Indiana University Mathematics Journal},
	publisher = {Indiana University Mathematics Department},
	author = {Jones, Peter W.},
	year = {1980},
	pages = {41--66},
}

@article{fefferman_hpspaces_1972,
	title = {Hpspaces of several variables},
	volume = {129},
	issn = {1871-2509},
	url = {https://doi.org/10.1007/BF02392215},
	doi = {10.1007/BF02392215},
	language = {en},
	number = {1},
	urldate = {2025-03-04},
	journal = {Acta Mathematica},
	author = {Fefferman, C. and Stein, E. M.},
	month = dec,
	year = {1972},
	pages = {137--193},
}

@article{carleson_interpolations_1962,
	title = {Interpolations by {Bounded} {Analytic} {Functions} and the {Corona} {Problem}},
	volume = {76},
	issn = {0003-486X},
	url = {https://www.jstor.org/stable/1970375},
	doi = {10.2307/1970375},
	number = {3},
	urldate = {2026-07-01},
	journal = {Annals of Mathematics},
	publisher = {[Annals of Mathematics, Trustees of Princeton University on Behalf of the Annals of Mathematics, Mathematics Department, Princeton University]},
	author = {Carleson, Lennart},
	year = {1962},
	pages = {547--559},
}

@article{hytonen_carleson_2013,
	title = {On the {Carleson} duality},
	volume = {51},
	issn = {0004-2080, 1871-2487},
	url = {https://projecteuclid.org/journals/arkiv-for-matematik/volume-51/issue-2/On-the-Carleson-duality/10.1007/s11512-012-0167-7.full},
	doi = {10.1007/s11512-012-0167-7},
	number = {2},
	urldate = {2025-09-17},
	journal = {Arkiv för Matematik},
	publisher = {Institut Mittag-Leffler},
	author = {Hytönen, Tuomas and Rosén, Andreas},
	month = oct,
	year = {2013},
	pages = {293--313},
}

@article{li_pipher_bmo,
      author = {Linhan Li and Jill Pipher},
        title = {Boundary behavior of solutions of elliptic operators in divergence form with a BMO anti-symmetric part},
        journal = {Communications in Partial Differential Equations},
        volume = {44},
        number = {2},
        pages = {156--204},
        year = {2019},
        publisher = {Taylor \& Francis},
        doi = {10.1080/03605302.2018.1542437},
}

@book{heinonen_nonlinear_1993,
  author    = {Juha Heinonen and Tero Kilpel{\"a}inen and Olli Martio},
  title     = {Nonlinear Potential Theory of Degenerate Elliptic Equations},
  series    = {Oxford Mathematical Monographs},
  publisher = {The Clarendon Press, Oxford University Press},
  address   = {New York},
  year      = {1993},
  note       = {Oxford Science Publications}
}

@article{bourgain_hausdorff_1987,
  author  = {Jean Bourgain},
  title   = {On the Hausdorff Dimension of Harmonic Measure in Higher Dimension},
  journal = {Inventiones Mathematicae},
  volume  = {87},
  number  = {3},
  pages   = {477--483},
  year    = {1987},
  doi     = {10.1007/BF01389238}
}

@article{seregin_divergence-free_2012,
	title = {On divergence-free drifts},
	volume = {252},
	issn = {0022-0396},
	url = {https://www.sciencedirect.com/science/article/pii/S0022039611003597},
	doi = {https://doi.org/10.1016/j.jde.2011.08.039},
	number = {1},
	journal = {Journal of Differential Equations},
	author = {Seregin, Gregory and Silvestre, Luis and {\v S}verák, Vladimír and Zlato{\v s}, Andrej},
	year = {2012},
	pages = {505--540},
}

@article{escauriaza_hofmann_kato_2018,
  author  = {Escauriaza, Luis and Hofmann, Steve},
  title   = {{Kato} square root problem with unbounded leading coefficients},
  journal = {Proceedings of the American Mathematical Society},
  volume  = {146},
  number  = {12},
  pages   = {5295--5310},
  year    = {2018},
  doi     = {10.1090/proc/14224}
}

@article{dong_kim_fundamental_2018,
  author  = {Dong, Hongjie and Kim, Seick},
  title   = {Fundamental solutions for second-order parabolic systems
             with drift terms},
  journal = {Proceedings of the American Mathematical Society},
  volume  = {146},
  number  = {7},
  pages   = {3019--3029},
  year    = {2018},
  doi     = {10.1090/proc/14004}
}

@article{qian_xi_parabolic_2019,
  author  = {Qian, Zhongmin and Xi, Guangyu},
  title   = {Parabolic equations with singular divergence-free drift
             vector fields},
  journal = {Journal of the London Mathematical Society},
  series  = {2},
  volume  = {100},
  number  = {1},
  pages   = {17--40},
  year    = {2019},
  doi     = {10.1112/jlms.12202}
}

@article{ulmer_solvability_2026,
  author  = {Ulmer, Martin},
  title   = {Solvability of the {D}irichlet problem for a new class of elliptic operators},
  journal = {Journal of Differential Equations},
  volume  = {485},
  pages   = {114740},
  year    = {2026},
	doi = {10.1016/j.jde.2026.114740},
}

@misc{ulmer_equivalence_2025,
	title = {Equivalence between solvability of the {Dirichlet} and {Regularity} problem under a mixed {$L^1-L^\infty$} condition on {$\partial_tA$}},
	url = {http://arxiv.org/abs/2509.10328},
	doi = {10.48550/arXiv.2509.10328},
	urldate = {2025-11-17},
	publisher = {arXiv},
	author = {Ulmer, Martin},
	month = oct,
	year = {2025},
	note = {arXiv:2509.10328 [math]},
}
\end{document}